\documentclass[12pt]{article}

\usepackage[margin=1.2in]{geometry}

\usepackage{amssymb}

\usepackage{lineno}

\usepackage{amsmath,amsthm,amsfonts}
\usepackage{tikz}
\usepackage[colorlinks=false,linkcolor=blue,citecolor=green]{hyperref}
\usepackage{enumerate}
\usepackage{booktabs}

\newtheorem{theorem}{Theorem}

\newtheorem{lemma}{Lemma}
\newtheorem{proposition}{Proposition}

\newtheorem{remark}{Remark}

\numberwithin{equation}{section}
\numberwithin{figure}{section}
\numberwithin{table}{section}

\renewcommand{\leq}{\leqslant}
\renewcommand{\geq}{\geqslant}

\newcommand{\rme}{\mathrm{e}}
\newcommand{\bfi}{\mathbf{i}}

\begin{document}

\title{Increasing stability of the first order linearized inverse Schr\"{o}dinger potential problem with integer power type nonlinearities%
\thanks{S. ZOU and S. LU is supported by Key-Area Research and Development Program of Guangdong Province (No.2021B0101190003), NSFC (No.11925104), Science and Technology Commission of Shanghai Municipality (21JC1400500). B. XU is supported by NSFC (No.12171301 and No.11801351).}
}


%
%
%

\author{%
Sen ZOU%
\footnote{ School of Mathematical Sciences, Fudan University. %
Email: \texttt{szou18@fudan.edu.cn}}, %
Shuai LU%
\footnote{ School of Mathematical Sciences, Fudan University. %
Email: \texttt{slu@fudan.edu.cn}}, %
Boxi XU%
\footnote{ Corresponding author. School of Mathematics, Shanghai University of Finance and Economics. %
Email: \texttt{xu.boxi@mail.sufe.edu.cn}}
}

\maketitle


\begin{abstract}
We investigate the increasing stability of the inverse Schr\"{o}dinger potential problem with integer power type nonlinearities at a large wavenumber. By considering the first order linearized system with respect to the unknown potential function, a combination formula of the first order linearization is proposed, which provides a Lipschitz type stability for the recovery of the Fourier coefficients of the unknown potential function in low frequency mode. These stability results highlight the advantage of nonlinearity in solving this inverse potential problem by explicitly quantifying the dependence to the wavenumber and the nonlinearities index. A reconstruction algorithm for general power type nonlinearities is also provided. Several numerical examples illuminate the efficiency of our proposed algorithm.
\end{abstract}






%
%
%







\section{Introduction}
\label{sec:introduction}

In this paper, we consider the inverse Schr\"{o}dinger potential problem with power type nonlinearities.
Specifically, the nonlinear Schr\"{o}dinger equation is given by
\begin{equation}\label{eqn:original_0}
\Delta u + k^{2} u - c(x) u^{m} =0 \quad \text{in\ } \Omega \subset \mathbb{R}^{n},
\end{equation}
where the wavenumber $k > 0$, the integer $m \geq 2$, $\Omega \subset \mathbb{R}^{n}$ is an open bounded domain with $C^{\infty}$ smooth boundary $\partial\Omega$ and the dimensionality $n \geq 2$.
The inverse potential problem considered in this paper is to recover the unknown potential function $c := c(x)$ from the linearized Dirichlet-to-Neumann (DtN) map with respect to the potential function, which will be defined later.
Assuming that the potential function $c(x)$ is sufficiently small and using the linearized DtN map as observation, we show that the recovery of the potential function $c(x)$ satisfies increasing stability, which also leads to a stable reconstruction algorithm depending explicitly on the wavenumber $k$ and the index $m$.

In view of related inverse problems of \eqref{eqn:original_0}, more is know when the forward model \eqref{eqn:original_0} is linear, for instance when $m = 1$. When the wavenumber $k = 0$, the proposed inverse Schr\"{o}dinger potential problem is closely related to the electrical impedance tomography problem where a logarithmic stability was proved in \cite{A1988} and further validated to be optimal in \cite{M2001}.
When $k > 0$, an increasing stability estimate was firstly obtained in \cite{I2011}, showing that the inversion resolution shall improve when the wavenumber $k$ increases.
For increasing stability of other inverse problems with different linear models, we refer to \cite{SI2007, NUW2013, I2015, CIL2016, LY2017, KUW2021} and references therein. Numerically, compared to the logarithmic stability when the wavenumber $k = 0$, a direct benefit of the increasing stability is to construct more stable inversion algorithms via the linearized DtN map. For instance, \cite{ILX2020} proposed a Fourier-based reconstruction algorithm for the first order linearized inverse Schr\"{o}dinger potential problem for \eqref{eqn:original_0} with $m = 1$ and $k \geq 1$.
It has been shown there that the Fourier coefficients of the unknown potential function can be stably recovered from the linearized DtN map in a range of Fourier modes, and this range increases with the growing wavenumber.
This first order linearization approach is further applied to other inverse potential problems within different settings and we refer to \cite{ILX2022, ZLX2022}.
Here we mention a rencent study \cite{CFO2022} about the rigidity of DtN map, which gives the local uniqueness result by showing the DtN map is locally convex with respect to the potential function near $0$.

When the forward problem is nonlinear, initiated by the study for inverse problems of nonlinear parabolic equations in \cite{I1993}, the linearization approach with respect to boundary data has been widely applied in solving various inverse problems of nonlinear models, for example in \cite{IS1994, SU1997, S2010, SZ2012, LL2019} where the first order linearization of the DtN map is an (indirect) DtN map of the linearized equation.
Recently, the higher order linearization approach was introduced in \cite{KLU2018} for parameter identification of the wave equation on Lorentzian manifolds and further extended in \cite{FO2020, LLLS2021} to solve inverse problems of \eqref{eqn:original_0} when $k = 0$ and $m \in \mathbb{N}_{+}$.
Precisely, by considering the $m$-th order derivative with respect to the Dirichlet boundary data with different small (scale) variables for the nonlinear forward problem, one reduces the original nonlinear equation to a linear one.
Therefore uniqueness and stability results can be derived to solve certain inverse problems arising from the nonlinear equations \eqref{eqn:original_0}.
This higher order linearization method has been further developed to solve various types of inverse problems for nonlinear equations \cite{LLPT2020, FLL2021, LLLS2021b, LLST2022, CFKKU2021}.
Here we mention a recent numerical study \cite{LLPT2022} on the inverse problem for the nonlinear wave equations using higher order linearization of DtN map.

Both above two linearization methods, i.e., the first order linearization (with respect to small potential function) and higher order linearization (with respect to small Dirichlet boundary data), have been well discussed and compared in \cite{LSX2022} for solving the inverse Schr\"{o}dinger potential problem of \eqref{eqn:original_0} with an integer nonlinear index $m \geq 2$.
Using the higher order linearization method by deriving the $m$-th derivative of the DtN map, the increasing stability for recovery of the potential function $c$ is proved in \cite{LSX2022} for an arbitrary finite integer $m \geq 2$. Meanwhile, only for $m = 2$, the first order linearization method with respect to the potential function for the nonlinear problem \eqref{eqn:original_0} can guarantee the increasing stability as shown in \cite[Section 3]{LSX2022}. Whether such an approach can provide the increasing stability for an arbitrary finite integer $m \geq 3$ is not answered yet.
In this paper, by integrating the principle of inclusion-exclusion (PIE) in combinatorics, we show that for $m \geq 3$ the increasing stability holds true for the first order linearization method with respect to the potential function, which is consistent with the results by the higher order linearization method in \cite{LSX2022} and links the intrinsic connection between two linearization approaches.

The paper is organized in the following manner.
In Section \ref{sec:problem} we introduce the (first order) linearized DtN map for \eqref{eqn:original_0} and prove an Alessandrini-PIE type identity which allows us to recover the unknown potential function $c$ by a combination of the first order linearized DtN map with different boundary data.
In Section \ref{sec:main_results} we present the increasing stability results and prove those theorems.
The reconstruction algorithm and numerical results are presented in Section \ref{sec:algorithm} and \ref{sec:numerical} respectively.


\section{The first order linearized inverse Schr\"{o}dinger potential problem and the Alessandrini-PIE type identity}
\label{sec:problem}

We first recall the inverse Schr\"{o}dinger potential problem with integer power type nonlinearities as follows.
In this paper, let an integer $m \geq 2$, we consider the nonlinear Schr\"{o}dinger equation with Dirichlet boundary condition, that is
\begin{equation}\label{eqn:problem}
(I) ~
\left\{\begin{aligned}
\Delta u + k^{2} u - c(x) u^{m} &= 0 & & \text{in\ } \Omega \subset \mathbb{R}^{n}, \\
u &= f & & \text{on\ } \partial\Omega.
\end{aligned}\right.
\end{equation}
Here $k > 1$ is the wavenumber, $\Omega \subset \mathbb{R}^{n}$ is an open bounded domain with $C^{\infty}$ smooth boundary $\partial\Omega$ and the dimension $n \geq 2$. For convenience, denote $c := c(x)$ as the potential function with $\operatorname{supp}(c) \subset \Omega$. Our aim is to recover the unknown potential function $c$ from many boundary measurements, or more precisely the linearized Dirichlet-to-Neumann (DtN) map.

To begin with, we first state the well-posedness of the original problem \eqref{eqn:problem} with small Dirichlet boundary data.

\begin{proposition}\label{prp:wellposed}
Let $m, n \geq 2$ be integers, and $\Omega \subset \mathbb{R}^{n}$ be an open bounded domain with $C^{\infty}$ boundary $\partial\Omega$. Assume that the potential function $c \in L^{\infty}(\Omega)$ with $\operatorname{supp}(c) \subset \Omega$, and the wavenumber $k > 1$ such that $k^{2}$ is not a Dirichlet eigenvalue of $-\Delta$ in $\Omega$.

Denote $p := \frac{n}{2} + 1$, there exist two small constants $\eta > 0$ and $\eta_{c} > 0$, if any Dirichlet boundary data $f$ and potential function $c$ satisfy that
\begin{equation*}
f \in B_{\eta}
:= \left\{ f \in W^{2-\frac{1}{p},p}(\partial\Omega) \,:\, \|f\|_{W^{2-\frac{1}{p},p}(\partial\Omega)} < \eta \right\},
\end{equation*}
and
\begin{equation*}
\|c\|_{L^{\infty}(\Omega)} < \eta_{c},
\end{equation*}
then the boundary value problem \eqref{eqn:problem} has a unique solution $u \in W^{2,p}(\Omega)$, which satisfies the a-priori estimate
\begin{equation}\label{eqn:nonlinear_apriori}
\|u\|_{W^{2,p}(\Omega)}
\leq C(m,\Omega,k) \left( \|c\|_{L^{\infty}(\Omega)} + \|f\|_{W^{2-\frac{1}{p},p}(\partial\Omega)} \right).
\end{equation}

Furthermore, the nonlinear Dirichlet-to-Neumann map for \eqref{eqn:problem} is a $C^{\infty}$ map defined by
\begin{equation}\label{eqn:dtn}
\Lambda_{c} : B_{\eta} \to W^{1-\frac{1}{p},p}(\partial\Omega),
\quad f \mapsto \partial_{\nu} u |_{\partial\Omega}
\end{equation}
such that $\Lambda_{c}(f) := \partial_{\nu} u |_{\partial\Omega}$.
\end{proposition}

\begin{proof}
The proof is derived by the implicit function theorem for Banach spaces \cite[Theorem 10.6 and Remark 10.5]{RR2004}, and is similar to that in \cite[Proof of Theorem 2.1]{N2022}. See \ref{sec:wellposed} for the details.
\end{proof}

\begin{remark}
When $k = 0$, the well-posedness results for \eqref{eqn:problem} with small boundary data are proved in \cite{LLLS2021} and \cite{N2022} where $c$ is $C^{\infty}(\overline{\Omega})$ and $L^{\frac{n}{2}+\epsilon}(\Omega)$ correspondingly.
\end{remark}


\subsection{The linearization method and linearized DtN map $\Lambda'_{c}$}\label{sec:linearized}

Next, we state the first order linearized inverse Schr\"{o}dinger potential problem by introducing a linearization method to the DtN map $\Lambda_{c}$ in \eqref{eqn:dtn}.

Let $\gamma$ be a small constant. Given a small potential function $\gamma c \in L^{\infty}(\Omega)$ such that $\| \gamma c \|_{L^{\infty}(\Omega)} < \eta_{c}$ in Proposition \ref{prp:wellposed}, by computing the partial derivative of DtN map \eqref{eqn:dtn} with respect to $\gamma$, we can derive the linearized DtN map as below,
\begin{equation}\label{eqn:dtn_linear_derive}
\Lambda'_{c}(f)
:= \partial_{\gamma} \big( \Lambda_{\gamma c} \big) \big|_{\gamma = 0} (f)
= \partial_{\nu} u^{(1)} \big|_{\partial\Omega},
\end{equation}
where $f \in B_{\eta} \subset W^{2-\frac{1}{p},p}(\partial\Omega)$ and the solution $u^{(1)}$ is derived by solving the following linearized system,
\begin{equation}\label{eqn:problem_0}
(I_{0}) ~
\left\{\begin{aligned}
\Delta u^{(0)} + k^{2} u^{(0)} &= 0 & & \text{in\ } \Omega, \\
u^{(0)} &= f & & \text{on\ } \partial\Omega,
\end{aligned}\right.
\end{equation}
for $u^{(0)}$ of the unperturbed problem, 
and
\begin{equation}\label{eqn:problem_1}
(I_{1}) ~
\left\{\begin{aligned}
\Delta u^{(1)} + k^{2} u^{(1)} &= c(x) ( u^{(0)} )^{m} & & \text{in\ } \Omega, \\
u^{(1)} &= 0 & & \text{on\ } \partial\Omega,
\end{aligned}\right.
\end{equation}
for $u^{(1)}$ of the linearized problem. 
Note that, here we define $\Lambda'_{c}$ for small $f \in B_{\eta}$ to guarantee the well-posedness of the original problem \eqref{eqn:problem}, but later we will show one can extend the domain of $\Lambda'_{c}$ to $W^{2-\frac{1}{p},p}(\partial\Omega)$ continuously.

According to Proposition \ref{prp:wellposed}, by the assumption that $k^{2}$ is not a Dirichlet eigenvalue of $-\Delta$ in $\Omega$, there exists a unique solution $u^{(0)} \in W^{2,p}(\Omega)$ if $f \in B_{\eta} \subset W^{2-\frac{1}{p},p}(\partial\Omega)$ and $c \equiv 0$ for the unperturbed problem \eqref{eqn:problem_0}.
Furthermore, by the Sobolev embedding theorem \cite[Theorem 4.12]{AF2003}, we have $u^{(0)} \in C^{0,s}(\overline{\Omega})$ for some $s \in (0,1)$. Therefore, $c ( u^{(0)} )^{m} \in L^{\infty}(\Omega) \subset L^{p}(\Omega)$, which admits a unique solution $u^{(1)} \in W^{2,p}(\Omega)$ for the linearized problem \eqref{eqn:problem_1}.
Consequently, following by the linearized system \eqref{eqn:problem_0}--\eqref{eqn:problem_1}, we have
\begin{equation}\label{eqn:dtn_linear}
\Lambda'_{c} : B_{\eta}
 \to W^{1-\frac{1}{p},p}(\partial\Omega),
\quad f \mapsto \partial_{\nu} u^{(1)} \big|_{\partial\Omega}.
\end{equation}
Here we note that, according to the definition \eqref{eqn:dtn_linear}, the linearized DtN map $\Lambda'_{c}$ is given by solving the system \eqref{eqn:problem_0}--\eqref{eqn:problem_1}, in which the potential function $c$ is not necessary to be small. Indeed, the linearized DtN map $\Lambda'_{c}$ can be defined for any $c \in L^{\infty}(\Omega)$.

However, the DtN map $\Lambda_{c}$ in \eqref{eqn:dtn} is well-defined only if the potential function $c$ satisfies the smallness assumption $\|c\|_{L^{\infty}(\Omega)} < \eta_{c}$ in Proposition \ref{prp:wellposed}. In this situation, the following proposition shows the linearized DtN map $\Lambda'_{c}$ is indeed the first order linearization of DtN map $\Lambda_{c}$ with respect to small potential function $c$.

\begin{proposition}\label{prp:linearized}
Under the assumptions and notations in Proposition \ref{prp:wellposed}.
For any $f \in B_{\eta}$, it holds that
\begin{equation}\label{eqn:linearized_err}
\| \Lambda_{c}(f) - \Lambda_{0}(f) - \Lambda'_{c}(f) \|_{W^{1-\frac{1}{p},p}(\partial\Omega)}
\sim O(\|c\|_{L^{\infty}(\Omega)}^{2}).
\end{equation}
\end{proposition}

\begin{proof}
Following the definitions of the DtN map $\Lambda_{c}$ in \eqref{eqn:dtn} and the linearized DtN map $\Lambda'_{c}$ in \eqref{eqn:dtn_linear}, we have
\begin{equation*}
\Lambda_{c}(f) - \Lambda_{0}(f) - \Lambda'_{c}(f)
= \big( \partial_{\nu} u - \partial_{\nu} u^{(0)} - \partial_{\nu} u^{(1)} \big) \big|_{\partial\Omega}, \quad \text{for\ } f \in B_{\eta},
\end{equation*}
where $u$, $u^{(0)}$ and $u^{(1)}$ solve the original problem $(I)$ in \eqref{eqn:problem}, the unperturbed problem $(I_{0})$ in \eqref{eqn:problem_0} and the linearized problem $(I_{1})$ in \eqref{eqn:problem_1}, respectively.

Notice that when $p = \frac{n}{2} + 1$, by Proposition \ref{prp:wellposed} and Sobolev embedding theorem \cite[Theorem 4.12]{AF2003}, the solution $u$ of the original problem \eqref{eqn:problem} satisfies $u \in W^{2,p}(\Omega) \subset C^{0,s}(\overline{\Omega}) \subset L^{\infty}(\Omega)$ for some $s \in (0,1)$.
More precisely it holds that
\begin{equation*}
\|u\|_{L^{\infty}(\Omega)}
\leq C(\Omega) \|u\|_{C^{0,s}(\overline{\Omega})}
\leq C(\Omega) \|u\|_{W^{2,p}(\Omega)}.
\end{equation*}
Therefore $c u^{m} \in L^{p}(\Omega)$ with the estimate
\begin{equation*}
\|c u^{m}\|_{L^{p}(\Omega)}
\leq C(\Omega) \|c\|_{L^{\infty}(\Omega)} \|u\|_{L^{\infty}(\Omega)}^{m}
\leq C(\Omega) \|c\|_{L^{\infty}(\Omega)} \|u\|_{W^{2,p}(\Omega)}^{m}.
\end{equation*}
By the assumption that $k^{2}$ is not a Dirichlet eigenvalue of $-\Delta$ in $\Omega$, then $u^{(0)} \in W^{2,p}(\Omega)$ given $f \in B_{\eta} \subset W^{2-\frac{1}{p},p}(\partial\Omega)$.
In the same way, we can derive
\begin{equation*}
\|u^{(0)}\|_{L^{\infty}(\Omega)}
\leq C(\Omega) \|u^{(0)}\|_{W^{2,p}(\Omega)},
\end{equation*}
and
\begin{equation*}
\|c (u^{(0)})^{m}\|_{L^{p}(\Omega)}
\leq C(\Omega) \|c\|_{L^{\infty}(\Omega)} \|u^{(0)}\|_{W^{2,p}(\Omega)}^{m}.
\end{equation*}

Subtracting the original problem \eqref{eqn:problem} with the unperturbed problem \eqref{eqn:problem_0}, we have
\begin{equation*}
\left\{\begin{aligned}
( \Delta + k^{2} ) ( u - u^{(0)} ) &= c u^{m} & & \text{in\ } \Omega, \\
u - u^{(0)} &= 0 & & \text{on\ } \partial\Omega,
\end{aligned}\right.
\end{equation*}
and the estimate
\begin{equation*}
\begin{aligned}
\| u - u^{(0)} \|_{W^{2,p}(\Omega)}
&\leq C(\Omega,k) \|c u^{m}\|_{L^{p}(\Omega)}
\leq C(\Omega,k) \|c\|_{L^{\infty}(\Omega)} \|u\|_{W^{2,p}(\Omega)}^{m}.
\end{aligned}
\end{equation*}
In a similar manner, by the linearized problem \eqref{eqn:problem_1}, we have
\begin{equation*}
\begin{aligned}
\| u - u^{(0)} - u^{(1)} \|_{W^{2,p}(\Omega)}
&\leq C(\Omega,k) \|c\|_{L^{\infty}(\Omega)} \| u^{m} - ( u^{(0)} )^{m} \|_{L^{\infty}(\Omega)}.
\end{aligned}
\end{equation*}

By the a-priori estimate \eqref{eqn:nonlinear_apriori} in Proposition \ref{prp:wellposed} and the binomial theorem, the last term can be bounded by
\begin{equation*}
\begin{aligned}
&\| u^{m} - ( u^{(0)} )^{m} \|_{L^{\infty}(\Omega)} \\
&\qquad \leq C(\Omega) \| u - u^{(0)} \|_{W^{2,p}(\Omega)} \sum_{j=1}^{m} \|u\|_{L^{\infty}(\Omega)}^{m-j} \|u^{(0)}\|_{L^{\infty}(\Omega)}^{j-1} \\
&\qquad \leq C(\Omega,k) \|c\|_{L^{\infty}(\Omega)} \sum_{j=1}^{m} \|u\|_{W^{2,p}(\Omega)}^{2m-j} \|u^{(0)}\|_{W^{2,p}(\Omega)}^{j-1} \\
&\qquad \leq C(m,\Omega,k) \|c\|_{L^{\infty}(\Omega)} \sum_{j=1}^{m} \left( \|c\|_{L^{\infty}(\Omega)} + \|f\|_{W^{2-\frac{1}{p},p}(\partial\Omega)} \right)^{2m-j} \|f\|_{W^{2-\frac{1}{p},p}(\partial\Omega)}^{j-1} \\
&\qquad \leq C(m,\Omega,k) \sum_{j=1}^{m} \sum_{\ell=0}^{2m-j} \binom{2m-j}{\ell} \|f\|_{W^{2-\frac{1}{p},p}(\partial\Omega)}^{2m-\ell-1} \|c\|_{L^{\infty}(\Omega)}^{\ell+1} \\
\end{aligned}
\end{equation*}
Hence we have
\begin{equation*}
\| u - u^{(0)} - u^{(1)} \|_{W^{2,p}(\Omega)}
\leq C(m,\Omega,k) \sum_{j=1}^{m} \sum_{\ell=0}^{2m-j} \binom{2m-j}{\ell} \|f\|_{W^{2-\frac{1}{p},p}(\partial\Omega)}^{2m-\ell-1} \|c\|_{L^{\infty}(\Omega)}^{\ell+2}.
\end{equation*}

By trace theorem \cite[Theorem 7.39]{AF2003}, $\partial_{\nu} : W^{2,p}(\Omega) \to W^{1-\frac{1}{p},p}(\partial\Omega)$ is a bounded linear operator, which leads to
\begin{equation*}
\big\| \big( \partial_{\nu} u - \partial_{\nu} u^{(0)} - \partial_{\nu} u^{(1)} \big) \big|_{\partial\Omega} \big\|_{W^{1-\frac{1}{p},p}(\partial\Omega)}
\lesssim C(m,\Omega,k) \|f\|_{W^{2-\frac{1}{p},p}(\partial\Omega)}^{2m-1} \|c\|_{L^{\infty}(\Omega)}^{2}.
\end{equation*}
It concludes the proof of Proposition \ref{prp:linearized}.
\end{proof}

Noticing that the linearized DtN map $\Lambda'_{c}$ is a nonlinear operator with respect to the Dirichlet boundary data $f$.
By checking the linearized system \eqref{eqn:problem_0}--\eqref{eqn:problem_1}, it can be verified that $\Lambda'_{c}$ is $m$-homogeneous, i.e.
\begin{equation}\label{eqn:m_homogeneous}
\Lambda'_{c}(\beta f)
= \beta^{m} \Lambda'_{c}(f), \quad \text{for $f \in B_{\eta}$ and $\beta \geq 0$},
\end{equation}
which means the domain of $\Lambda'_{c}$ can be continuously extended to the whole space $W^{2-\frac{1}{p},p}(\partial\Omega)$ by scaling. Therefore, in the rest of this paper, we use the same notation $\Lambda'_{c} : W^{2-\frac{1}{p},p}(\partial\Omega) \to W^{1-\frac{1}{p},p}(\partial\Omega)$ to denote the extended operator.
Furthermore. we can define the ``nonlinear operator norm'' of $\Lambda'_{c}$ by
\begin{equation}\label{eqn:dtn_linear_norm}
\| \Lambda'_{c} \|_{\mathcal{N}}
:= \sup_{f \neq 0} \frac{ \| \Lambda'_{c}(f) \|_{W^{1-\frac{1}{p},p}(\partial\Omega)} }{ \|f\|_{W^{2-\frac{1}{p},p}(\partial\Omega)}^{m} }.
\end{equation}

Now we propose the linearized inverse Schr\"{o}dinger potential problem: \textbf{recovering the unknown potential function $c$ from the knowledge of the (first order) linearized DtN map $\Lambda'_{c}$ defined in \eqref{eqn:dtn_linear}.}


\subsection{The Alessandrini-PIE type identity}\label{sec:identity}

As a general framework for inverse medium problems, in order to build the relation between the boundary observation and unknown information about the interior medium, it's common to derive the so-called Alessandrini type identity \cite[Lemma 1]{A1988}. While for the problem studied in our current work, the main difficulty is the nonlinearity from the DtN map, which we are going to demonstrate in the following subsection.

Multiplying the equation \eqref{eqn:problem} from both sides with test function $\varphi$ solving $\Delta \varphi + k^{2} \varphi = 0$ in $\Omega$, we obtain the Alessandrini type identity for the DtN map $\Lambda_{c}$, such that
\begin{equation}\label{eqn:identity_dtn}
\begin{aligned}
\int_{\Omega} c(x) u^{m} \varphi \,\mathrm{d}x
&= \int_{\partial\Omega} \big( ( \partial_{\nu} u ) \varphi - u ( \partial_{\nu} \varphi ) \big) \,\mathrm{d}S_{x} \\
&= \int_{\partial\Omega} \big( \Lambda_{c}(f) \, \varphi - f \, \Lambda_{0}(\varphi) \big) \,\mathrm{d}S_{x}.
\end{aligned}
\end{equation}
Meanwhile, for the linearized DtN map $\Lambda'_{c}$ and the linearized problem $(I_{1})$ in \eqref{eqn:problem_1}, by noticing that $u^{(1)} \big|_{\partial\Omega} = 0$, we obtain another Alessandrini type identity for the linearized DtN map $\Lambda'_{c}$,
\begin{equation}\label{eqn:identity}
\begin{aligned}
\int_{\Omega} c(x) ( u^{(0)} )^{m} \varphi \,\mathrm{d}x
&= \int_{\partial\Omega} ( \partial_{\nu} u^{(1)} ) \varphi \,\mathrm{d}S_{x} 
= \int_{\partial\Omega} \Lambda'_{c}(f) \, \varphi \,\mathrm{d}S_{x}.
\end{aligned}
\end{equation}

By choosing certain detect and test functions $u^{(0)}$ and $\varphi$ in \eqref{eqn:identity}, the uniqueness result for $m$-power nonlinearity was proved in \cite[Theorem 2.4]{LSX2022}. However such special solutions cannot lead to stable reconstruction of the potential function $c$ from the identity \eqref{eqn:identity} directly; see the discussion in \cite[Remark 2.5]{LSX2022}.
Inspired by the idea for $m = 2$ in \cite[Theorem 3.1]{LSX2022}, we derive an Alessandrini type identity with multiple detect functions for $m \geq 2$ cases, which is crucial in our current work. The key tool is the following identity derived from the principle of inclusion-exclusion (PIE) in combinatorics.

\begin{lemma}\label{lmm:identity_PIE}
Let $m \in \mathbb{N}_{+}$, define the index sets $U := \{1,2,\dots,m\}$ and $S$ be the subset of $U$. Then, given a set of numbers or functions $\{ w_{j} \}_{j = 1}^{m}$, we have
\begin{equation}\label{eqn:identity_PIE}
\prod_{j \in U} w_{j}
= \frac{1}{m!} \sum_{\emptyset \subsetneqq S \subseteq U} (-1)^{| U \setminus S |} \Big( \sum_{j \in S} w_{j} \Big)^{m}.
\end{equation}
Here, $U \setminus S$ is the relative complement of $S$ in $U$, and $|\cdot|$ is the cardinality of a set.
\par
Furthermore, given another $\ell \in \mathbb{N}_{+}$ such that $0 < \ell < m$, it holds that
\begin{equation}\label{eqn:identity_PIE_2}
0 = \frac{1}{m!} \sum_{\emptyset \subsetneqq S \subseteq U} (-1)^{| U \setminus S |} \Big( \sum_{j \in S} w_{j} \Big)^{\ell}.
\end{equation}
\end{lemma}

\begin{proof}
The detailed proof is given in \ref{sec:PIE}, with an overview of PIE.
We note that the identity \eqref{eqn:identity_PIE} is also proved by checking cancellation of each term carefully in \cite[Lemma 9]{LLPT2020}.
\end{proof}

In view of Lemma \ref{lmm:identity_PIE}, we denote $u_{j}^{(0)}$ as a solution of Helmholtz equation $(I_{0})$ in \eqref{eqn:problem_0} with Dirichlet boundary condition
\begin{equation*}
u_{j}^{(0)} = f_{j} \quad \text{on\ } \partial\Omega
\end{equation*}
for each $j \in U = \{1,2,\dots,m\}$, and denote $u_{S}^{(0)}$ as the combined solution of Helmholtz equation $(I_{0})$ in \eqref{eqn:problem_0} with Dirichlet boundary condition
\begin{equation*}
u_{S}^{(0)} = \sum_{j \in S} f_{j} \quad \text{on\ } \partial\Omega
\end{equation*}
for each non-empty subset $S \subseteq U$. By the linearity of Helmholtz equation $(I_{0})$ in \eqref{eqn:problem_0}, we have
\begin{equation*}
u_{S}^{(0)} = \sum_{j \in S} u_{j}^{(0)} \quad \text{in\ } \Omega.
\end{equation*}
Thus, from Lemma \ref{lmm:identity_PIE}, we get
\begin{equation}\label{eqn:temp_PIE}
\begin{aligned}
\int_{\Omega} c(x) \prod_{j \in U} u_{j}^{(0)} \varphi \,\mathrm{d}x
&= \frac{1}{m!} \sum_{\emptyset \subsetneqq S \subseteq U} (-1)^{| U \setminus S |} \int_{\Omega} c(x) \Big( \sum_{j \in S} u_{j}^{(0)} \Big)^{m} \varphi \,\mathrm{d}x \\
&= \frac{1}{m!} \sum_{\emptyset \subsetneqq S \subseteq U} (-1)^{| U \setminus S |} \int_{\Omega} c(x) ( u_{S}^{(0)} )^{m} \varphi \,\mathrm{d}x,
\end{aligned}
\end{equation}
with a test function $\varphi$ solving $\Delta \varphi + k^{2} \varphi = 0$ in $\Omega$.

Therefore, denote $u_{S}^{(1)}$ as the corresponding solution of the linearized equation
\begin{equation}\label{eqn:problem_1_S}
(I_{S}^{(1)}) ~
\left\{\begin{aligned}
\Delta u_{S}^{(1)} + k^{2} u_{S}^{(1)} &= c(x) ( u_{S}^{(0)} )^{m} & & \text{in\ } \Omega, \\
u_{S}^{(1)} &= 0 & & \text{on\ } \partial\Omega.
\end{aligned}\right.
\end{equation}
Then, for each $S$, we obtain
\begin{equation}\label{eqn:identity_combined}
\int_{\Omega} c(x) ( u_{S}^{(0)} )^{m} \varphi \,\mathrm{d}x
= \int_{\partial\Omega} ( \partial_{\nu} u_{S}^{(1)} ) \varphi \,\mathrm{d}S_{x},
\end{equation}
by the Alessandrini type identity \eqref{eqn:identity}. Following the definition of the linearized DtN map $\Lambda'_{c}$ \eqref{eqn:dtn_linear}, it means that, for each non-empty subset $S \subseteq U$,
\begin{equation}\label{eqn:dtn_linear_S}
\partial_{\nu} u_{S}^{(1)} \big|_{\partial\Omega}
= \Lambda'_{c} \big( u_{S}^{(0)} \big|_{\partial\Omega} \big)
= \Lambda'_{c} \Big( \sum_{j \in S} f_{j} \Big).
\end{equation}

Combining \eqref{eqn:temp_PIE}, \eqref{eqn:identity_combined} and \eqref{eqn:dtn_linear_S}, we have the following Alessandrini-PIE type identity, which is a combination formula of first order linearization.

\begin{lemma}\label{lmm:identity_S}
Let the sets $U$ and $S$, the solutions $\varphi$, $u_{j}^{(0)}$ and $u_{S}^{(1)}$ be defined as above. Given $m$ detect functions $\big\{ f_{j} = u_{j}^{(0)} \big|_{\partial\Omega} \big\}_{j=1}^{m}$, it holds that
\begin{equation}\label{eqn:identity_S}
\begin{aligned}
\int_{\Omega} c(x) \prod_{j \in U} u_{j}^{(0)} \varphi \,\mathrm{d}x
&= \frac{1}{m!} \sum_{\emptyset \subsetneqq S \subseteq U} (-1)^{| U \setminus S |} \int_{\partial\Omega} ( \partial_{\nu} u_{S}^{(1)} ) \varphi \,\mathrm{d}S_{x} \\
&= \frac{1}{m!} \sum_{\emptyset \subsetneqq S \subseteq U} (-1)^{| U \setminus S |} \int_{\partial\Omega}  \Lambda'_{c} \Big( \sum_{j \in S} f_{j} \Big) \, \varphi \,\mathrm{d}S_{x}.
\end{aligned}
\end{equation}
\end{lemma}

\begin{remark}
Take $m = 2$ for example and denote
\begin{equation*}
\begin{aligned}
u_{0} &= u_{\{1\}}^{(0)} = u_{1}^{(0)}, &\quad
v_{0} &= u_{\{2\}}^{(0)} = u_{2}^{(0)}, &\quad
w_{0} &= u_{\{1,2\}}^{(0)} = u_{1}^{(0)} + u_{2}^{(0)} = u_{0} + v_{0}, \\
u_{1} &= u_{\{1\}}^{(1)}, &\quad
v_{1} &= u_{\{2\}}^{(1)}, &\quad
w_{1} &= u_{\{1,2\}}^{(1)}.
\end{aligned}
\end{equation*}
Then, given detect functions $f_{1} = u_{0}|_{\partial\Omega}$ and $f_{2} = v_{0}|_{\partial\Omega}$, the identity \eqref{eqn:identity_S} becomes
\begin{equation*}
\begin{aligned}
\int_{\Omega} c(x) u_{0} v_{0} \varphi \,\mathrm{d}x
&= \frac{1}{2} \Big( \int_{\partial\Omega} ( \partial_{\nu} w_{1} ) \varphi \,\mathrm{d}S_{x}
- \int_{\partial\Omega} ( \partial_{\nu} u_{1} ) \varphi \,\mathrm{d}S_{x}
- \int_{\partial\Omega} ( \partial_{\nu} v_{1} ) \varphi \,\mathrm{d}S_{x} \Big) \\
&= \frac{1}{2} \Big( \int_{\partial\Omega} \Lambda'_{c}(f_{1}+f_{2}) \, \varphi \,\mathrm{d}S_{x}
- \int_{\partial\Omega} \Lambda'_{c}(f_{1}) \, \varphi \,\mathrm{d}S_{x}
- \int_{\partial\Omega} \Lambda'_{c}(f_{2}) \, \varphi \,\mathrm{d}S_{x} \Big),
\end{aligned}
\end{equation*}
by which the increasing stability and the algorithm for recovering the potential function $c$ are shown for the quadratic type nonlinearity, see \cite[Equation (3.3)]{LSX2022}.
\end{remark}

\begin{figure}[!htb]
\centering
\begin{tikzpicture}[>=latex,scale=1.0]
\draw[rounded corners=8,dashed,fill=gray!15!white]
 (0.6,-0.8) rectangle (2.6,0.8) node at ++(-0.3,-0.3) {$\Omega$};
\node  (f) at (0,0) {$f$};
\node  (u) at (2,0) {$u$};
\node (du) at (4,0) {$\partial_{\nu} u$};
\draw[->,thick,blue] (f) -- node[above] {$c$}  (u) node[pos=0.5,below] {\footnotesize$(I)$};
\draw[->,thick] (u) -- (du) node[pos=0.5,below] {\footnotesize$\partial_{\nu}$};
\draw[rounded corners=15,->,thick,blue] (f) |- (2,1.2) node[above] {$\Lambda_{c}$} -| (du);
\node at (2,-1.3) {\textbf{(i)} original problem};
\end{tikzpicture}
\hspace{4ex}
\begin{tikzpicture}[>=latex,scale=1.0]
\draw[rounded corners=8,dashed,fill=gray!15!white]
 (0.5,-0.8) rectangle (4.6,0.8) node at ++(-0.3,-0.3) {$\Omega$};
\node   (f) at (0,0) {$f$};
\node  (u0) at (2,0) {$u^{(0)}$};
\node  (u1) at (4,0) {$u^{(1)}$};
\node (du1) at (6,0) {$\partial_{\nu} u^{(1)}$};
\draw[->,thick]  (f) --  (u0) node[pos=0.5,below] {\footnotesize$(I_{0})$};
\draw[->,thick,blue] (u0) -- node[above] {$c$}  (u1) node[pos=0.5,below] {\footnotesize$(I_{1})$};
\draw[->,thick] (u1) -- (du1) node[pos=0.5,below] {\footnotesize$\partial_{\nu}$};
\draw[rounded corners=15,->,thick,blue] (f) |- (3,1.2) node[above] {$\Lambda'_{c}$} -| (du1);
\node at (3,-1.3) {\textbf{(ii)} linearized system};
\end{tikzpicture}\\
\vspace{1ex}
\begin{tikzpicture}[>=latex,scale=1.0]
\draw[rounded corners=8,dashed,fill=gray!15!white]
 (0.5,-1) rectangle (6.6,1) node at ++(-0.3,-0.3) {$\Omega$};
\node   (f) at (0,0) {$f_{j}$};
\draw[rounded corners=2,red] (-0.5,0.7) rectangle (2.45,-0.7) node at ++(-0.4,0.2) {\tiny$j \in S$};
\node  (u0) at (2,0) {$u_{j}^{(0)}$};
\node  (v0) at (4,0) {$u_{S}^{(0)}$};
\node  (v1) at (6,0) {$u_{S}^{(1)}$};
\node (dv1) at (8,0) {$\partial_{\nu} u_{S}^{(1)}$};
\draw[->,thick]  (f) --  (u0) node[pos=0.5,below] {\footnotesize$(I_{0})$};
\draw[->,thick,red] (u0) --  (v0) node[pos=0.5,below] {\footnotesize$\sum\limits_{j \in S}$};
\draw[->,thick,blue] (v0) -- node[above] {$c$}  (v1) node[pos=0.5,below] {\footnotesize$(I_{S}^{(1)})$};
\draw[->,thick] (v1) -- (dv1) node[pos=0.5,below] {\footnotesize$\partial_{\nu}$};
\draw[rounded corners=2] (-0.7,1.3) rectangle (8.7,-1.3) node at ++(-0.4,0.2) {\tiny$S \subseteq U$};
\draw[rounded corners=15,->,thick,blue] (f)++(0,0.7) |- (4,1.6) node[above] {$\Lambda'_{c}$} -| (dv1);
\node at (4,-1.8) {\textbf{(iii)} PIE type linearized system};
\end{tikzpicture}
\caption{\textbf{(i)} The Alessandrini type identity \eqref{eqn:identity_dtn} for DtN map $\Lambda_{c}$ \eqref{eqn:dtn}.
\textbf{(ii)} The Alessandrini type identity \eqref{eqn:identity} and
\textbf{(iii)} the Alessandrini-PIE type identity \eqref{eqn:identity_S} for linearized DtN map $\Lambda'_{c}$ \eqref{eqn:dtn_linear}.}
\label{fig:identity}
\end{figure}
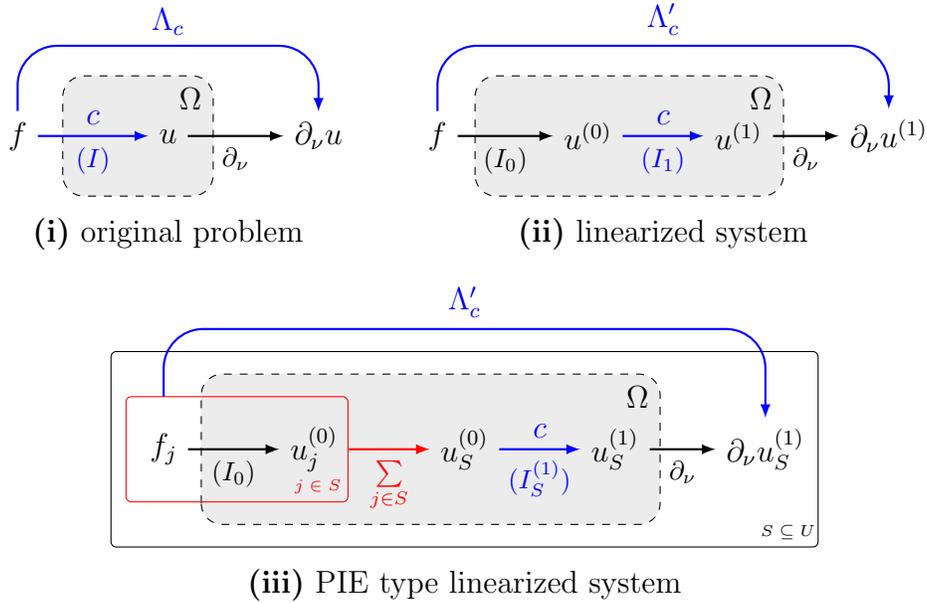

To succinctly explain the differences among these identities \eqref{eqn:identity_dtn}, \eqref{eqn:identity} and \eqref{eqn:identity_S}, there are three diagrams in Figure \ref{fig:identity}. The subfigure \textbf{(i)} shows the DtN map $\Lambda_{c}$ \eqref{eqn:dtn} for the original problem $(I)$ in \eqref{eqn:problem} and the Alessandrini type identity \eqref{eqn:identity_dtn}. The subfigure \textbf{(ii)} shows the linearized DtN map $\Lambda'_{c}$ \eqref{eqn:dtn_linear} for the linearized system $(I_{0})$--$(I_{1})$ in \eqref{eqn:problem_0}--\eqref{eqn:problem_1} and the correpsonding Alessandrini type identity \eqref{eqn:identity}. The subfigure \textbf{(iii)} explains the combined solutions $u_{S}^{(0)}$ and $u_{S}^{(1)}$ used in the Alessandrini-PIE type identity \eqref{eqn:identity_S} for the power type nonlinearities. Indeed, the identity \eqref{eqn:identity_S} is a combination formula of first order linearization (linearized DtN map $\Lambda'_{c}$). In Section \ref{sec:combined}, we will provide a numerical explanation for the combined solution $u_{S}^{(0)}$ and the identity \eqref{eqn:identity_S}, and explain the combination of boundary data via the principle of inclusion-exclusion.

\subsection{The reconstruction formula}\label{sec:reconstruct}

In this subsection, we introduce the special choice of detect and test functions $\big\{ u_{j}^{(0)} \big\}_{j=1}^{m}$ and $\varphi$, which will lead to a reconstruction formula for the inverse Schr\"{o}dinger potential problem by the Alessandrini-PIE type identity \eqref{eqn:identity_S}.

Given a vector $\xi \in \mathbb{R}^{n}$ in Fourier frequency space with $\xi \neq 0$.
In view of the identity \eqref{eqn:identity_S}, we set the solutions $\varphi$ and $u_{j}^{(0)}$ of the Helmholtz equation $(I_{0})$ in \eqref{eqn:problem_0} to be the complex exponential (CE) solutions in $\Omega$,
\begin{equation}\label{eqn:ce_sol}
\varphi(x) := \rme^{\bfi \zeta_{0} \cdot x}, \qquad
u_{j}^{(0)}(x) := \rme^{\bfi \zeta_{j} \cdot x}, \quad j \in U = \{1,2,\dots,m\},
\end{equation}
where the auxiliary complex vectors $\zeta_{0} \in \mathbb{C}^{n}$ and $\zeta_{j} \in \mathbb{C}^{n}$ satisfy
\begin{equation}\label{eqn:condition_k}
\zeta_{0} \cdot \zeta_{0} = k^{2}, \qquad
\zeta_{j} \cdot \zeta_{j} = k^{2}, \quad j \in U,
\end{equation}
and
\begin{equation}\label{eqn:condition_xi}
\zeta_{0} + \sum_{j \in U} \zeta_{j} = \xi.
\end{equation}

Noticing that
\begin{equation*}
\varphi \prod_{j \in U} u_{j}^{(0)}
= \rme^{\bfi ( \zeta_{0} + \sum_{j \in U} \zeta_{j} ) \cdot x}
= \rme^{\bfi \xi \cdot x}
\end{equation*}
in this stiuation and by the assumption $\operatorname{supp}(c) \subset \Omega$, the left hand side of the identity \eqref{eqn:identity_S} becomes the Fourier transform of the potential function $c$ at frequency $\xi$,
\begin{equation*}
\mathcal{F}[c](\xi)
:= \int_{\mathbb{R}^{n}} c(x) \rme^{\bfi \xi \cdot x} \,\mathrm{d}x
= \int_{\Omega} c(x) \rme^{\bfi \xi \cdot x} \,\mathrm{d}x,
\end{equation*}
which leads to the reconstruction formula via (inverse) Fourier transform.

\begin{lemma}\label{lmm:reconstruct}
The reconstruction formula for potential $c$ using the linearized DtN map $\Lambda_{c}'$,
\begin{equation}\label{eqn:reconstruct}
\begin{aligned}
\mathcal{F}[c](\xi)
&= \frac{1}{m!} \sum_{\emptyset \subsetneqq S \subseteq U} (-1)^{| U \setminus S |} \int_{\partial\Omega} ( \partial_{\nu} u_{S}^{(1)} ) \varphi \,\mathrm{d}S_{x} \\
&= \frac{1}{m!} \sum_{\emptyset \subsetneqq S \subseteq U} (-1)^{| U \setminus S |} \int_{\partial\Omega} \Lambda'_{c} \Big( \sum_{j \in S} f_{j} \Big) \, \varphi \,\mathrm{d}S_{x}.
\end{aligned}
\end{equation}
Here, the test and detect functions $\varphi$ and $u_{j}^{(0)}$ are defined in \eqref{eqn:ce_sol}, $f_{j} = u_{j}^{(0)} \big|_{\partial\Omega}$, $u_{S}^{(0)} = \sum_{j \in S} u_{j}^{(0)}$, and $u_{S}^{(1)}$ is the solution of the linearized problem $(I_{S}^{(1)})$ in \eqref{eqn:problem_1_S}, for each $j \in U$ and non-empty subset $S \subseteq U$.
\end{lemma}

The above Lemma \ref{lmm:reconstruct} allows us to carry out the stability estimate (Section \ref{sec:main_results}) and reconstruction algorithm (Section \ref{sec:algorithm}) of the inverse Schr\"{o}dinger potential problem \eqref{eqn:problem} with integer power type nonlinearities. A particular choice of the CE solutions \eqref{eqn:ce_sol} which satisfy the conditions \eqref{eqn:condition_k} and \eqref{eqn:condition_xi} wil be given below. Indeed, such a set of solutions has been also used in \cite{LSX2022} to prove the increasing stability in the higher order linearization method. Given $\xi \in \mathbb{R}^{n} \setminus \{0\}$, define an orthonormal basis of $\mathbb{R}^{n}$ by
\begin{equation*}
\left\{ e_{1} := \frac{\xi}{|\xi|}, e_{2}, \dots, e_{n} \right\},
\end{equation*}
for $n \geq 2$.
Thus, for different index $m$, the complex vectors $\zeta_{0}$ and $\{ \zeta_{j} \}_{j=1}^{m}$ are chosen as below.
\begin{enumerate}[(1)]
\item \textbf{Odd $m$:} for each $j = 0,1,2,\dots,m$,
\begin{equation}\label{eqn:vec_odd_m}
\zeta_{j} := \left\{\begin{aligned}
&\frac{1}{m+1} \Big( |\xi| e_{1} + \Xi_{\textrm{odd}} e_{2} \Big), & & \text{for odd $j$}, \\
&\frac{1}{m+1} \Big( |\xi| e_{1} - \Xi_{\textrm{odd}} e_{2} \Big), & & \text{for even $j$ or $j = 0$}, \\
\end{aligned}\right.
\end{equation}
where the notation $\Xi_{\textrm{odd}}$ is given by
\begin{equation}\label{eqn:Xi_odd_m}
\Xi_{\textrm{odd}} := \left\{\begin{aligned}
\sqrt{(m+1)^{2} k^{2} - |\xi|^{2}},& &\quad |\xi| \leq (m+1)k, \\
\bfi \sqrt{|\xi|^{2} - (m+1)^{2} k^{2}},& &\quad |\xi| > (m+1)k.
\end{aligned}\right.
\end{equation}
\item \textbf{Even $m$:} for each $j = 1,2,\dots,m$,
\begin{equation}\label{eqn:vec_even_m}
\begin{aligned}
\zeta_{0} &:= k e_{1}, \\
\zeta_{j} &:= \left\{\begin{aligned}
&\frac{1}{m} \Big( (|\xi|-k) e_{1} + \Xi_{\textrm{even}} e_{2} \Big), & & \text{for odd $j$}, \\
&\frac{1}{m} \Big( (|\xi|-k) e_{1} - \Xi_{\textrm{even}} e_{2} \Big), & & \text{for even $j$}, \\
\end{aligned}\right.
\end{aligned}
\end{equation}
where the notation $\Xi_{\textrm{even}}$ is given by
\begin{equation}\label{eqn:Xi_even_m}
\Xi_{\textrm{even}} := \left\{\begin{aligned}
\sqrt{m^{2} k^{2} - (|\xi|-k)^{2}},& &\quad |\xi| \leq (m+1)k, \\
\bfi \sqrt{(|\xi|-k)^{2} - m^{2} k^{2}},& &\quad |\xi| > (m+1)k.
\end{aligned}\right.
\end{equation}
\end{enumerate}

We will see in the following sections that this set of CE solutions leads to the increasing stability estimate and reconstruction algorithm for our linearized inverse Schr\"{o}dinger potential problem.

\section{The increasing stability theorems}\label{sec:main_results}

Under the definitions and assumptions in Section \ref{sec:problem}, we state and prove the main results of increasing stability in this paper.

The first theorem shows that, via the measurements of the linearized DtN map $\Lambda'_{c}$, one can estimate the Fourier coefficient $\mathcal{F}[c](\xi)$ of the potential function $c$ at frequency $\xi \in \mathbb{R}^{n}$ in a Lipschitz stable way, if $|\xi| \leq (m+1)k$.

\begin{theorem}\label{thm:lipschitz}
Let $m, n \geq 2$ be integers, and $\Omega \subset \mathbb{R}^{n}$ be an open bounded domain with $C^{\infty}$ boundary $\partial\Omega$. Assume that the potential function $c \in L^{\infty}(\Omega)$ with $\operatorname{supp}(c) \subset \Omega$, and the wavenumber $k > 1$ such that $k^{2}$ is not a Dirichlet eigenvalue of $-\Delta$ in $\Omega$.
Denote
\begin{equation*}
\varepsilon := \| \Lambda'_{c} \|_{\mathcal{N}}
\end{equation*}
as the operator norm \eqref{eqn:dtn_linear_norm} of the linearized DtN map $\Lambda'_{c}$. Then the Fourier transform of the potential function $c$ at frequency $\xi \neq 0$ satisfies
\begin{equation}\label{eqn:lipschitz}
| \mathcal{F}[c](\xi) |
\leq C(m,\Omega) k^{2m} \varepsilon, \quad \text{for\ } |\xi| \leq (m+1)k.
\end{equation}
\end{theorem}

\begin{remark}
Comparing to \cite{ILX2020} for linear Schr\"{o}dinger equation, in which the Lipschitz stability of recovering the Fourier coefficients holds only for $|\xi| \leq 2k$, 
Theorem \ref{thm:lipschitz} yields a larger interval $|\xi|\leq (m+1)k$ in which Fourier coefficients of the unknown potential function can be recovered with a Lipschitz stability. This observation is in consistent with the result in \cite[Theorem 2.1]{LSX2022} by the higher order linearization method.
\end{remark}

The second theorem gives an increasing stability estimate of the potential function $c$ by linearized DtN map $\Lambda'_{c}$, which generalizes the previous result \cite[Theorem 3.1]{LSX2022} from the quadratic nonlinear term to general integer power type ones.

\begin{theorem}\label{thm:stability}
Under the same notations and assumptions in Theorem \ref{thm:lipschitz}, and without loss of generality, we furthermore assume that $0 \in \Omega$ and $R := \sup_{x \in \overline{\Omega}} |x| \leq \frac{1}{2}$.
If $\varepsilon < 1$ and $\|c\|_{H^{1}(\Omega)} \leq M$, the following stability estimate holds
\begin{equation}\label{eqn:stability}
\|c\|_{L^{2}(\Omega)}^{2}
\leq C(m,\Omega) \left( k^{\tau} \varepsilon^{2} + E^{\tau} \varepsilon + \frac{M}{1+m^{2}k^{2}+E^{2}} \right)
\end{equation}
with $E := -\ln\varepsilon$ and $\tau := 4m+n$.
\end{theorem}

\begin{remark}
Note that in \cite[Theorem 2.1 and Theorem 2.2]{LSX2022} a similar result holds for the higher order linearization with respect to small Dirichlet boundary data $f$. We highlight that the difference is that we use the linearized DtN map $\Lambda'_{c}$ with respect to small potential function $c$ in current work. In particular, the Dirichlet boundary data of current work is not associated with different varying small variables, which is necessary for the higher order linearization method.
\end{remark}


\subsection{Proof of Theorem \ref{thm:lipschitz}}

We notice that, by choosing the complex vectors $\zeta_{0}$ and $\{ \zeta_{j} \}_{j=1}^{m}$ as \eqref{eqn:vec_odd_m} and \eqref{eqn:vec_even_m}, the complex exponential (CE) solutions
\begin{equation*}
\varphi(x) = \rme^{\bfi \zeta_{0} \cdot x}, \quad\text{and}\quad
u_{j}^{(0)}(x) = \rme^{\bfi \zeta_{j} \cdot x}, \quad\text{for each\ } j \in U = \{1,2,\dots,m\},
\end{equation*}
defined in \eqref{eqn:ce_sol} are indeed the plane waves if $|\xi| \leq (m+1)k$.
This property leads to the stable recovery in low Fourier frequency mode as shown in Theorem \ref{thm:lipschitz} and we prove below.

\begin{proof}[Proof of Theorem \ref{thm:lipschitz}]
Given Fourier frequency $\xi \in \mathbb{R}^{n}$ with $0 < |\xi| \leq (m+1)k$.
Recall $p = \frac{n}{2} + 1$ in Proposition \ref{prp:wellposed}.
When the CE solutions $\varphi$ and $\{ u_{j}^{(0)} \}_{j=1}^{m}$ are defined by \eqref{eqn:ce_sol}, the complex vectors $\zeta_{0}$ and $\{ \zeta_{j} \}_{j=1}^{m}$ are chosen as \eqref{eqn:vec_odd_m} and \eqref{eqn:vec_even_m}, we obtain their $W^{2,p}$-norm estimate under the condition $k \geq 1$,
\begin{equation*}
\big\| u_{j}^{(0)} \big\|_{W^{2,p}(\Omega)}
\leq \sum_{|\alpha| \leq 2} \big\| \mathrm{D}^{\alpha} u_{j}^{(0)} \big\|_{L^{p}(\Omega)}
\leq C(\Omega) (1+k+k^{2})
\leq C(\Omega) k^{2}
\end{equation*}
for any $j \in U$, where $C(\Omega)$ is a constant only depending on $\Omega$.
Therefore, for any non-empty subset $S \subseteq U$,
\begin{equation*}
\big\| u_{S}^{(0)} \big\|_{W^{2,p}(\Omega)}
\leq \sum_{j \in S} \big\| u_{j}^{(0)} \big\|_{W^{2,p}(\Omega)}
\leq C(\Omega) k^{2} |S|.
\end{equation*}
Meanwhile, it holds that
\begin{equation*}
\left\| \varphi \right\|_{L^{\infty}(\partial\Omega)}
= \sup_{x \in \partial\Omega} \left| \rme^{\bfi \zeta_{0} \cdot x} \right|
= 1.
\end{equation*}
According to the definition of the linearized DtN map $\Lambda'_{c}$ in \eqref{eqn:dtn_linear_norm}, by using the H\"{o}lder inequality \cite[Theorem 2.4]{AF2003}, the Sobolev embedding theorem \cite[Theorem 4.12]{AF2003} and the trace theorem \cite[Theorem 7.39]{AF2003}, we have
\begin{equation}\label{eqn:linearcomb}
\begin{aligned}
\left| \int_{\partial\Omega} \Lambda'_{c} \Big( \sum_{j \in S} f_{j} \Big) \, \varphi \,\mathrm{d}S_{x} \right|
&\leq \Big\| \Lambda'_{c} \Big( \sum_{j \in S} f_{j} \Big) \Big\|_{L^{p}(\partial\Omega)} \| \varphi \|_{L^{\frac{p}{p-1}}(\partial\Omega)} \\
&\leq C(\Omega) \Big\| \Lambda'_{c} \Big( \sum_{j \in S} f_{j} \Big) \Big\|_{W^{1-\frac{1}{p},p}(\partial\Omega)} \| \varphi \|_{L^{\infty}(\partial\Omega)} \\
&\leq C(\Omega) \, \varepsilon \, \Big\| \sum_{j \in S} f_{j} \Big\|_{W^{2-\frac{1}{p},p}(\partial\Omega)}^{m} \| \varphi \|_{L^{\infty}(\partial\Omega)} \\
&\leq C(m,\Omega) \, \varepsilon \big\| u_{S}^{(0)} \big\|_{W^{2,p}(\Omega)}^{m} \| \varphi \|_{L^{\infty}(\partial\Omega)} \\
&\leq C(m,\Omega) \, \varepsilon \, k^{2m} |S|^{m}.
\end{aligned}
\end{equation}
Here $C(m,\Omega)$ is another constant depending on $m$ and $\Omega$.
From Lemma \ref{lmm:reconstruct}, by summing up over all non-empty subsets $S \subseteq U$, we have
\begin{equation*}
\begin{aligned}
| \mathcal{F}[c](\xi) |
&\leq \frac{1}{m!} \sum_{\emptyset \subsetneqq S \subseteq U} \left| \int_{\partial\Omega} \Lambda'_{c} \Big( \sum_{j \in S} f_{j} \Big) \, \varphi \,\mathrm{d}S_{x} \right|
\leq C(m,\Omega) \,\varepsilon \, \frac{k^{2m}}{m!} \sum_{\emptyset \subsetneqq S \subseteq U} |S|^{m} \\
&\leq C(m,\Omega) \, \varepsilon \, \frac{k^{2m}}{m!} \sum_{j=1}^{m} \binom{m}{j} j^{m}
\leq C(m,\Omega) k^{2m} \varepsilon,
\end{aligned}
\end{equation*}
which concludes the proof of Theorem \ref{thm:lipschitz}.
\end{proof}


\subsection{Proof of Theorem \ref{thm:stability}}

\begin{proof}[Proof of Theorem \ref{thm:stability}]
Recall that we assume $0 \in \Omega$ and denote $R := \sup_{x \in \Omega} |x|$.

While Theorem \ref{thm:lipschitz} gives the point-wise estimate for Fourier coefficients of the potential function $c$ in low Fourier frequency mode $\{ \xi \in \mathbb{R}^{n} \,:\, |\xi| \leq (m+1)k \}$,
the a-priori assumption $\|c\|_{H^{1}(\Omega)} \leq M$ leads to its $L^{2}$-norm estimate in high Fourier frequency mode.
That is, given any $\rho > 0$ and $\{ \xi \in \mathbb{R}^{n} \,:\, |\xi| > \rho \}$, we have
\begin{equation}\label{eqn:high_freq}
\begin{aligned}
\int_{|\xi| > \rho} | \mathcal{F}[c](\xi) |^{2} \,\mathrm{d}\xi
&\leq \frac{1}{1+\rho^{2}} \int_{|\xi| > \rho} (1+|\xi|^{2}) | \mathcal{F}[c](\xi) |^{2} \,\mathrm{d}\xi \\
&\leq \frac{1}{1+\rho^{2}} \int_{\mathbb{R}^{n}} (1+|\xi|^{2}) |\mathcal{F}[c](\xi) |^{2} \,\mathrm{d}\xi \\
&\leq \frac{C(\Omega)}{1+\rho^{2}} \|c\|_{H^{1}(\Omega)}^{2}
\leq C(\Omega) \frac{M^{2}}{1+\rho^{2}},
\end{aligned}
\end{equation}
by the extension theorem \cite[Theorem 5.28]{AF2003} and the assumption $\operatorname{supp}(c) \subset \Omega$.

To combine these two estimates \eqref{eqn:lipschitz} and \eqref{eqn:high_freq}, denote that $E := -\ln \varepsilon > 0$ ($\varepsilon < 1$), the following two cases are considered,
\begin{enumerate}[(A)]
\item $k > E$ (i.e. $\varepsilon > \rme^{-k}$), and
\item $k \leq E$ (i.e. $\varepsilon \leq \rme^{-k}$).
\end{enumerate}

For \textbf{Case (A)}, since $k > E$ is relatively large, by setting $\rho = (m+1)k$, we have
\begin{equation}\label{eqn:case_A}
\begin{aligned}
\|c\|_{L^{2}(\Omega)}^{2}
&\leq \left( \int_{|\xi| \leq (m+1)k} + \int_{|\xi| > (m+1)k} \right) | \mathcal{F}[c](\xi) |^{2} \,\mathrm{d}\xi \\
&\leq C(m,\Omega) \left( (m+1)^{n} k^{n} k^{4m} \varepsilon^{2} + \frac{M^{2}}{1+(m+1)^{2}k^{2}} \right) \\
&\leq C(m,\Omega) \left( k^{4m+n} \varepsilon^{2} + \frac{M^{2}}{1+m^{2}k^{2}+E^{2}} \right).
\end{aligned}
\end{equation}

For \textbf{Case (B)}, given $\xi \in \mathbb{R}^{n}$ with $|\xi| > (m+1)k$, we apply Lemma \ref{lmm:reconstruct} again to estimate Fourier coefficient $\mathcal{F}[c](\xi)$.
Notice that, the corresponding CE solutions $\varphi$ and $u_{j}^{(0)}$ with $\zeta_{0}$ and $\{ \zeta_{j} \}_{j=1}^{m}$ chosen as \eqref{eqn:vec_odd_m} and \eqref{eqn:vec_even_m} are growing exponentially when $|\xi| > (m+1)k$.

Recall $p = \frac{n}{2} + 1$ in Proposition \ref{prp:wellposed}.
Here we discuss in detail for different index $m$.

\noindent (B.1) \textbf{Odd $m$:} In this case $\Xi_{\textrm{odd}} = \bfi \sqrt{|\xi|^{2} - (m+1)^{2} k^{2}}$ in \eqref{eqn:Xi_odd_m}, which leads to the estimate
\begin{equation*}
\big\| u_{j}^{(0)} \big\|_{W^{2,p}(\Omega)}
\leq \sum_{|\alpha| \leq 2} \big\| \mathrm{D}^{\alpha} u_{j}^{(0)} \big\|_{L^{p}(\Omega)}
\leq C(\Omega) k^{2} \rme^{\frac{|\Xi_{\textrm{odd}}|}{m+1}R},
\end{equation*}
for any $j \in U = \{1,2,\dots,m\}$ and $k \geq 1$.
Therefore, for any non-empty subset $S \subseteq U$,
\begin{equation*}
\big\| u_{S}^{(0)} \big\|_{W^{2,p}(\Omega)}
\leq \sum_{j \in S} \big\| u_{j}^{(0)} \big\|_{W^{2,p}(\Omega)}
\leq C(\Omega) k^2 |S| \rme^{\frac{|\Xi_{\textrm{odd}}|}{m+1}R}.
\end{equation*}
Meanwhile, the estimate of test function $\varphi$ holds that
\begin{equation*}
\left\| \varphi \right\|_{L^{\infty}(\partial\Omega)}
= \sup_{x \in \partial\Omega} \left| \rme^{\mathbf{i} \zeta_{0} \cdot x} \right|
= \rme^{\frac{|\Xi_{\textrm{odd}}|}{m+1}R}.
\end{equation*}
Then, analogous to \eqref{eqn:linearcomb}, we have
\begin{equation*}
\begin{aligned}
\left| \int_{\partial\Omega} \Lambda'_{c} \Big( \sum_{j \in S} f_{j} \Big) \, \varphi \,\mathrm{d}S_{x} \right|
&\leq C(m,\Omega) \, \varepsilon \big\| u_{S}^{(0)} \big\|_{W^{2,p}(\Omega)}^{m} \| \varphi \|_{L^{\infty}(\partial\Omega)} \\
&\leq C(m,\Omega) \, \varepsilon \, k^{2m} |S|^{m} \rme^{|\Xi_{\textrm{odd}}|R} .
\end{aligned}
\end{equation*}
From Lemma \ref{lmm:reconstruct}, by summing up over all non-empty subsets $S \subseteq U$, we have
\begin{equation*}
| \mathcal{F}[c](\xi) |
\leq C(m,\Omega) k^{2m} \rme^{|\Xi_{\textrm{odd}}|R} \varepsilon, \quad \text{for\ } |\xi| > (m+1)k.
\end{equation*}

By setting $\rho = \sqrt{(m+1)^{2}k^{2}+\frac{E^{2}}{4R^{2}}} > (m+1)k$, we have
\begin{equation*}
\sqrt{ \rho^{2} - (m+1)^{2} k^{2} } = \frac{E}{2R}
\end{equation*}
and
\begin{equation*}
\rho^{n} \leq \left( (m+1)^{2} E^{2} + \frac{E^{2}}{4R^{2}} \right)^\frac{n}{2}
\leq C(m,\Omega) E^{n},
\end{equation*}
under the condition that $k \leq E$.
Thus, the $L^{2}$-norm estimate for the Fourier coefficients of the potential function $c$ in intermediate Fourier frequency mode $\{ \xi \in \mathbb{R}^{n} \,:\, (m+1)k < |\xi| \leq \rho \}$ satisfies
\begin{equation}\label{eqn:holder_odd_m}
\begin{aligned}
\int_{(m+1)k < |\xi| \leq \rho} | \mathcal{F}[c](\xi) |^{2} \,\mathrm{d}\xi
&\leq C(m,\Omega) k^{4m} \varepsilon^{2} \int_{(m+1)k < |\xi| \leq \rho} \rme^{2|\Xi_{\textrm{odd}}|R} \,\mathrm{d}\xi \\
&\leq C(m,\Omega) k^{4m} \rme^{2R\sqrt{\rho^{2}-(m+1)^{2}k^{2}}} \rho^{n} \varepsilon^{2} \\
&\leq C(m,\Omega) k^{4m} \rme^{E} E^{n} \varepsilon^{2} \\
&\leq C(m,\Omega) E^{4m+n} \varepsilon.
\end{aligned}
\end{equation}
Combining \eqref{eqn:lipschitz}, \eqref{eqn:high_freq} and \eqref{eqn:holder_odd_m}, we obtain
\begin{equation}\label{eqn:case_B_odd_m}
\begin{aligned}
\|c\|_{L^{2}(\Omega)}^{2}
&\leq \left( \int_{|\xi| \leq (m+1)k} + \int_{(m+1)k < |\xi| \leq \rho} + \int_{|\xi| > \rho} \right) | \mathcal{F}[c](\xi) |^{2} \,\mathrm{d}\xi \\
&\leq C(m,\Omega) \left( (m+1)^{n} k^{n} k^{4m} \varepsilon^{2} + E^{4m+n} \varepsilon + \frac{M^{2}}{1+(m+1)^{2}k^{2}+\frac{E^{2}}{4R^{2}}} \right) \\
&\leq C(m,\Omega) \left( k^{4m+n} \varepsilon^{2} + E^{4m+n} \varepsilon + \frac{M^{2}}{1+m^{2}k^{2}+E^{2}} \right)
\end{aligned}
\end{equation}
since the condition $R \leq \frac{1}{2}$.

\noindent (B.2) \textbf{Even $m$:} The proof follows the same routine as the odd $m$ case, with a substitution of $\Xi_{\textrm{odd}}$ by $\Xi_{\textrm{even}} = \bfi \sqrt{(|\xi|-k)^{2} - m^{2} k^{2}}$ in \eqref{eqn:Xi_even_m}.
But we note that $\varphi(x) = \rme^{\bfi k e_{1} \cdot x}$ is always a plane wave.%

In this case, by setting a different $\rho = k + \sqrt{m^{2}k^{2}+\frac{E^{2}}{4R^{2}}} > (m+1)k$, we have
\begin{equation*}
\sqrt{ ( \rho - k )^{2} - m^{2} k^{2} } = \frac{E}{2R}
\end{equation*}
and
\begin{equation*}
\rho^{n} < \left( k + m k + \frac{E}{2R} \right)^{n}
\leq C(\Omega,m) E^{n},
\end{equation*}
under the condition that $k \leq E$. Thus, the estimate for Fourier coefficients in intermediate Fourier frequency domain reads that
\begin{equation}\label{eqn:holder_even_m}
\begin{aligned}
\int_{(m+1)k < |\xi| \leq \rho} | \mathcal{F}[c](\xi) |^{2} \,\mathrm{d}\xi
&\leq C(m,\Omega) k^{4m} \varepsilon^{2} \int_{(m+1)k < |\xi| \leq \rho} \rme^{2|\Xi_{\textrm{even}}|R} \,\mathrm{d}\xi \\
&\leq C(m,\Omega) k^{4m} \rme^{2R\sqrt{(\rho-k)^{2}-m^{2}k^{2}}} \rho^{n} \varepsilon^{2} \\
&\leq C(m,\Omega) k^{4m} \rme^{E} E^{n} \varepsilon^{2} \\
&\leq C(m,\Omega) E^{4m+n} \varepsilon.
\end{aligned}
\end{equation}
Combining the estimates \eqref{eqn:lipschitz}, \eqref{eqn:high_freq}, \eqref{eqn:holder_even_m}, and noticing that
\begin{equation*}
\rho^{2} > m^{2} k^{2} + \frac{E^{2}}{4R^{2}}
\geq m^{2} k^{2} + E^{2},
\end{equation*}
we derive the same $L^{2}$-norm estimate for even $m$ as that in \eqref{eqn:case_B_odd_m} for odd $m$.

Therefore we conclude the proof by combining \eqref{eqn:case_A} and \eqref{eqn:case_B_odd_m}.
\end{proof}

\begin{remark}
We note that Theorem \ref{thm:stability} is an increasing stability result: when $k$ is sufficiently large, the stability estimate improves from a logarithmic type to a H\"{o}lder one.
Here we only consider the case when $k > E$.
From \eqref{eqn:case_A}, we have the following upper bound estimate
\begin{equation*}
\|c\|_{L^{2}(\Omega)}^{2}
\leq \omega(k) := C \left( k^{4m+n} \varepsilon^{2} + \frac{1}{k^2} \right),
\end{equation*}
for all $k > E$.
Notice that $\omega(k)$ achieves its minimum in $(0,\infty)$ at $k^{\star} = \varepsilon^{-\frac{2}{4m+n+2}}$.
Therefore if $E \leq k^{\star}$, we have
\begin{equation*}
\|c\|_{L^{2}(\Omega)}^{2}
\leq \omega(k^{\star}) = C \varepsilon^{\frac{4}{4m+n+2}},
\end{equation*}
and if $E > k^{\star}$, we have
\begin{equation*}
\|c\|_{L^2(\Omega)}^{2}
\leq \omega(E) = C \left( E^{4m+n} \varepsilon^{2} + E^{-2} \right)
\leq C \left( E^{4m+n} \varepsilon^{2} + \varepsilon^{\frac{4}{4m+n+2}} \right).
\end{equation*}
We conclude that when $k$ is sufficiently large, the stability estimate is of H\"{o}lder type, which improves the logarithmic result when $k = 0$, see \cite[Remark 3.2]{LLLS2021}.
\end{remark}


\section{The reconstruction algorithm}\label{sec:algorithm}

Recalling Lemma \ref{lmm:reconstruct}, given the CE solutions $\varphi$, $u_{j}^{(0)}$ in \eqref{eqn:ce_sol} and the choice of complex vectors $\zeta_{0}$, $\zeta_{j}$ in \eqref{eqn:vec_odd_m}--\eqref{eqn:Xi_even_m}, $j \in U = \{ 1,2,\dots,m \}$, we obtain the reconstruction formula \eqref{eqn:reconstruct} for the Fourier coefficient $\mathcal{F}[c](\xi)$ of the potential function $c$ at frequency $\xi$. That is,
\begin{equation*}
\begin{aligned}
\mathcal{F}[c](\xi)
= \int_{\Omega} c(x) \rme^{\bfi \xi \cdot x} \,\mathrm{d}x
&= \frac{1}{m!} \sum_{\emptyset \subsetneqq S \subseteq U} (-1)^{| U \setminus S |} \int_{\partial\Omega} ( \partial_{\nu} u_{S}^{(1)} ) \rme^{\bfi \zeta_{0} \cdot x} \,\mathrm{d}S_{x} \\
&= \frac{1}{m!} \sum_{\emptyset \subsetneqq S \subseteq U} (-1)^{| U \setminus S |} \int_{\partial\Omega} \Lambda'_{c} \Big( \sum_{j \in S} f_{j} \Big) \, \varphi \,\mathrm{d}S_{x},
\end{aligned}
\end{equation*}
where the solution $u_{S}^{(1)}$ satisfies the linearized problem \eqref{eqn:problem_1_S}, and $\partial_{\nu} u_{S}^{(1)} \big|_{\partial\Omega}$ is the corresponding Neumann boundary data on $\partial\Omega$.
Then, by choosing $\xi \neq 0$, we aim to recover all the Fourier coefficients $\mathcal{F}[c](\xi)$ of the potential function $c$ satisfying $|\xi| \leq (m+1)k$; see Theorem \ref{thm:lipschitz}. The larger the wavenumber $k$ is, the more Fourier coefficients can be recovered.

For any $0 \neq \xi \in \mathbb{R}^{n}$ in Fourier frequency space, we denote the vector $\xi := \kappa \hat{y}$ with the length $\kappa$ and angle $\hat{y}$ (a unit vector), i.e.
\begin{equation*}
\kappa := |\xi| \in \mathbb{R}, \quad
\hat{y} := \frac{\xi}{|\xi|} = e_{1} \in \mathbb{S}^{n-1}.
\end{equation*}
Then there exists another unit vector $\hat{z} := e_{2} \in \mathbb{S}^{n-1}$ such that $\hat{y} \cdot \hat{z} = 0$.
Similar to \cite{ILX2020, LSX2022}, we define the following discrete sets of lengths and angles of the vectors in the phase space. The discrete and finite length set is defined by
\begin{equation*}
\{\kappa_{r}\}_{r=1}^{\mathcal{R}} \subset (0, (m+1)k ]
\quad \text{for any fixed $k$}.
\end{equation*}
Here $(m+1)k$ is the maximum length of the vector $\xi$. Two angle (unit vector) sets are defined by
\begin{equation*}
\{\hat{y}_{\theta}\}_{\theta=1}^{\Theta} \subseteq \mathbb{S}^{n-1}
\quad \text{and} \quad
\{\hat{z}_{\theta}\}_{\theta=1}^{\Theta} \subseteq \mathbb{S}^{n-1},
\end{equation*}
which satisfy $\hat{y}_{\theta} \cdot \hat{z}_{\theta} = 0$ for each $\theta \in \{ 1,2,\dots,\Theta \}$.
For the sake of convenience, we also denote
\begin{equation*}
\xi^{\langle r;\theta \rangle} := \kappa_{r} \hat{y}_{\theta}
\quad \text{with} \quad
\kappa_{r} := |\xi^{\langle r;\theta \rangle}|
\quad \text{and} \quad
\hat{y}_{\theta} := \frac{\xi^{\langle r;\theta \rangle}}{|\xi^{\langle r;\theta \rangle}|}.
\end{equation*}
To address the inverse Fourier transform, we can construct a numerical quadrature rule by a suitable choice of the weights $\omega^{\langle r;\theta \rangle}$ according to the vector $\xi^{\langle r;\theta \rangle}$.

Referring to Lemma \ref{lmm:reconstruct}, we choose the complex vectors $\zeta_{0}$ and $\zeta_{j}$, $j \in U$ as follows.
\begin{enumerate}[(1)]
\item \textbf{Odd $m$:} for each $j = 0,1,2,\dots,m$, the complex vectors are given by
\begin{equation*}
\zeta_{j}^{\langle r;\theta \rangle} = \left\{~\begin{aligned}
& \frac{1}{m+1} \big( \kappa_{r} \hat{y}_{\theta} + \sqrt{(m+1)^{2} k^{2} - \kappa_{r}^{2}} \, \hat{z}_{\theta} \big), & & \text{for odd $j$}, \\
& \frac{1}{m+1} \big( \kappa_{r} \hat{y}_{\theta} - \sqrt{(m+1)^{2} k^{2} - \kappa_{r}^{2}} \, \hat{z}_{\theta} \big), & & \text{for even $j$ or $j = 0$}. \\
\end{aligned}\right.
\end{equation*}
\item \textbf{Even $m$:} for each $j = 1,2,\dots,m$, the complex vectors are given by
\begin{equation*}
\begin{aligned}
\zeta_{0}^{\langle r;\theta \rangle} &= k \hat{y}_{\theta}, \\
\zeta_{j}^{\langle r;\theta \rangle} &= \left\{~\begin{aligned}
& \frac{1}{m} \big( (\kappa_{r}-k) \hat{y}_{\theta} + \sqrt{m^{2} k^{2} - (\kappa_{r}-k)^{2}} \, \hat{z}_{\theta} \big), & & \text{for odd $j$}, \\
& \frac{1}{m} \big( (\kappa_{r}-k) \hat{y}_{\theta} - \sqrt{m^{2} k^{2} - (\kappa_{r}-k)^{2}} \, \hat{z}_{\theta} \big), & & \text{for even $j$}. \\
\end{aligned}\right.
\end{aligned}
\end{equation*}
\end{enumerate}

Finally, we summarize our reconstruction algorithm below (\textbf{Algorithm 1}), which generalizes those in \cite{ILX2020, LSX2022} towards general integer power type nonlinearities. Noticing that, since the nonlinearity index is $m$, we have to solve the nonlinear Schr\"{o}dinger potential problem $2^{m}-1$ times at each iteration, which causes considerably high computational costs in computing the forward and linearized problems and we refer to \cite{FT2005, WZ2018, XB2010, YL2017} for extended discussion on varies numerical schemes.

\begin{table}[!htb]
\centering
\begin{tabular}{p{\textwidth}}
\toprule
\textbf{Algorithm 1: Reconstruction Algorithm for the Linearized Schr\"{o}dinger Potential Problem} \\ \textbf{(Integer power type nonlinearity)} \\%
\midrule
\textbf{Input:} %
$k$, %
$m$, %
$U = \{ 1,2,\dots,m \}$, %
$\beta$, %
$\{\kappa_{r}\}_{r=1}^{\mathcal{R}}$, %
$\{\hat{y}_{\theta}\}_{\theta=1}^{\Theta}$, %
$\{\hat{z}_{\theta}\}_{\theta=1}^{\Theta}$ and %
$\{\omega^{\langle r;\theta \rangle}\}$. \\[5pt]%
\textbf{Output:} %
Approximated Potential $c_{\textrm{inv}} = c^{\langle \mathcal{R}+1;1 \rangle}$. \\[-15pt]%
\begin{enumerate}[1:]
  \item Set $c^{\langle 1;1 \rangle} := 0$; %
  \item \textbf{For} $r = 1,2,\dots,\mathcal{R}$, (length~updating) %
  \item ~ \textbf{For} $\theta = 1,2,\dots,\Theta$, (angle~updating) %
  \item ~ ~ \textbf{For} each non-empty subset $S \subseteq U$, (subset~updating)%
  \item ~ ~ ~ Choose $u_{S}^{(0)} := \sum_{j \in S} u_{j}^{(0)}$ with $u_{j}^{(0)} := \beta \exp\{ \bfi \zeta_{j}^{\langle r;\theta \rangle} \cdot x \}$ and $j \in S$; %
  \item ~ ~ ~ Measure Neumann data $\partial_{\nu} u_{S} |_{\partial\Omega}$ of the original problem \eqref{eqn:problem} %
  \item[] ~ ~ ~ ~ while Dirichlet data $u_{S} |_{\partial\Omega} := u_{S}^{(0)} |_{\partial\Omega}$ are given; %
  \item ~ ~ ~ Calculate approximated linearized Neumann data $( \partial_{\nu} u_{S} - \partial_{\nu} u_{S}^{(0)} ) |_{\partial\Omega}$; 
  \item ~ ~ \textbf{End}; %
  \item ~ ~ Choose $\varphi := \exp\{ \bfi \zeta_{0}^{\langle r;\theta \rangle} \cdot x \}$ and $\psi := \exp\{ -\bfi \xi^{\langle r;\theta \rangle} \cdot x \}$; 
  \vspace{1ex}
  \item ~ ~ Compute $\mathcal{F}[c](\xi^{\langle r;\theta \rangle}) := \frac{1}{m!} \sum\limits_{\scriptscriptstyle \emptyset \subsetneqq S \subseteq U} {\scriptstyle (-1)^{| U \setminus S|}} \int\limits_{\scriptscriptstyle \partial\Omega} \beta^{-m} ( \partial_{\nu} u_{S} - \partial_{\nu} u_{S}^{(0)} ) \varphi \,\mathrm{d}S_{x}$; %
  \vspace{0.5ex}
  \item ~ ~ Update $c^{\langle r;\theta+1 \rangle} := c^{\langle r;\theta \rangle} + \mathcal{F}[c](\xi^{\langle r;\theta \rangle}) \, \psi \, \omega^{\langle r;\theta \rangle}$, for $\kappa_{r} \leq (m+1)k$; %
  \item ~ \textbf{End}; %
  \item ~ Set $c^{\langle r+1;1 \rangle} := c^{\langle r;\Theta+1 \rangle}$; %
  \item \textbf{End}. %
\end{enumerate} \\[-10pt]%
\bottomrule
\end{tabular}
\end{table}

The proposed \textbf{Algorithm 1} contains three iterations from Step 2 to Step 12 aiming to recover the (truncated) Fourier coefficients of the potential function. The iterations of $r$ and $\theta$ define the Fourier frequency $\xi^{\langle r;\theta \rangle} \in \mathbb{R}^{n}$ near which the Fourier coefficient $\mathcal{F}[c](\xi^{\langle r;\theta \rangle})$ will be reconstructed by the third iteration of $S \subseteq U$. We emphasize that the last iteration is based on the equality \eqref{eqn:identity_S} in Lemma \ref{lmm:identity_S}, see Step 10.

To better understand \textbf{Algorithm 1}, some remarks are listed as follows.
For the sake of simplicity, we explain these issues by ignoring the subscript notations $j$ or $S$ of $f_{j}$, $u_{j}$, $u_{S}$, $u_{S}^{(0)}$, $u_{S}^{(1)}$, and using $f$, $u$, $u^{(0)}$, $u^{(1)}$ for instead.



\begin{remark}\label{rmk:alg_scale}
The rescaling parameter $\beta > 0$ in Step 5 of \textbf{Algorithm 1}, comes from the well-posedness condition for the original problem \eqref{eqn:problem}. And thanks to the $m$-homogeneous property \eqref{eqn:m_homogeneous}, the corresponding linearized DtN map can be computed via rescaling in Step 10 of \textbf{Algorithm 1}.
By Proposition \ref{prp:wellposed}, the well-posedness of this nonlinear problem \eqref{eqn:problem} is guaranteed only when the Dirichlet boundary data $f$ is relatively small.
This feature of nonlinearity is a major difference compared to the linear problem \cite{ILX2020}.
However, Proposition \ref{prp:wellposed} is a qualitative result derived from the implicit function theorem for Banach spaces \cite[Theorem 10.6 and Remark 10.5]{RR2004}, which guarantees the well-posedness of the original problem \eqref{eqn:problem} only if the potential function $c$ and the Dirichlet boundary data $f$ satisfy the smallness assumptions, i.e. $\|c\|_{L^{\infty}(\Omega)} < \eta_{c}$ and $f \in B_{\eta}$, a small neighbourhood of $0$.
A concrete criteria for the smallness of the boundary data is not provided.
These fundamental issues are far beyond the scope of inverse problems in this paper.
Numerically, we could choose a small scale value $\beta \ll 1$, and regard the well-posedness holding for the original problem, that is $\beta f \in B_{\eta}$ sufficiently.
\end{remark}

\begin{remark}\label{rmk:alg_approx}
In a realistic scenario, Step 7 of \textbf{Algorithm 1}, one can measure the difference between Neumann boundary values of the original problem \eqref{eqn:problem} and the unperturbed problem \eqref{eqn:problem_1}, but the linearized Neumann boundary data of the linearized problem \eqref{eqn:problem_1} is not available.
Therefore, we do not generate the linearized Neumann boundary data directly, since it depends on the unknown potential function $c$ referring to the linearized problem \eqref{eqn:problem_1}.
Fortunately, in view of Proposition \ref{prp:linearized}, when the well-posedness of the original problem \eqref{eqn:problem} is guaranteed, the linearized Neumann boundary data can be approximated well by
\begin{equation*}
\partial_{\nu} u^{(1)} \big|_{\partial\Omega}
\approx \big( \partial_{\nu} u - \partial_{\nu} u^{(0)} \big) \big|_{\partial\Omega}
\end{equation*}
according to \eqref{eqn:linearized_err}. Recall that
\begin{equation}\label{eqn:alg_approx}
\begin{aligned}
&\big\| \partial_{\nu} u^{(1)} \big|_{\partial\Omega} - \big( \partial_{\nu} u - \partial_{\nu} u^{(0)} \big) \big|_{\partial\Omega} \big\|_{W^{1-\frac{1}{p},p}(\partial\Omega)} \\
&\qquad = \big\| \Lambda'_{c}(f) - \big( \Lambda_{c}(f) - \Lambda_{0}(f) \big) \big\|_{W^{1-\frac{1}{p},p}(\partial\Omega)}
\sim O(\|c\|_{L^{\infty}(\Omega)}^{2})
\end{aligned}
\end{equation}
for any $f \in B_{\eta}$; see Proposition \ref{prp:linearized} and its proof.
The similar approximations were mentioned in \cite[Equation (4.3)]{ILX2020} and \cite[Equation (4.1)]{LSX2022}.
\end{remark}

\section{The numerical examples}\label{sec:numerical}

In current section, we provide several numerical examples in 2D setting, i.e. the dimension $n = 2$.

\begin{figure}[!htb]
\centering
\,\hfill \textbf{(i)} \hfill \textbf{(ii)} \hfill\,\\
\includegraphics[width=0.28\textwidth,trim=100 20 110 0,clip]{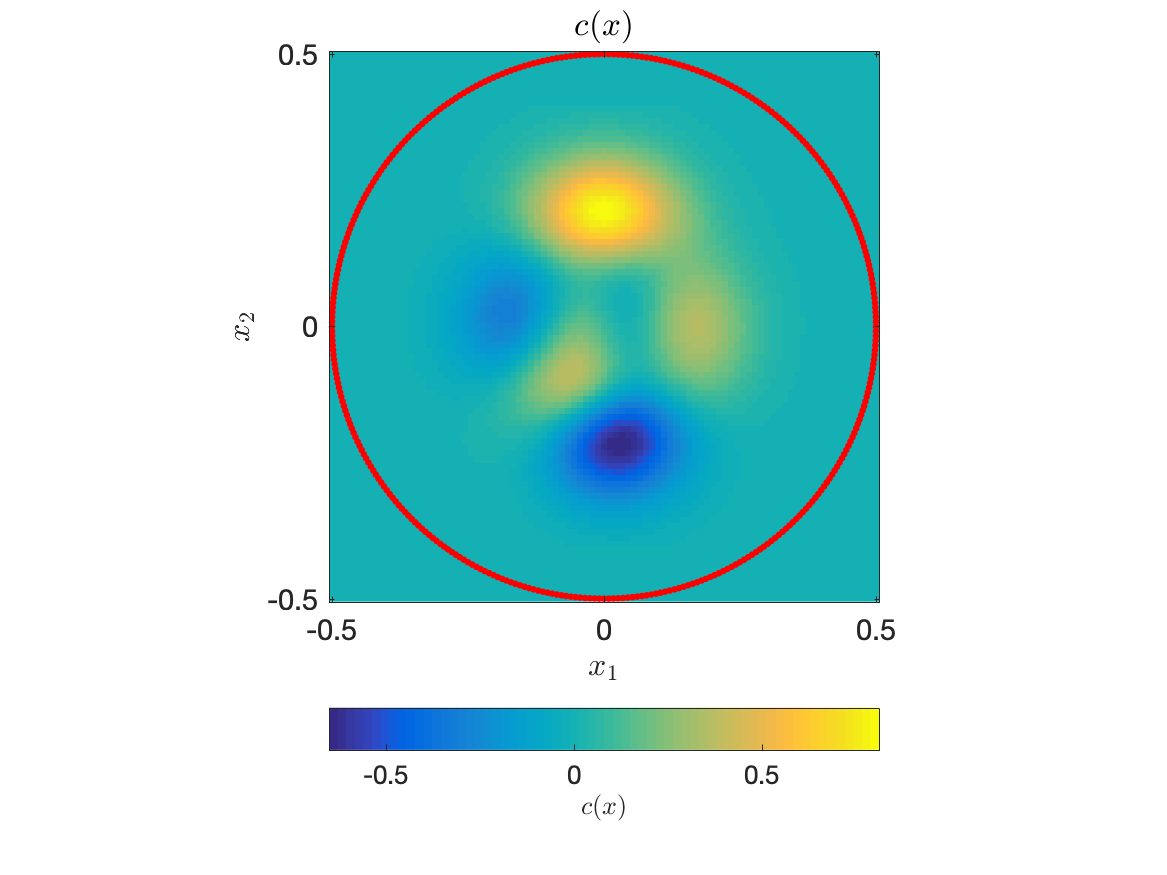}
\includegraphics[width=0.32\textwidth,trim=60 0 80 0,clip]{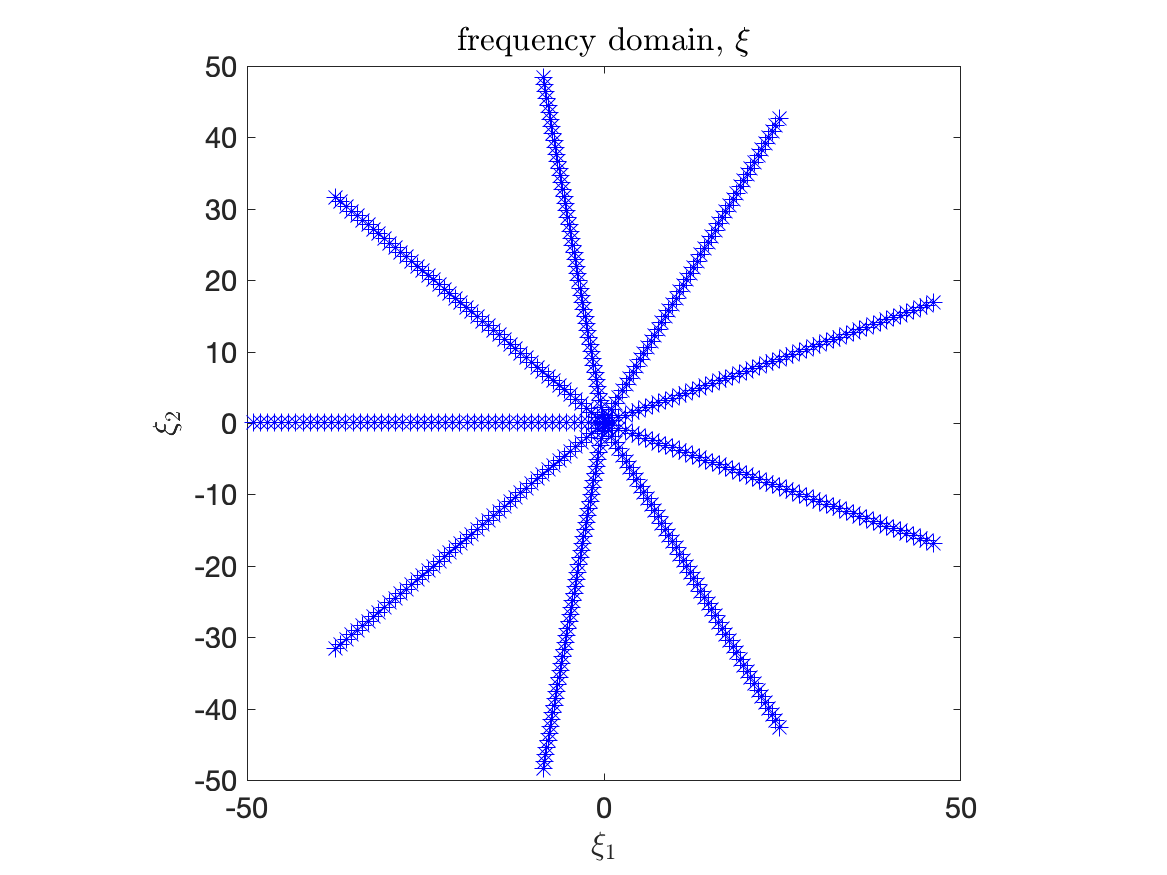}\\[0ex]
\,\hfill \textbf{(iii)} \hfill\,\\
\includegraphics[width=0.60\textwidth,trim=10 90 30 90,clip]{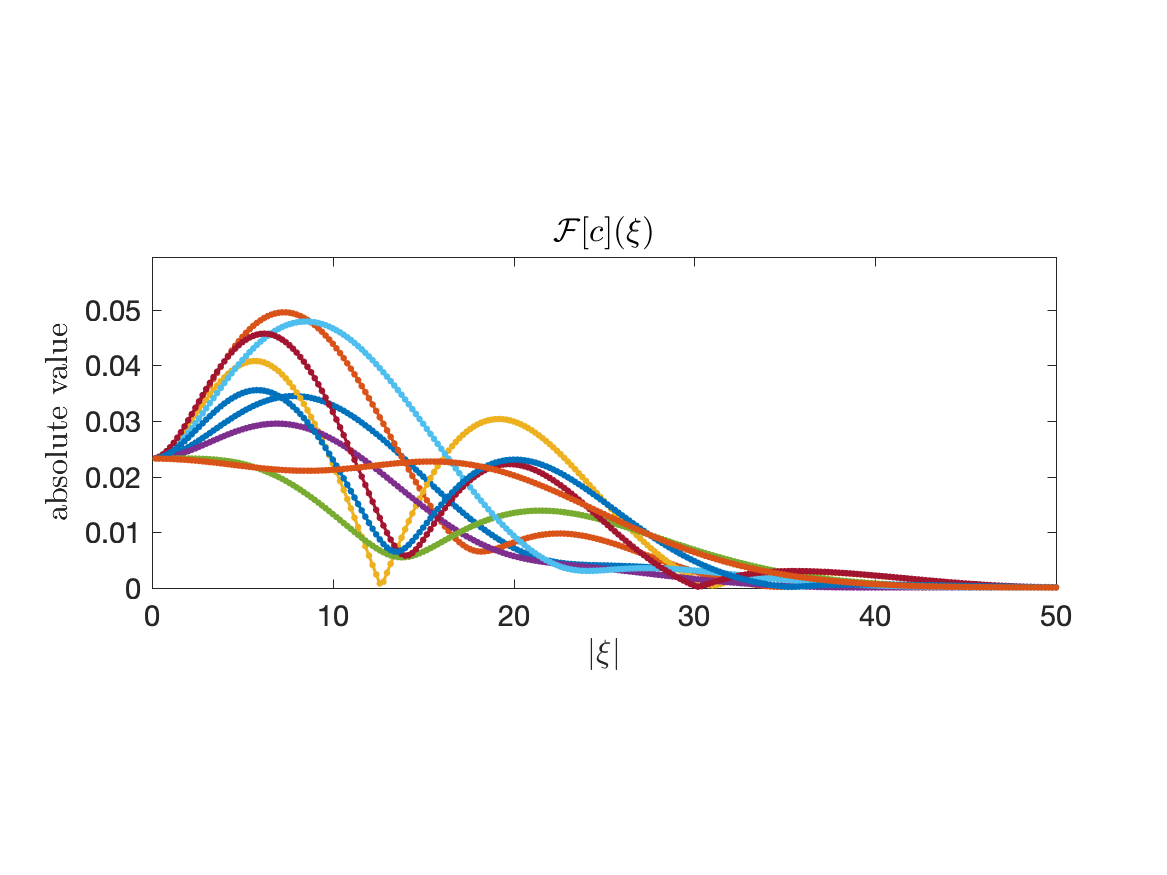}\\
\caption{\textbf{(i)} The exact potential function $c(x)$ in the domain $\Omega = B_{0.5}(0) \subset [{-0.5},0.5]^{2}$.
\textbf{(ii)} The sampling points $\xi$ with $|\xi| \leq 50$ in frequency domain.
\textbf{(iii)} The exact Fourier coefficients $\mathcal{F}[c](\xi)$ of $c$.
The horizontal axis shows the length $|\xi|$ of all $\xi$, and the vertical axis shows the absolute value $|\mathcal{F}[c](\xi)|$ of Fourier coefficients near the sampling points.}
\label{fig:potential_xi_Fc}
\end{figure}

The considered domain is a disk of radius $0.5$ centered at the origin s.t. $\Omega = B_{0.5}(0)$ and is contained by a square $[{-0.5},0.5]^{2}$; see Figure \ref{fig:potential_xi_Fc} \textbf{(i)}.
To avoid the inverse crime, we use a fine grid ($200 \times 200$ equal-distance points in the square) for the forward problem and a coarse grid ($90 \times 90$ equal-distance points in the square) for the inversion. The exact potential function $c := c(x)$ is also shown in Figure \ref{fig:potential_xi_Fc} \textbf{(i)}, $x = (x_{1},x_{2}) \in \mathbb{R}^{2}$, and the boundary $\partial\Omega$ is shown by a red circle.

The sampling points $\xi = (\xi_{1},\xi_{2}) \in \mathbb{R}^{2}$ with $|\xi| \leq 50$ in the frequency domain are shown in Figure \ref{fig:potential_xi_Fc} \textbf{(ii)}, marked by blue ``$\ast$'' near which all the Fourier coefficients $\mathcal{F}[c](\xi)$ will be recovered. In Figure \ref{fig:potential_xi_Fc} \textbf{(iii)}, the horizontal axis shows the length $|\xi|$ of all $\xi$, and the vertical axis shows the absolute value $|\mathcal{F}[c](\xi)|$ of Fourier coefficients of $c$ near the sampling points.

\subsection{Numerical example with $m = 4$}

\begin{figure}[!htb]
\centering
\textbf{Quartic nonlinearity:} $m = 4$. \\[2ex]
\,\hfill \textbf{(i)} $k = 5$ \hfill\,\\
\includegraphics[width=0.45\textwidth,trim=20 60 20 35,clip]{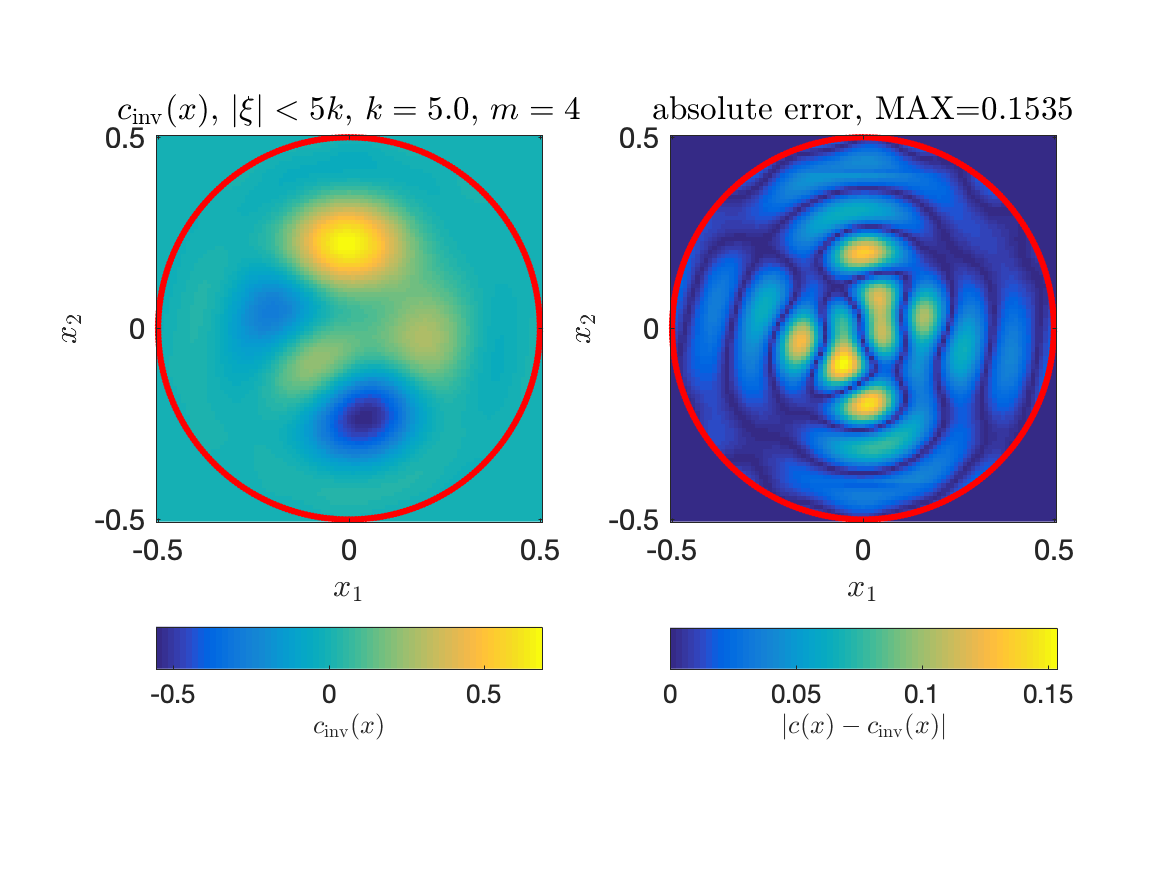}
\includegraphics[width=0.54\textwidth,trim=10 60 30 90,clip]{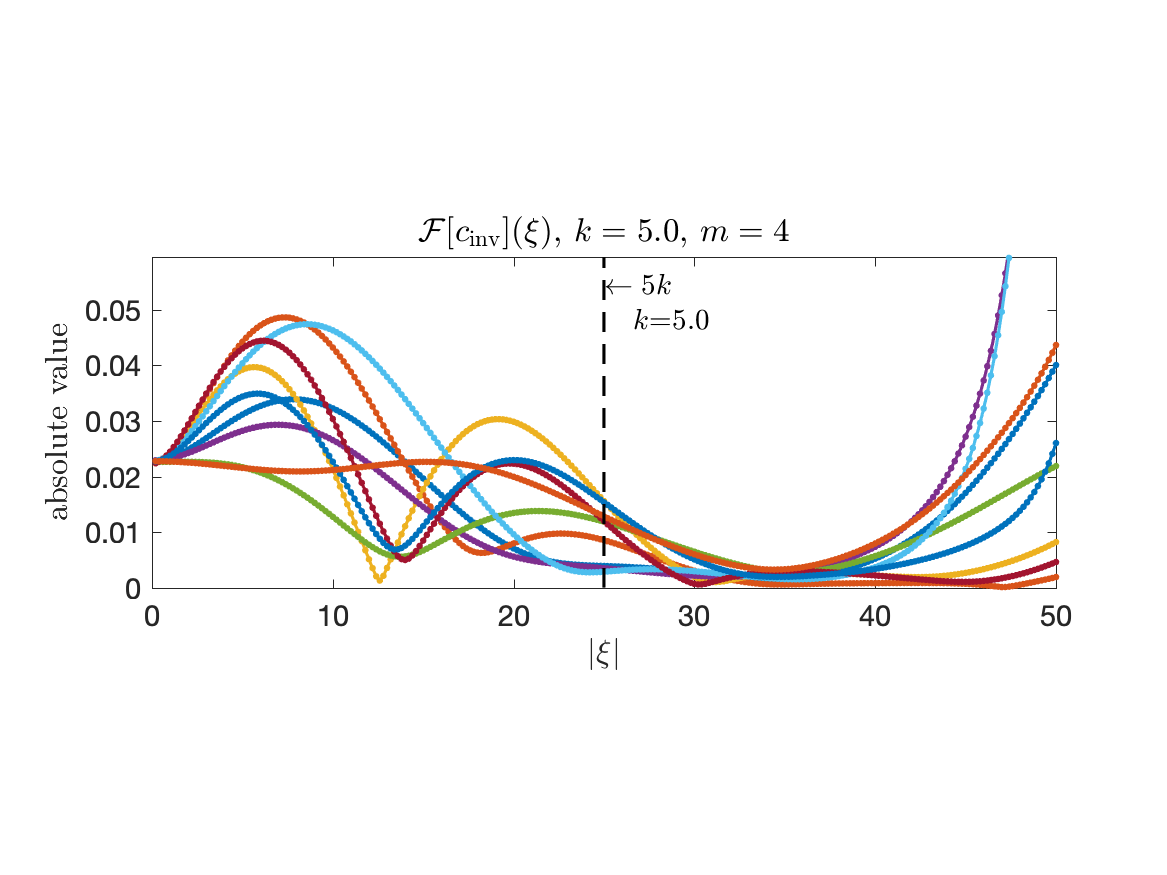}\\[2ex]
\,\hfill \textbf{(ii)} $k = 10$ \hfill\,\\
\includegraphics[width=0.45\textwidth,trim=20 60 20 35,clip]{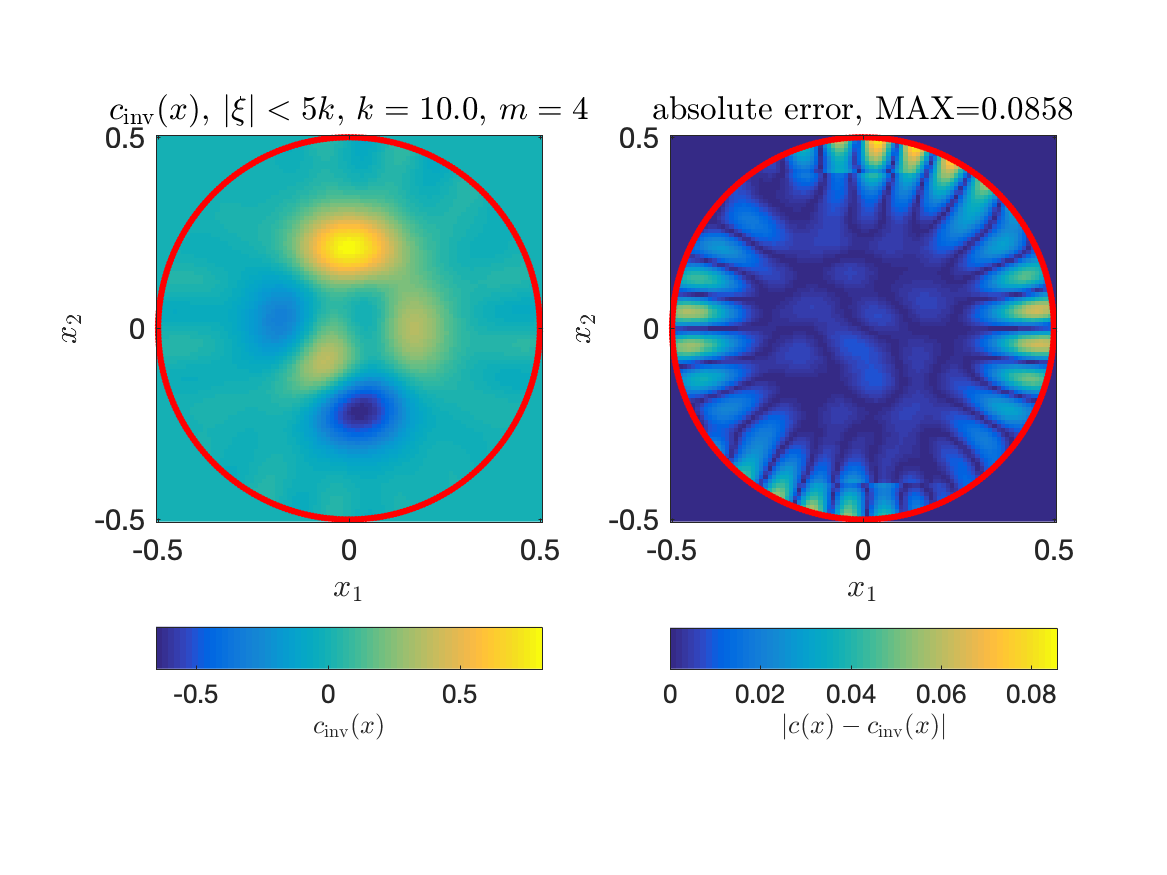}
\includegraphics[width=0.54\textwidth,trim=10 60 30 90,clip]{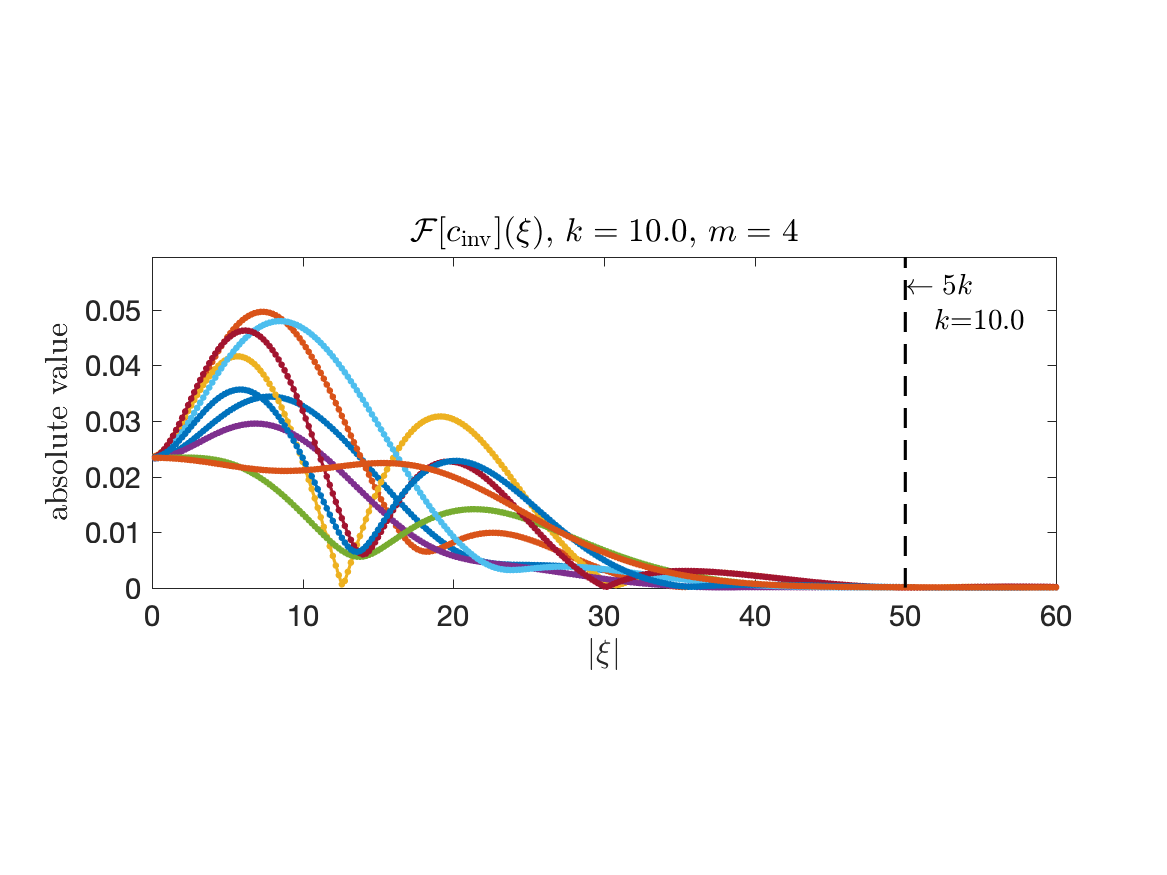}\\[2ex]
\,\hfill \textbf{(iii)} $k = 15$ \hfill\,\\
\includegraphics[width=0.45\textwidth,trim=20 60 20 35,clip]{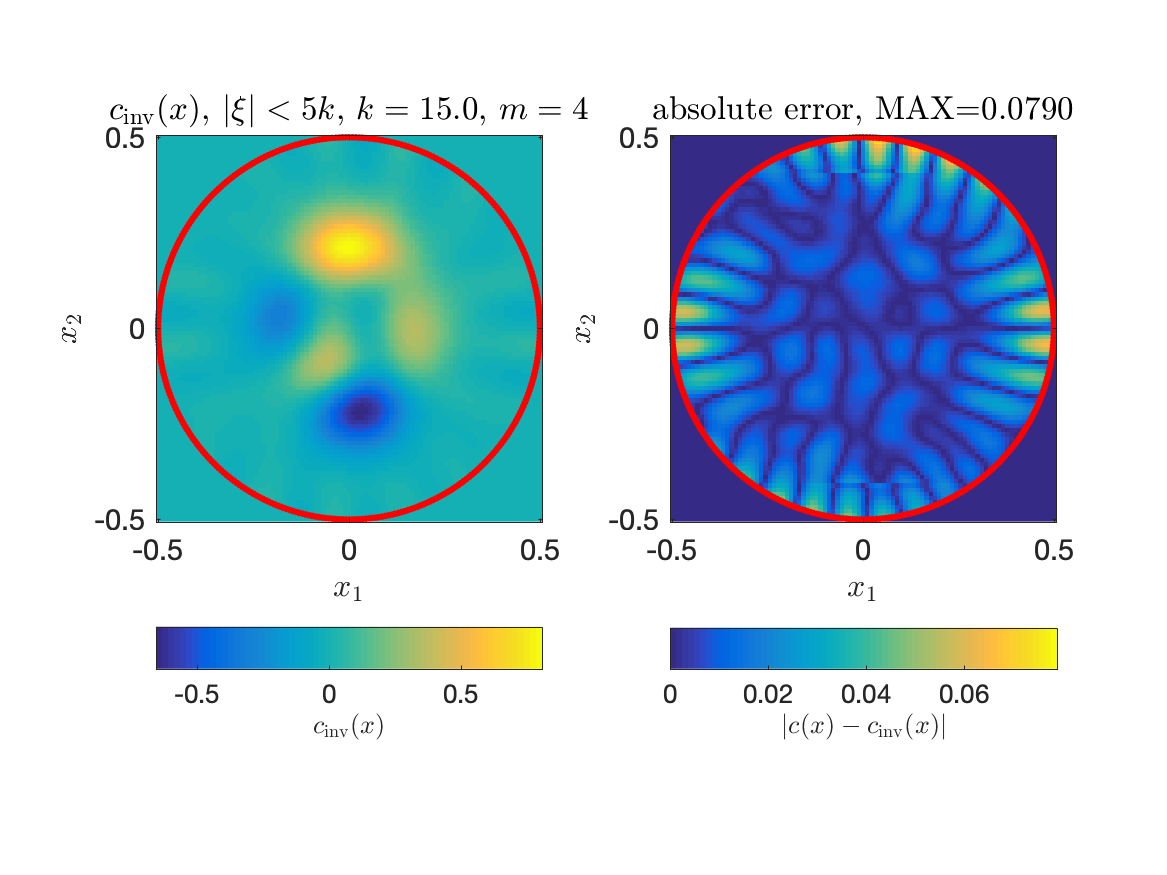}
\includegraphics[width=0.54\textwidth,trim=10 60 30 90,clip]{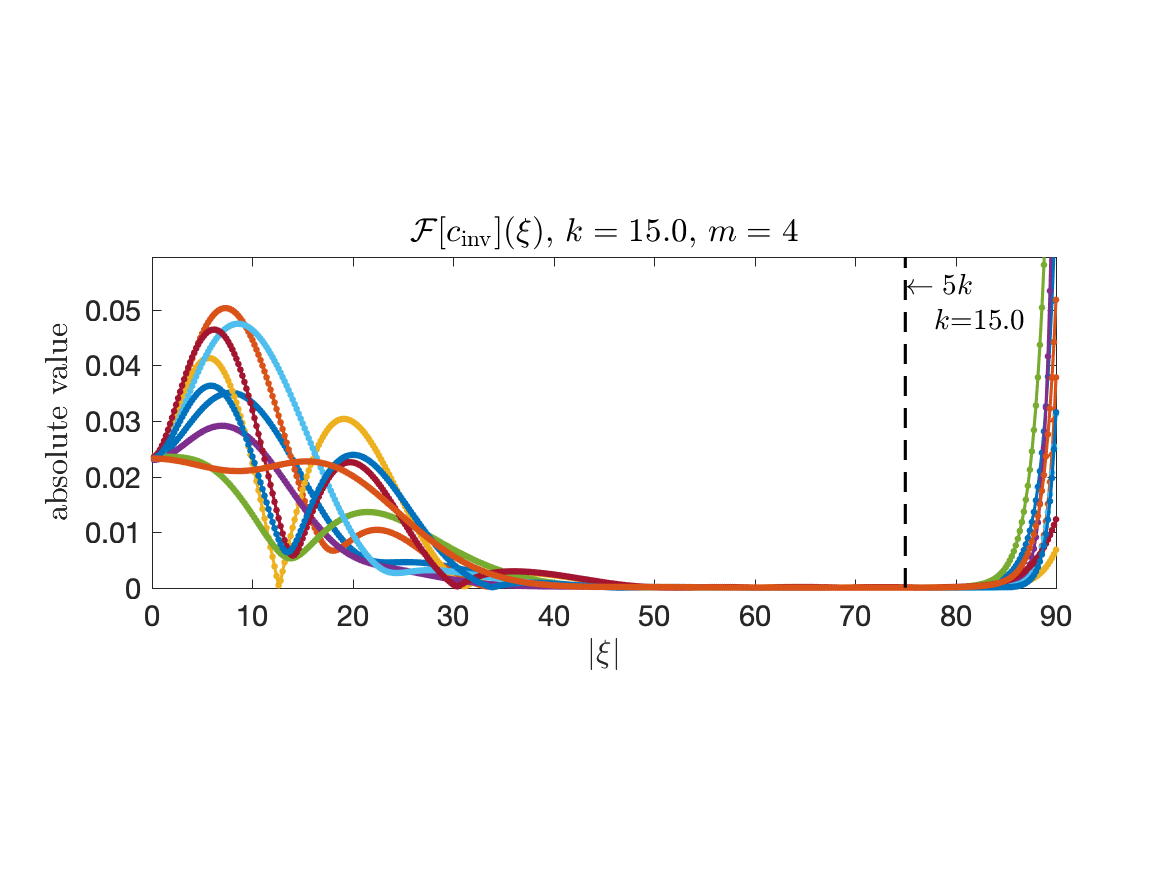}\\[2ex]
\caption{[The nonlinearity index $m = 4$.]
The reconstruction results, when the wavenumber \textbf{(i)} $k = 5$, \textbf{(ii)} $k = 10$ and \textbf{(iii)} $k = 15$.
In each sub-figure,
\textbf{Left column:} the recovered potential $c_{\textrm{inv}}(x)$ with $|\xi| \leq (m+1)k$.
\textbf{Middle column:} the absolute error between the exact and recovered potential functions $|c(x) - c_{\textrm{inv}}(x)|$ in $\Omega$.
\textbf{Right column:} the absolute value of the recovered Fourier coefficients $|\mathcal{F}[c_{\textrm{inv}}](\xi)|$.}
\label{fig:4_Fc}
\end{figure}

To numerically test \textbf{Algorithm 1}, let the nonlinearity index $m = 4$, we implement \textbf{Algorithm 1} to reconstruct the potential function $c(x)$, including all the recovered Fourier coefficients $\mathcal{F}[c_{\textrm{inv}}](\xi)$ with $|\xi| \leq (m+1)k$.
In \textbf{left column} of Figure \ref{fig:4_Fc}, we present the reconstructed potential functions $c_{\textrm{inv}}(x)$ with different wavenumbers: \textbf{(i)} $k = 5$, \textbf{(ii)} $k = 10$ and \textbf{(iii)} $k = 15$, respectively.
In \textbf{middle column} of Figure \ref{fig:4_Fc}, the absolute errors $|c(x) - c_{\textrm{inv}}(x)|$ between the exact and recovered potential functions are also shown, one can see that \textbf{Algorithm 1} reduces the maximum absolute error from $0.1535$ to $0.0790$ when $k$ increases from $5$ to $15$.
Moreover, the recovered Fourier coefficients $\mathcal{F}[c_{\textrm{inv}}](\xi)$ are shown in \textbf{right column} of Figure \ref{fig:4_Fc}.
It concludes that, while $k$ is larger, the more Fourier modes can be recovered stably, i.e. $\mathcal{F}[c_{\textrm{inv}}](\xi)$ with $|\xi| \leq (m+1)k$. These results numerically verify the increasing stability in Theorem \ref{thm:lipschitz} and Theorem \ref{thm:stability} while $k$ becomes large.


\subsection{Extended numerical tests: different index $m$ and non-smooth case}

\begin{figure}[!htb]
\centering
\textbf{The wavenumber:} $k = 5$. \\[2.5ex]
\,\hfill \textbf{(i)} $m = 3$ \hfill\,\\
\includegraphics[width=0.45\textwidth,trim=20 60 20 35,clip]{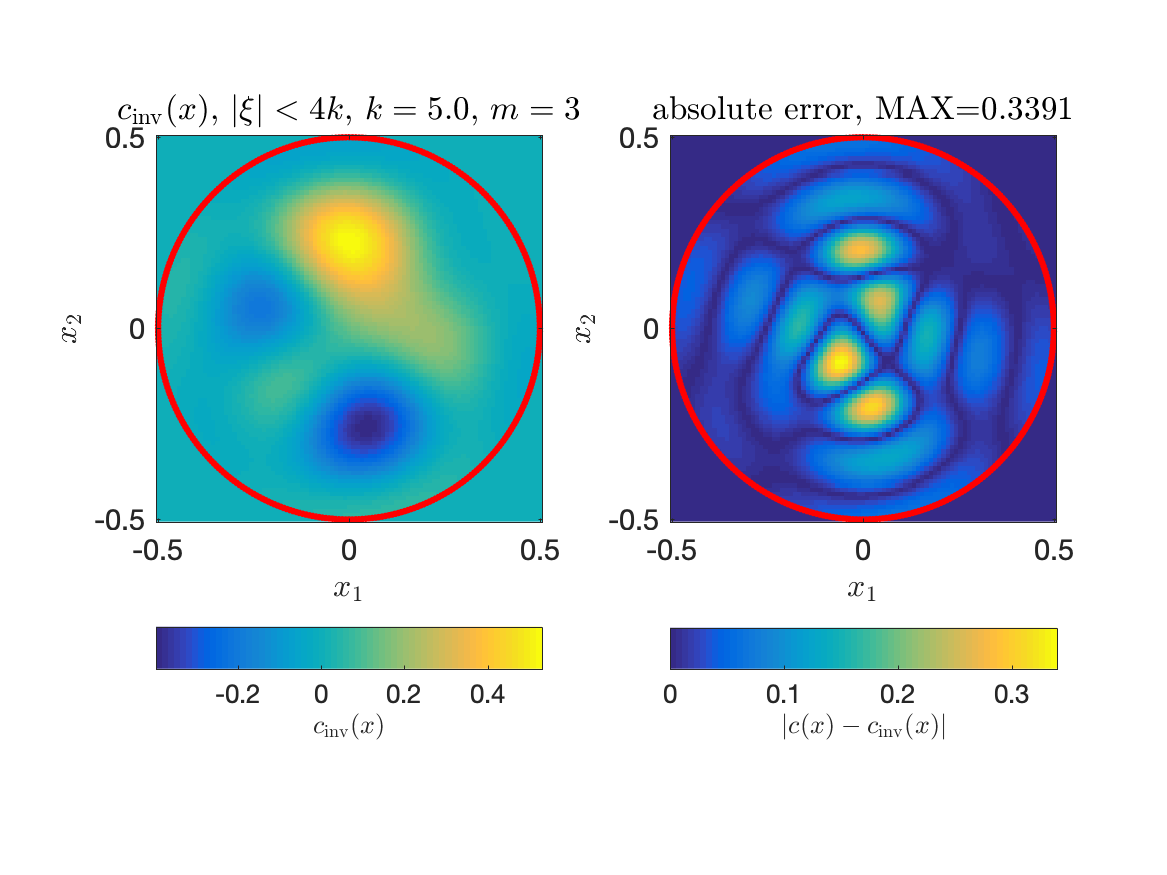}
\includegraphics[width=0.54\textwidth,trim=10 60 30 90,clip]{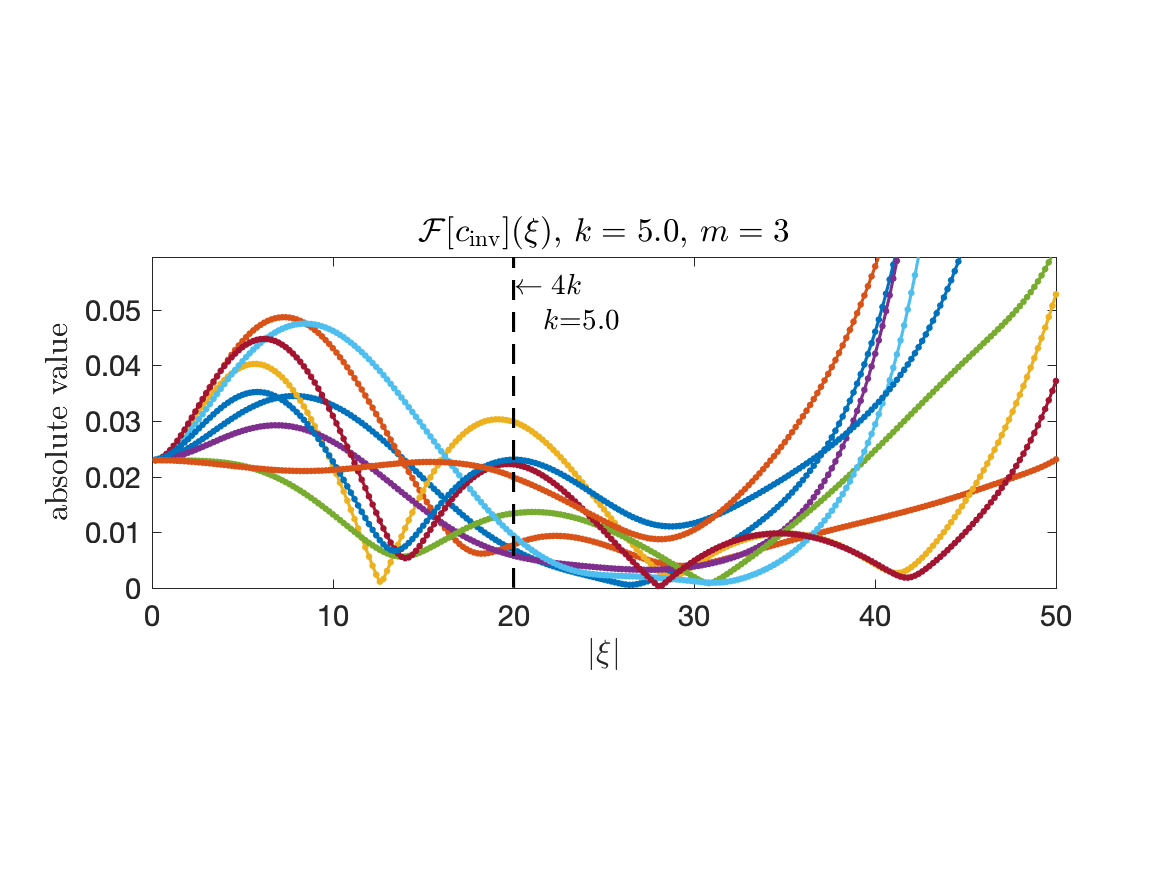}\\[2.5ex]
\,\hfill \textbf{(ii)} $m = 4$ \hfill\,\\
\includegraphics[width=0.45\textwidth,trim=20 60 20 35,clip]{./figs/4_Ic_0050.png}
\includegraphics[width=0.54\textwidth,trim=10 60 30 90,clip]{./figs/4_Fc_0050.png}\\[2.5ex]
\,\hfill \textbf{(iii)} $m = 5$ \hfill\,\\
\includegraphics[width=0.45\textwidth,trim=20 60 20 35,clip]{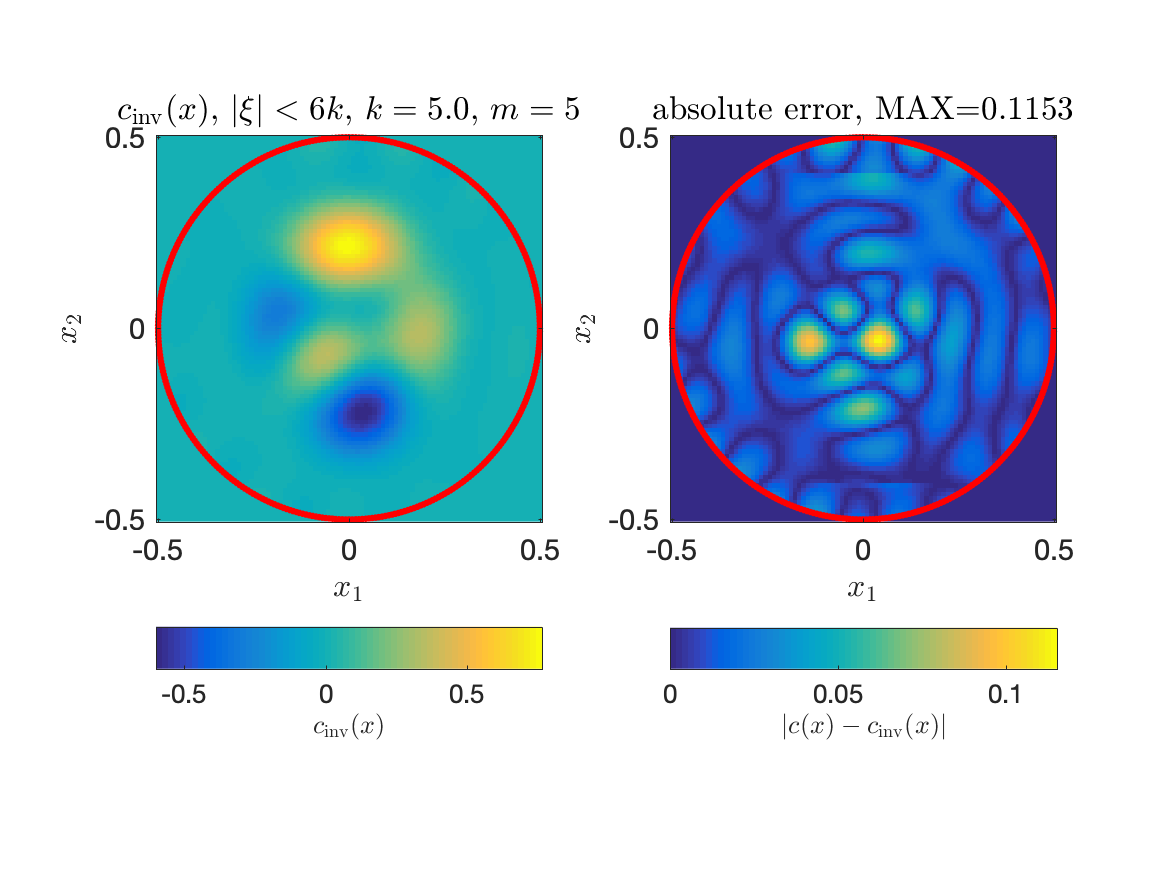}
\includegraphics[width=0.54\textwidth,trim=10 60 30 90,clip]{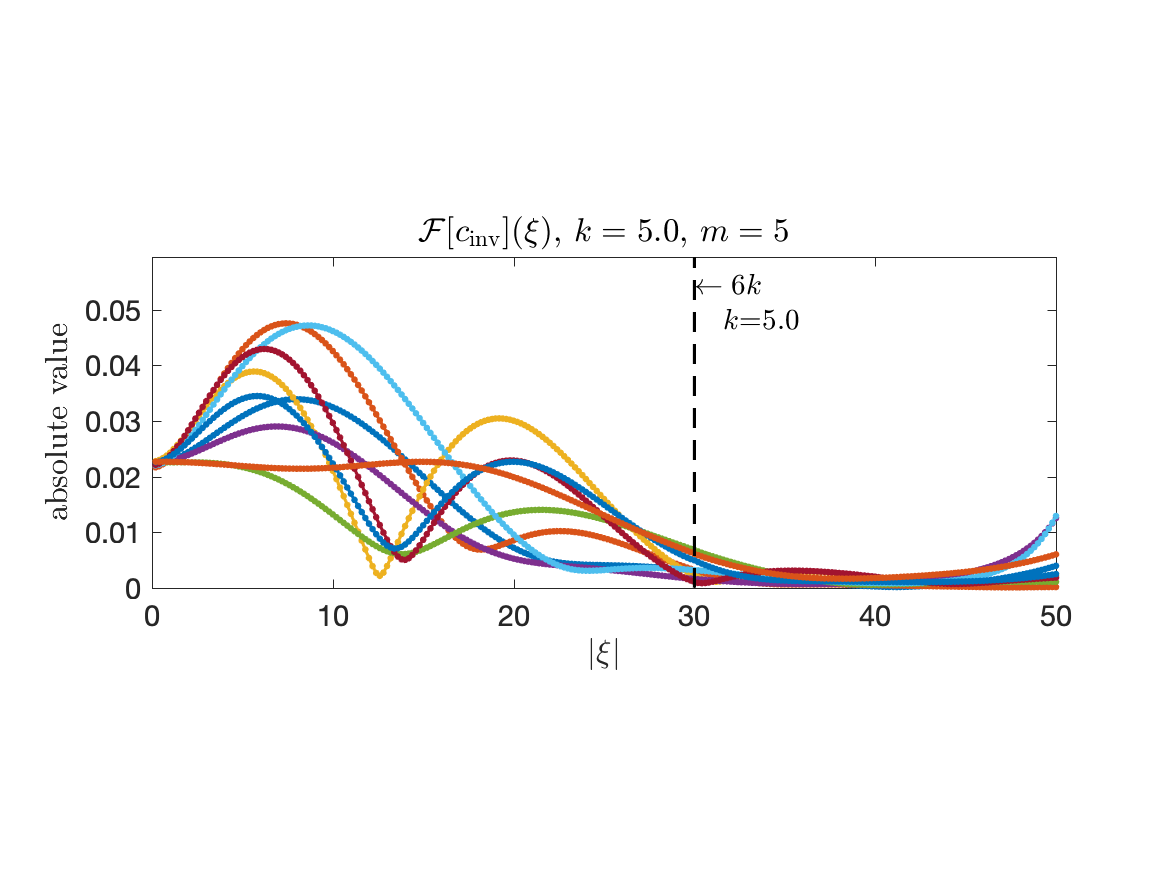}\\[2ex]
\caption{[The wavenumber $k = 5$.]
The reconstructed results, when the nonlinearity index \textbf{(i)} $m = 3$, \textbf{(ii)} $m = 4$ and \textbf{(iii)} $m = 5$.
In each sub-figure,
\textbf{Left column:} the recovered potential $c_{\textrm{inv}}(x)$ with $|\xi| \leq (m+1)k$,
\textbf{Middle column:} the absolute error between the exact and recovered potential functions $|c(x) - c_{\textrm{inv}}(x)|$ in $\Omega$.
\textbf{Right column:} the absolute value of the recovered Fourier coefficients $|\mathcal{F}[c_{\textrm{inv}}](\xi)|$.}
\label{fig:diff_m}
\end{figure}

In this subsection, we provide two comparability tests with different choice of the index $m$ and a non-smooth potential function.
We first compare the results for different nonlinearity index $m$ while the wavenumber $k = 5$; see Figure \ref{fig:diff_m}. The reconstructed potential functions $c_{\textrm{inv}}(x)$ (\textbf{left column}) are shown together with their absolute errors $|c(x) - c_{\textrm{inv}}(x)|$ (\textbf{middle column}) and recovered Fourier coefficients $\mathcal{F}[c_{\textrm{inv}}](\xi)$ (\textbf{right column}) in Figure \ref{fig:diff_m}, when the nonlinearity index is chosen to be \textbf{(i)} $m = 3$, \textbf{(ii)} $m = 4$ and \textbf{(iii)} $m = 5$.
Obviously, these results also verify the increasing stability, since the maximum absolute error is reduced from $0.5732$ to $0.1153$ and the range of Lipschitz type stability $|\xi| \leq (m+1)k$ becomes larger while $m$ increases.

\begin{figure}[!htb]
\centering
\,\hfill\; \textbf{(i)} \quad\hfill\hfill\; \textbf{(ii)} \hfill\hfill\,\\
\includegraphics[width=0.24\textwidth,trim=100 20 110 0,clip]{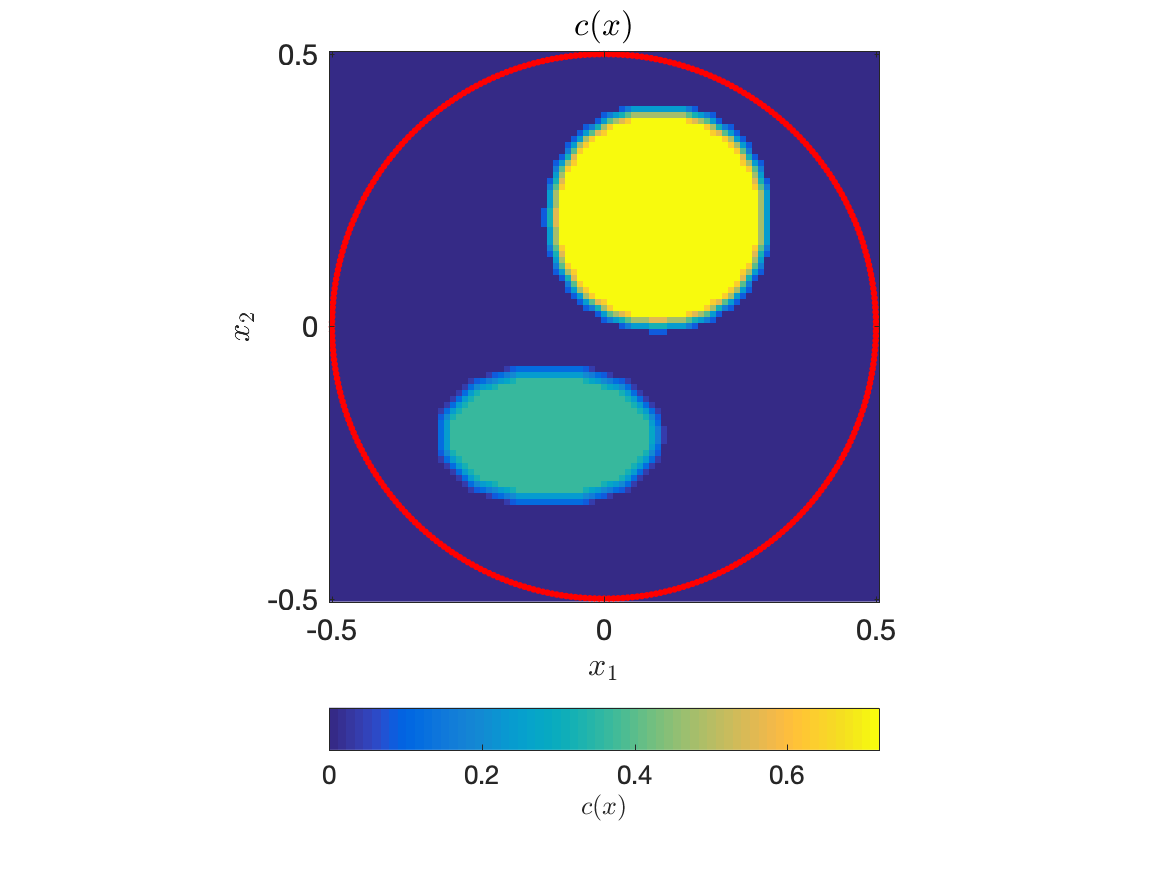}
\includegraphics[width=0.54\textwidth,trim=10 60 30 90,clip]{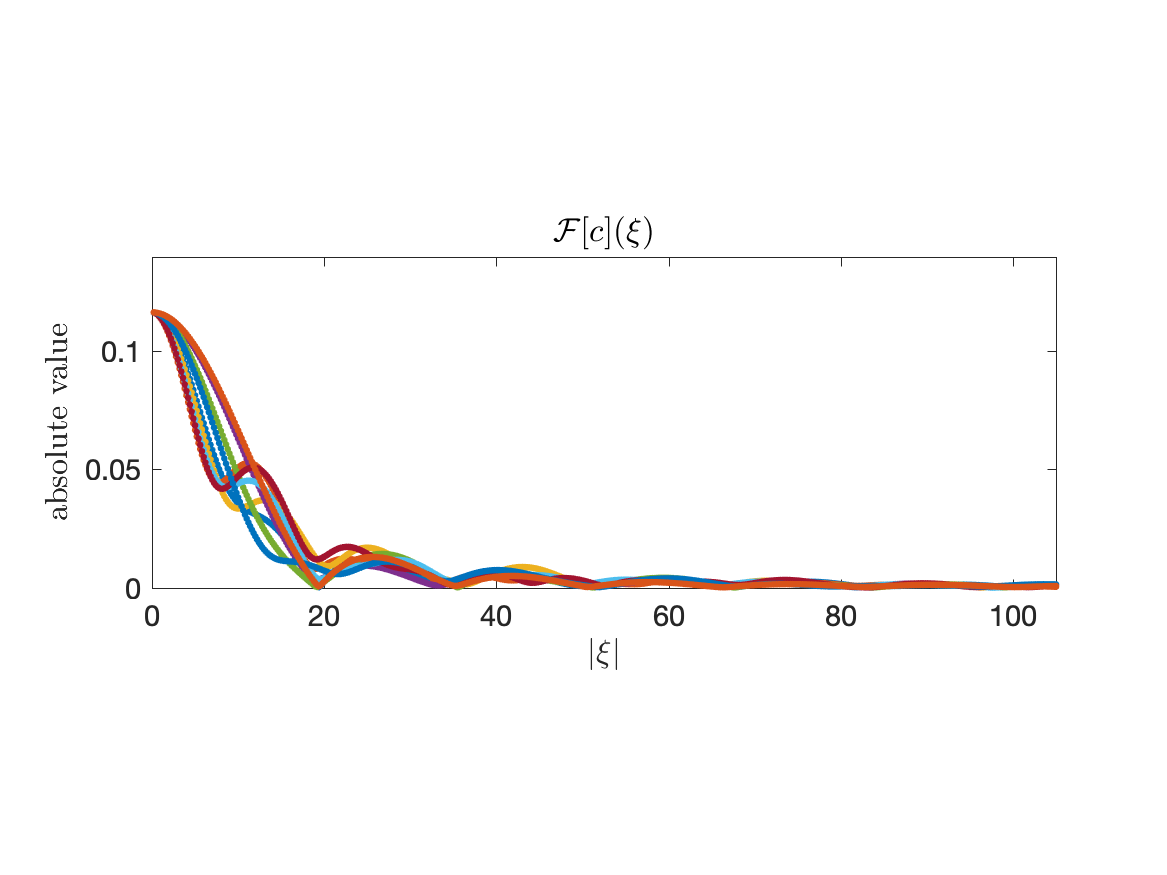}\\
\caption{\textbf{(i)} The exact potential function $c(x)$ in the domain $\Omega = B_{0.5}(0) \subset [{-0.5},0.5]^{2}$.
\textbf{(ii)} The exact Fourier coefficients $\mathcal{F}[c](\xi)$ of $c$.}
\label{fig:pwc_potential}
\end{figure}

\begin{figure}[!htb]
\centering
\textbf{Quintic nonlinearity:} $m = 5$ and
\textbf{the wavenumber:} $k = 15$. \\[1ex]
\,\hfill \textbf{(i)} \quad\hfill\; \textbf{(ii)} \hfill\hfill\hfill \textbf{(iii)} \hfill\hfill\,\\
\includegraphics[width=0.45\textwidth,trim=20 60 20 35,clip]{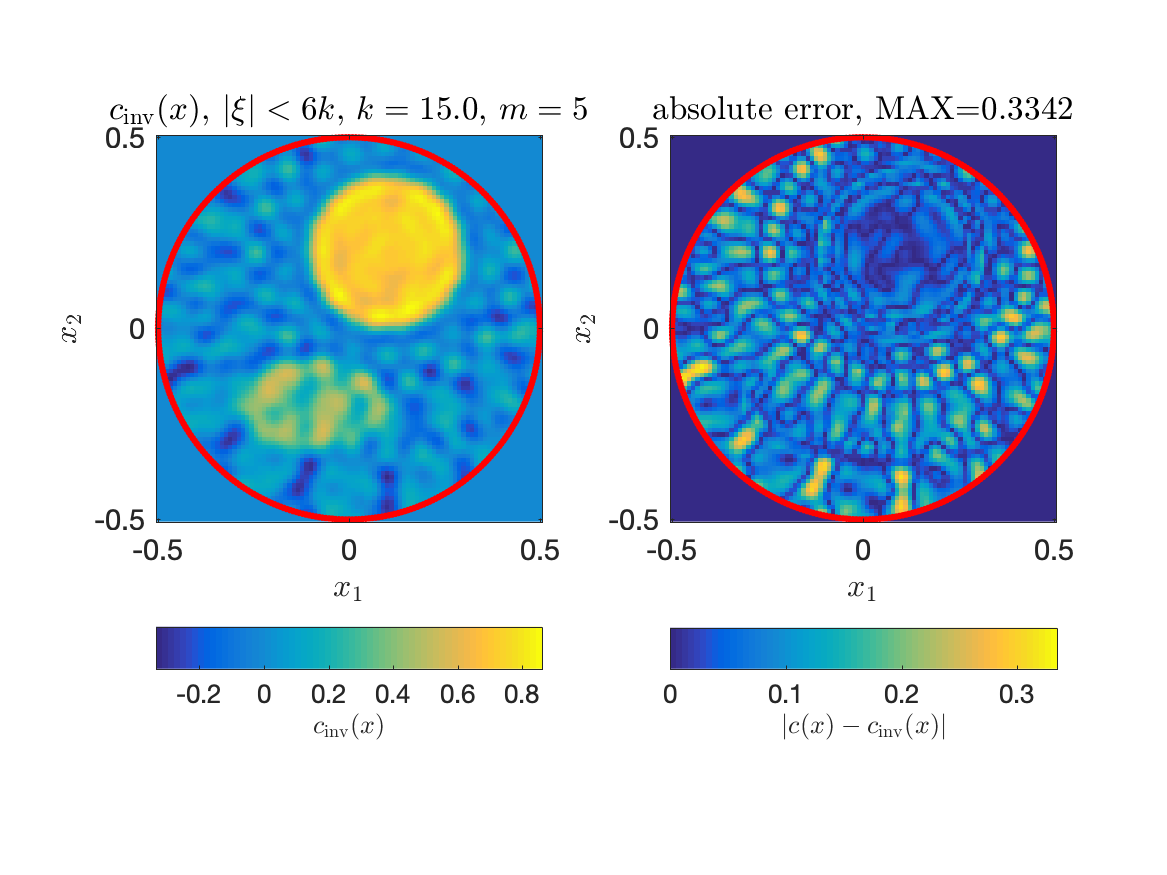}
\includegraphics[width=0.54\textwidth,trim=10 60 30 90,clip]{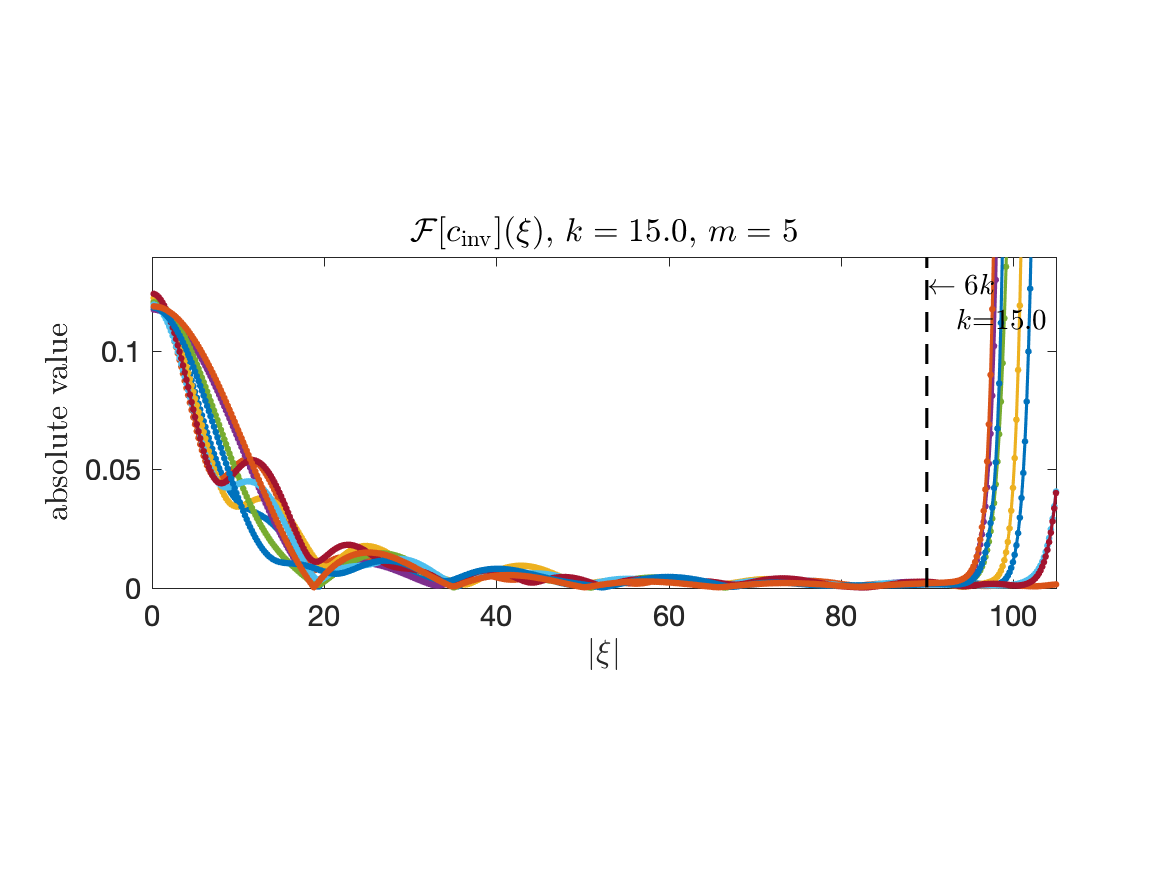}\\
\caption{[The nonlinearity index $m = 5$ and the wavenumber $k = 15$.]
\textbf{(i)} the recovered potential $c_{\textrm{inv}}(x)$ with $|\xi| \leq (m+1)k$,
\textbf{(ii)} the absolute error between the exact and recovered potential functions $|c(x) - c_{\textrm{inv}}(x)|$ in $\Omega$, and
\textbf{(iii)} the absolute value of the recovered Fourier coefficients $|\mathcal{F}[c_{\textrm{inv}}](\xi)|$.}
\label{fig:5_pwc}
\end{figure}

Although the assumption that $c \in H^{1}(\Omega)$ is needed in Theorem \ref{thm:stability}, we only recover the unknown potential function $c$ within the Lipschitz stable region $|\xi| \leq (m+1)k$; see Theorem \ref{thm:lipschitz}. Therefore, a non-smooth case is provided in this subsection to illustrate the applicability of the \textbf{Algorithm 1}.
The exact piecewise constant potential function $c(x)$ is shown in Figure \ref{fig:pwc_potential} \textbf{(i)}, the corresponding absolute value $|\mathcal{F}[c](\xi)|$ of Fourier coefficients of $c$ near the same sampling points as before are shown in Figure \ref{fig:pwc_potential} \textbf{(ii)}.
Then, we implement the \textbf{Algorithm 1} and present, in Figure \ref{fig:5_pwc} \textbf{(i)}, the reconstructed potential function $c_{\textrm{inv}}(x)$ using all the recovered Fourier coefficients $\mathcal{F}[c_{\textrm{inv}}](\xi)$ within $|\xi| \leq (m+1)k$ in Figure \ref{fig:5_pwc} \textbf{(iii)}.
The absolute error $|c(x) - c_{\textrm{inv}}(x)|$ between the exact and recovered potential functions is also provided in Figure \ref{fig:5_pwc} \textbf{(ii)}.


\subsection{Numerical explanation: combined solution $u_{S}^{(0)}$}\label{sec:combined}

In this subsection, we provide a numerical explanation for the combined solution $u_{S}^{(0)}$ and the combination in Alessandrini-PIE type identity \eqref{eqn:identity_S}.

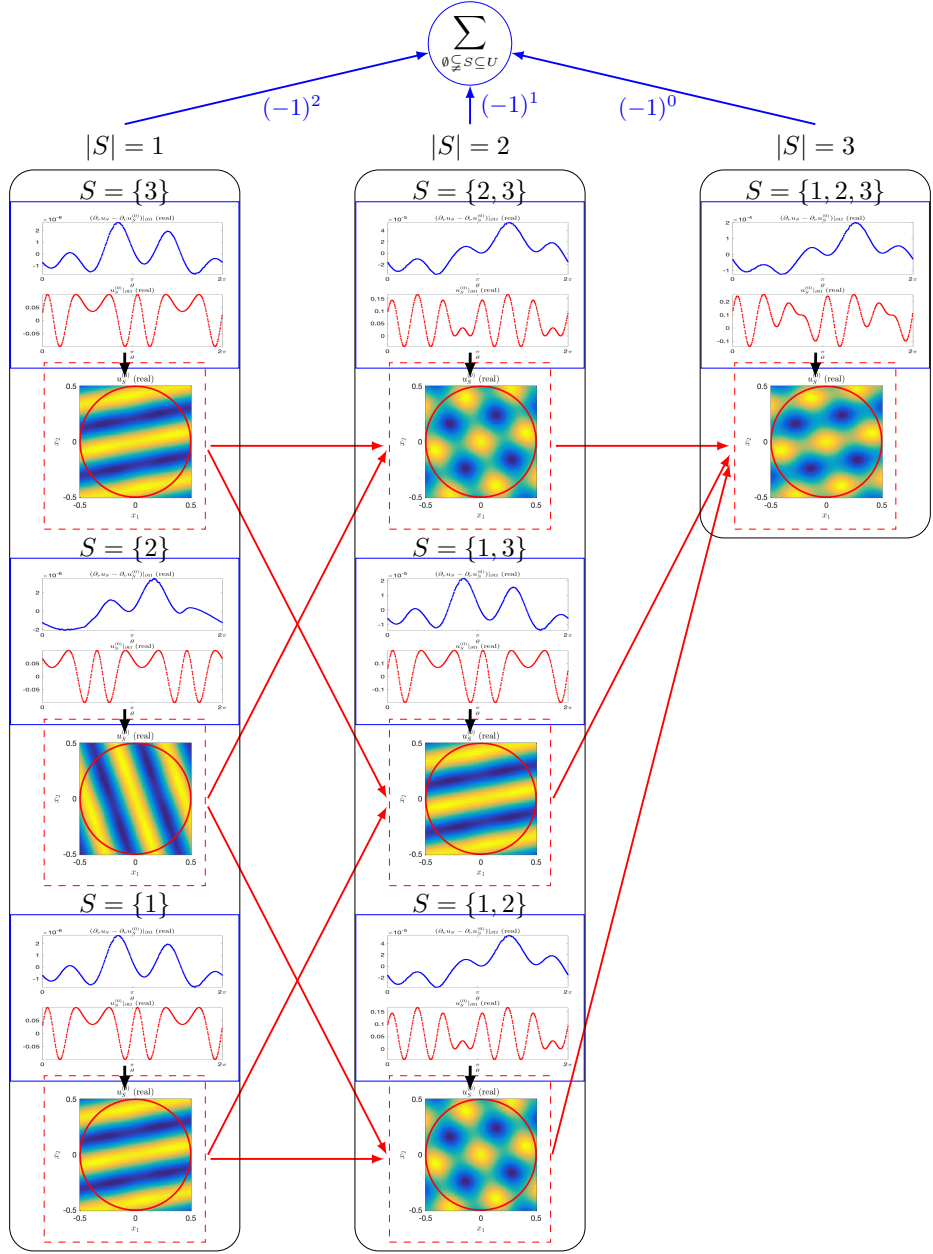
\begin{figure}[!htb]
\centering
\resizebox{0.9\width}{0.9\height}{
\begin{tikzpicture}[>=latex,scale=1.7]
\node[draw,blue]   (H1b) at (0,1.2) %
{\includegraphics[width=0.2\textwidth,height=0.1\textheight,trim=30 10 130 0,clip]{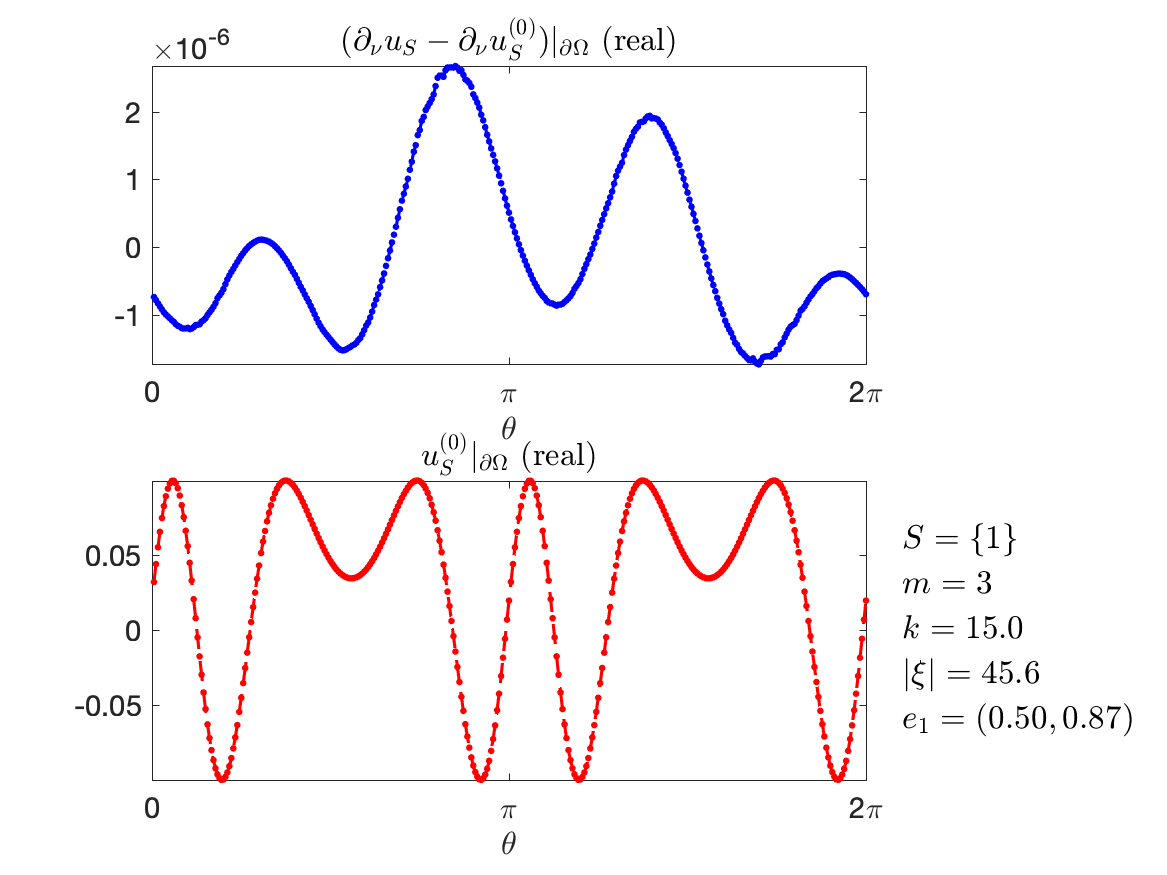}}
node at (0,2.0) {\small$S = \{1\}$};%
\node[draw,red ,dashed]   (H1)  at (0,-0.2)   %
{\includegraphics[height=0.1\textheight,trim=140 120 160 30,clip]{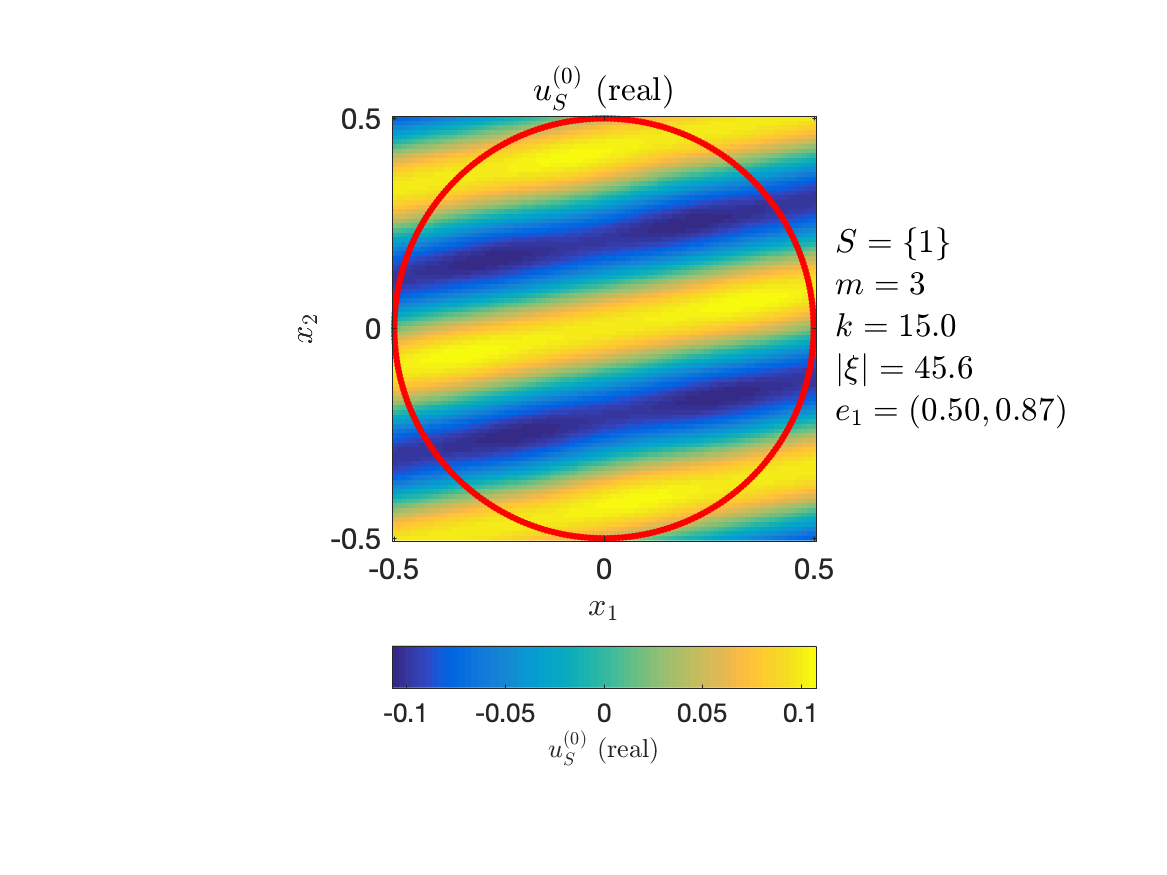}};
\node[draw,blue]   (H2b) at (0,4.3) %
{\includegraphics[width=0.2\textwidth,height=0.1\textheight,trim=30 10 130 0,clip]{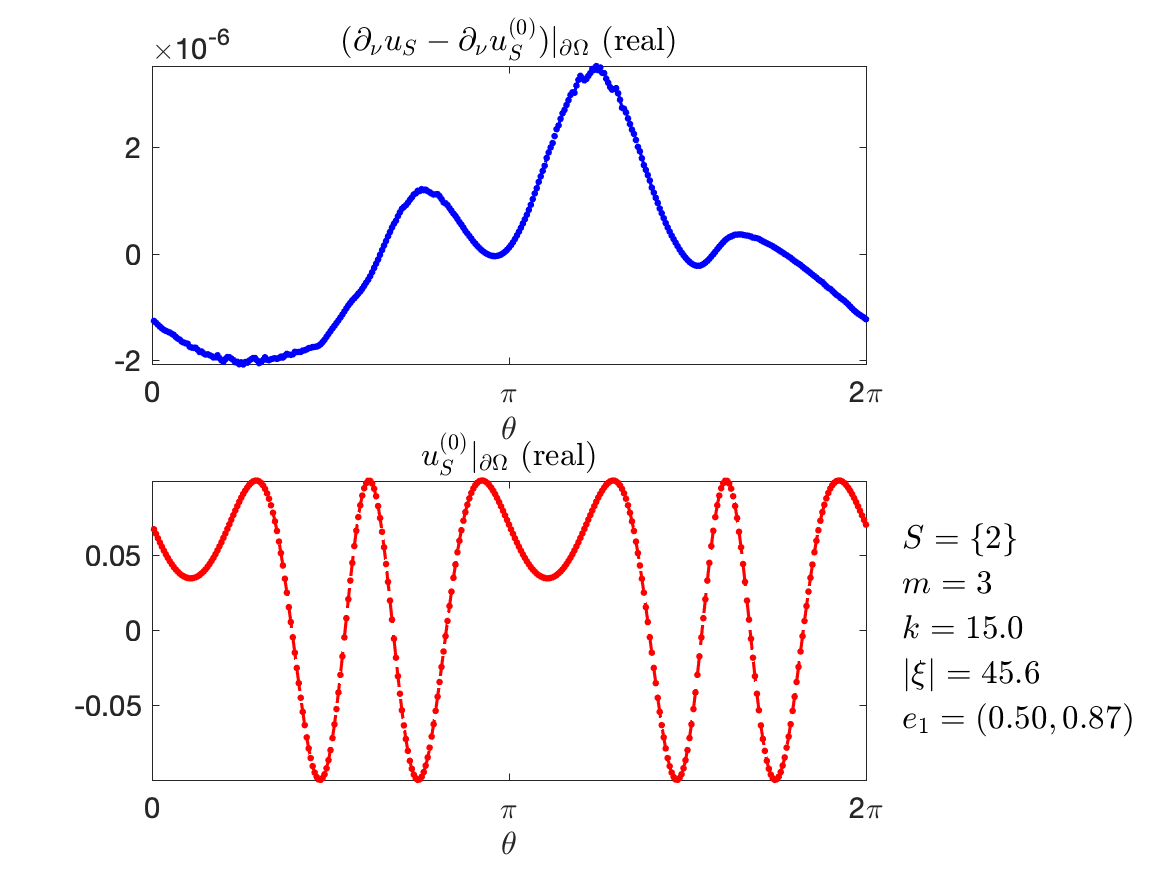}}
node at (0,5.1) {\small$S = \{2\}$};%
\node[draw,red ,dashed]   (H2)  at (0,2.9)   %
{\includegraphics[height=0.1\textheight,trim=140 120 160 30,clip]{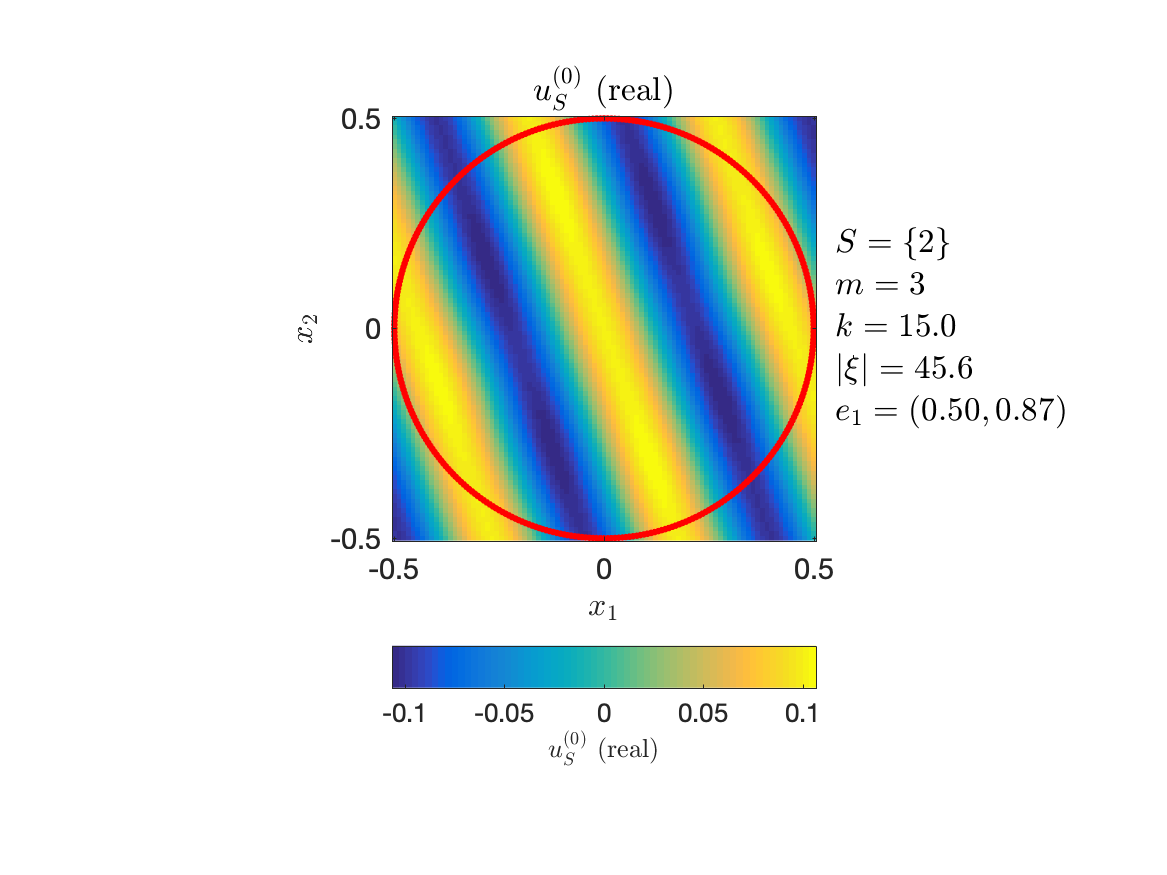}};
\node[draw,blue]   (H3b) at (0,7.4) %
{\includegraphics[width=0.2\textwidth,height=0.1\textheight,trim=30 10 130 0,clip]{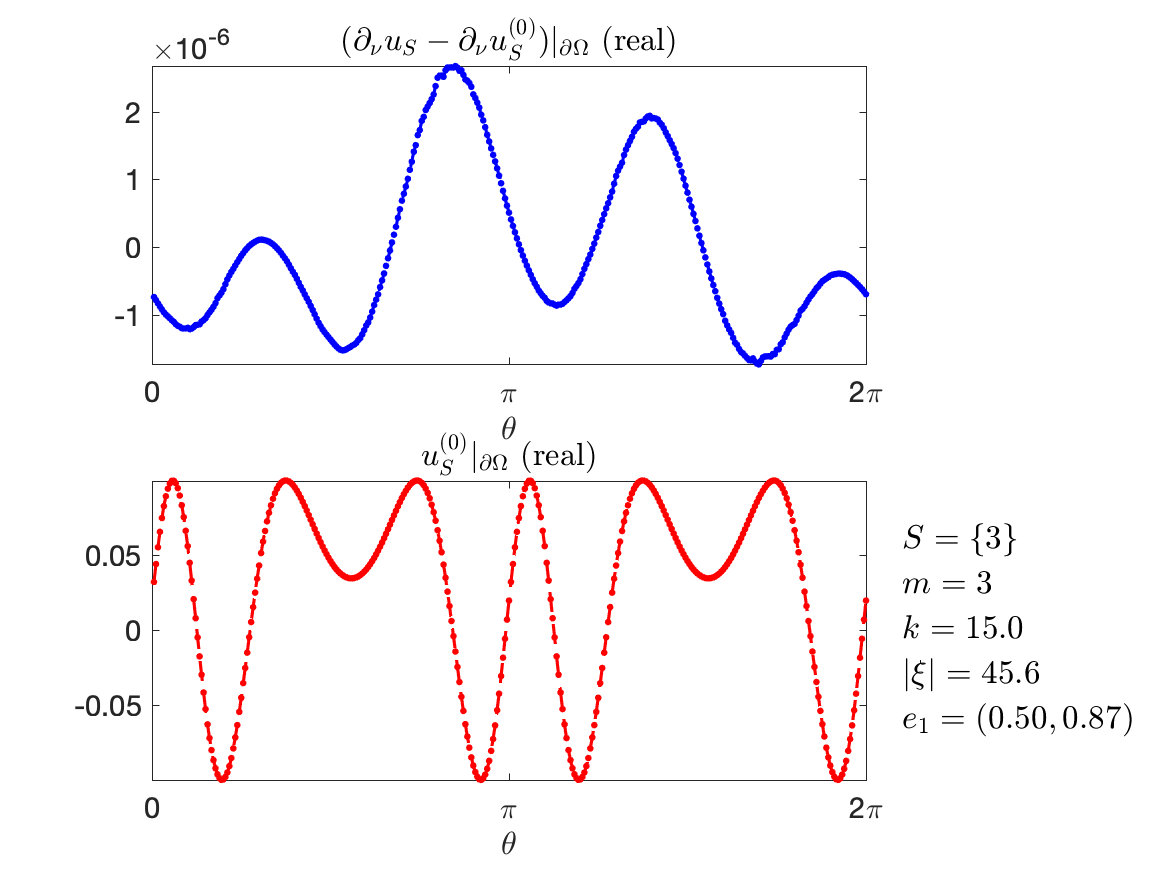}}
node at (0,8.2) {\small$S = \{3\}$};%
\node[draw,red ,dashed]   (H3)  at (0,6.0)   %
{\includegraphics[height=0.1\textheight,trim=140 120 160 30,clip]{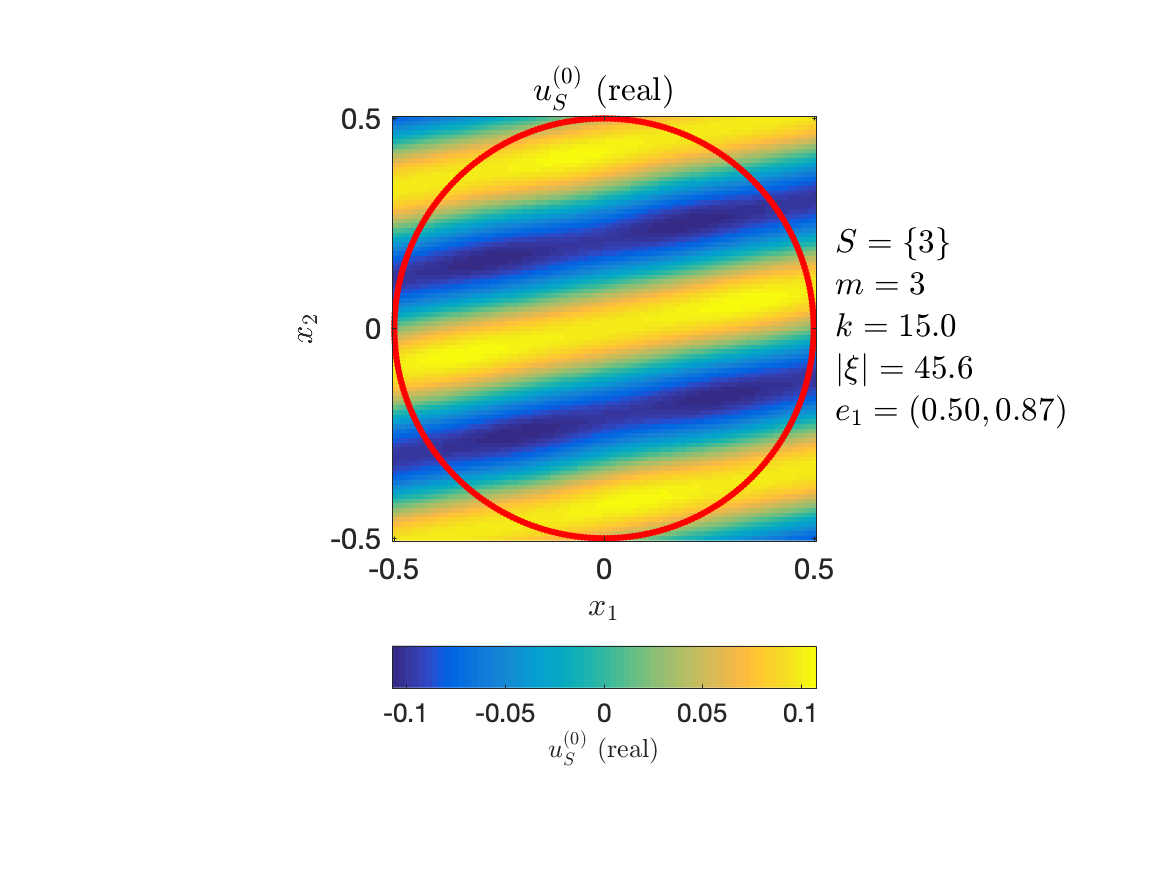}};
\node[draw,blue]  (H12b) at (3,1.2) %
{\includegraphics[width=0.2\textwidth,height=0.1\textheight,trim=30 10 130 0,clip]{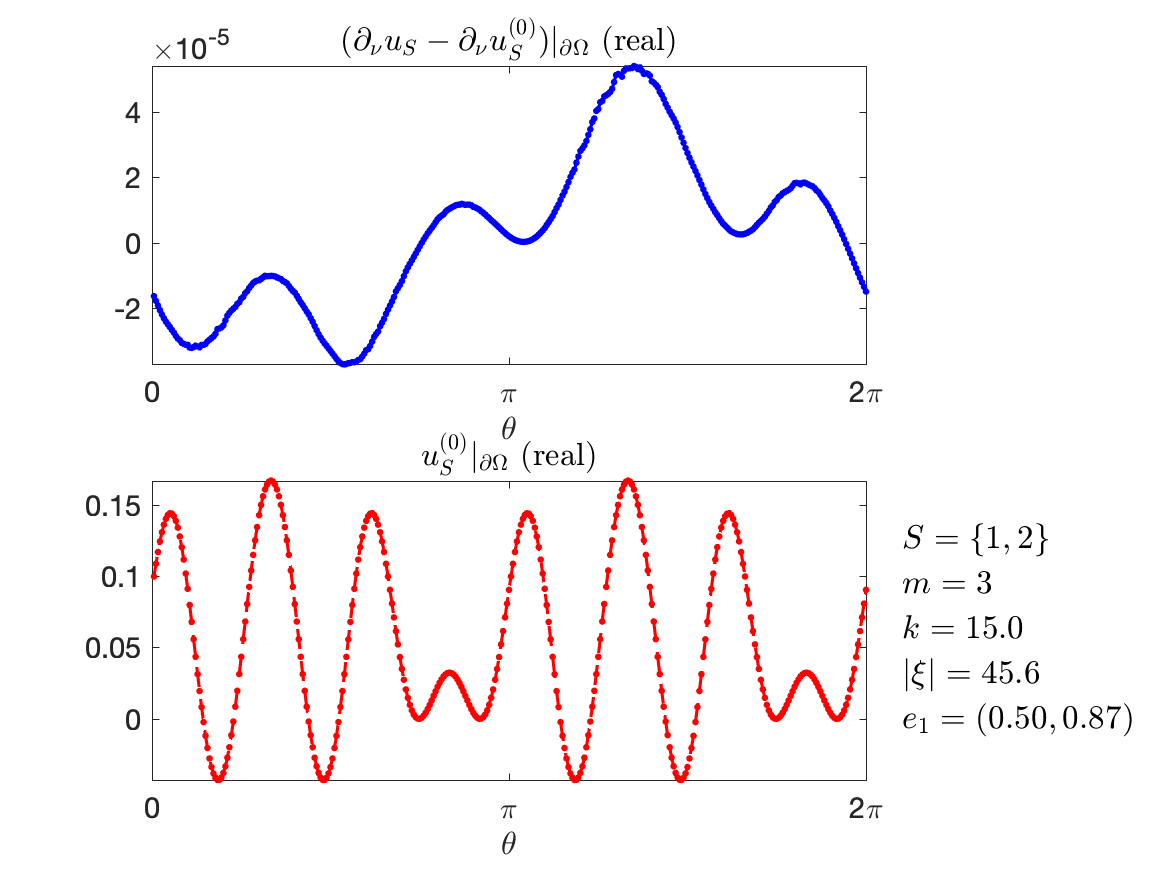}}
node at (3,2.0) {\small$S = \{1,2\}$};%
\node[draw,red ,dashed]  (H12)  at (3,-0.2)   %
{\includegraphics[height=0.1\textheight,trim=140 120 160 30,clip]{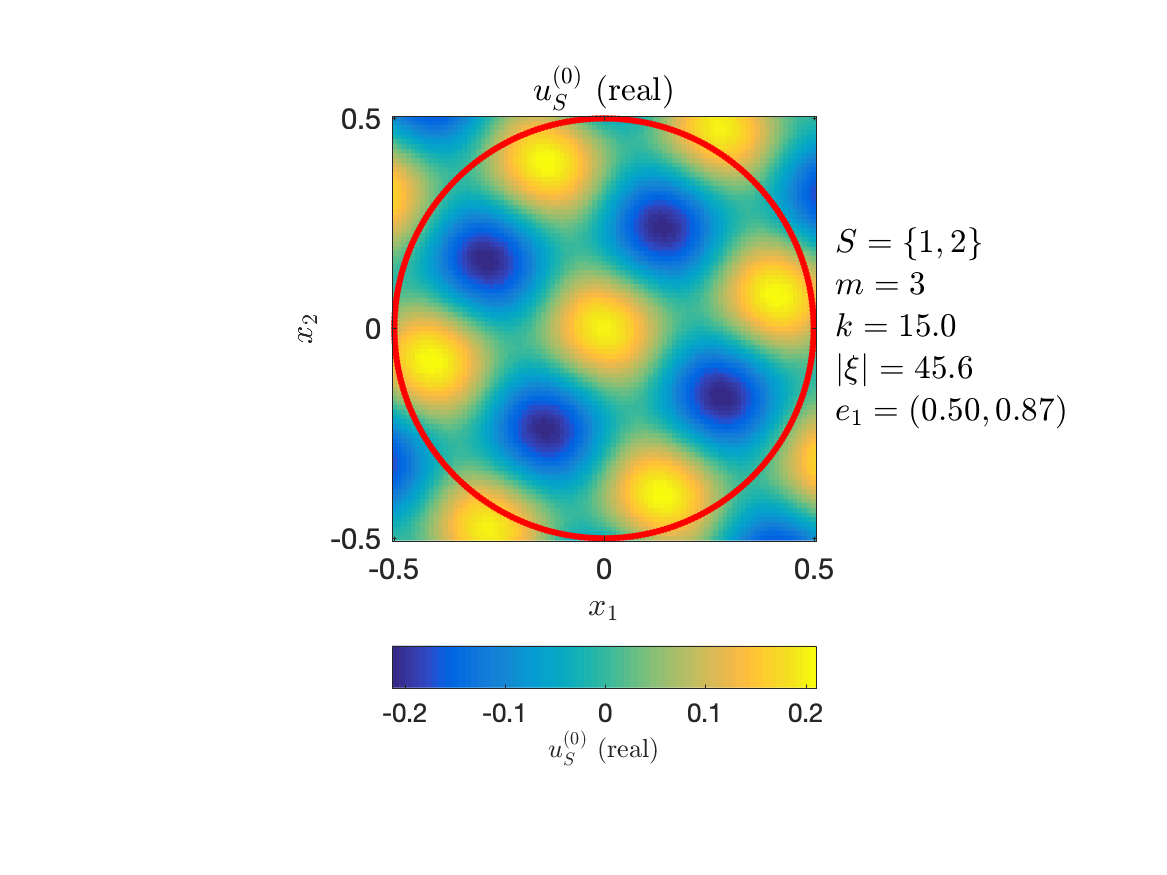}};
\node[draw,blue]  (H13b) at (3,4.3) %
{\includegraphics[width=0.2\textwidth,height=0.1\textheight,trim=30 10 130 0,clip]{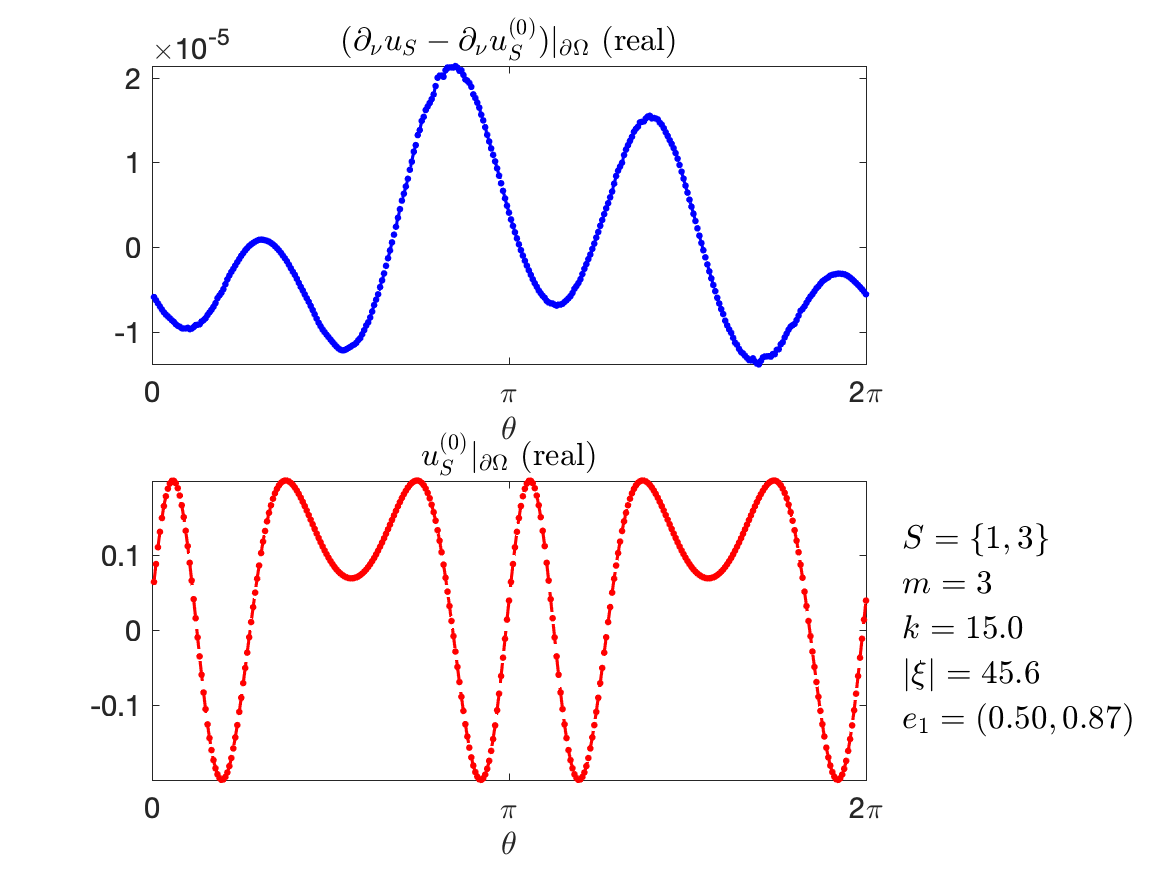}}
node at (3,5.1) {\small$S = \{1,3\}$};%
\node[draw,red ,dashed]  (H13)  at (3,2.9)   %
{\includegraphics[height=0.1\textheight,trim=140 120 160 30,clip]{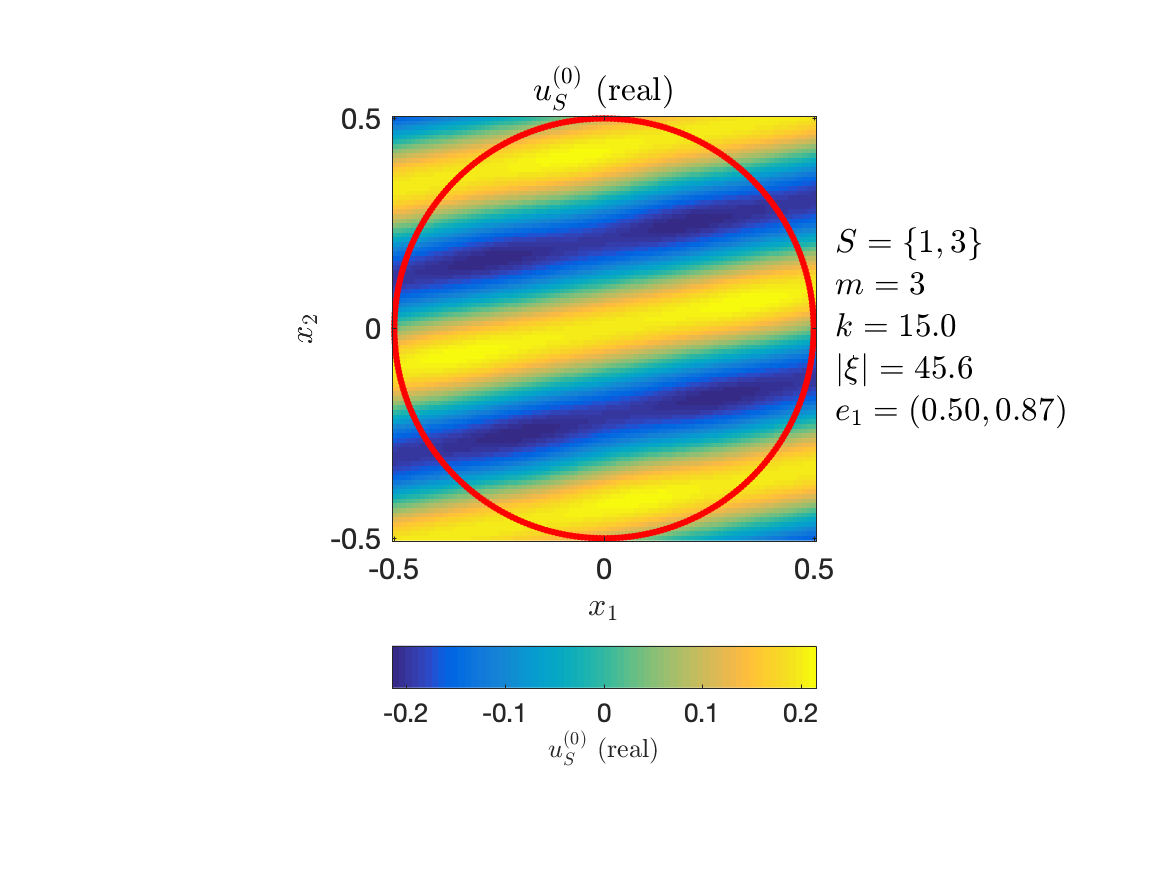}};
\node[draw,blue]  (H23b) at (3,7.4) %
{\includegraphics[width=0.2\textwidth,height=0.1\textheight,trim=30 10 130 0,clip]{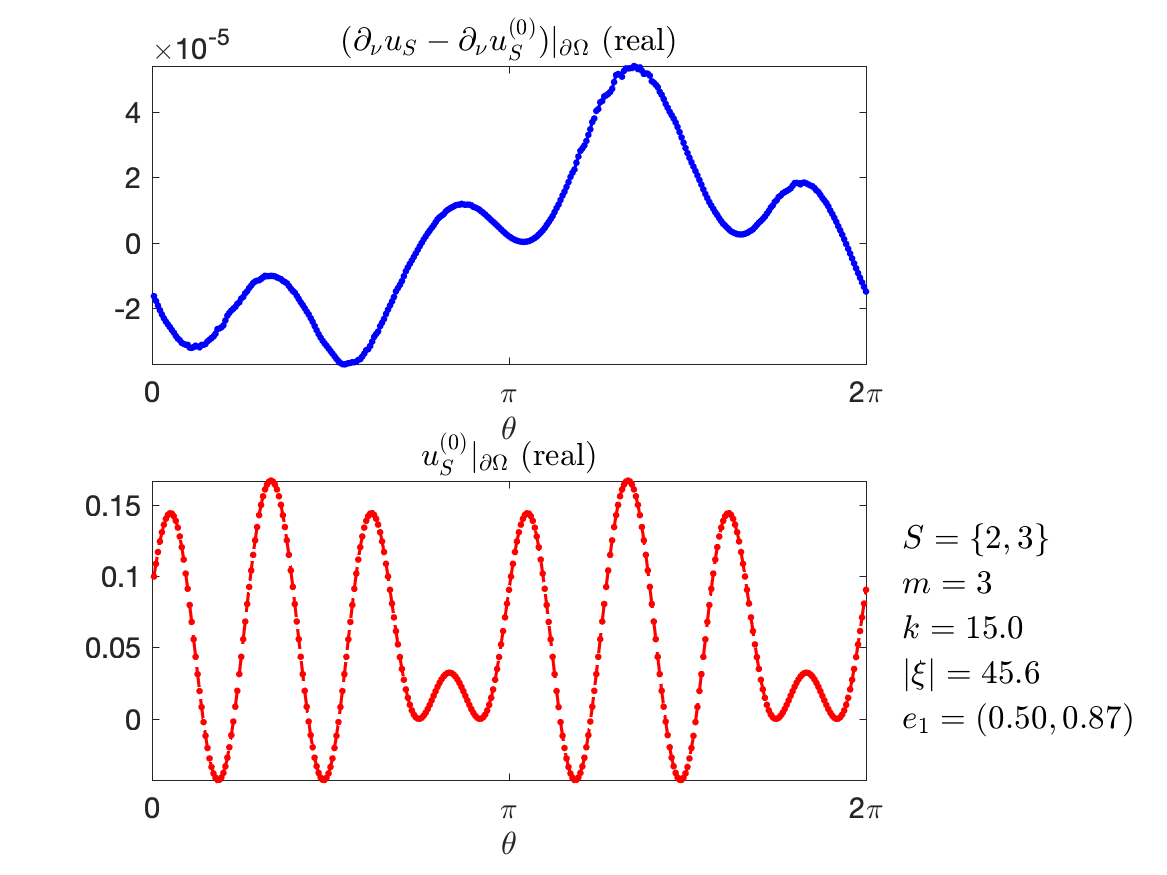}}
node at (3,8.2) {\small$S = \{2,3\}$};%
\node[draw,red ,dashed]  (H23)  at (3,6.0)   %
{\includegraphics[height=0.1\textheight,trim=140 120 160 30,clip]{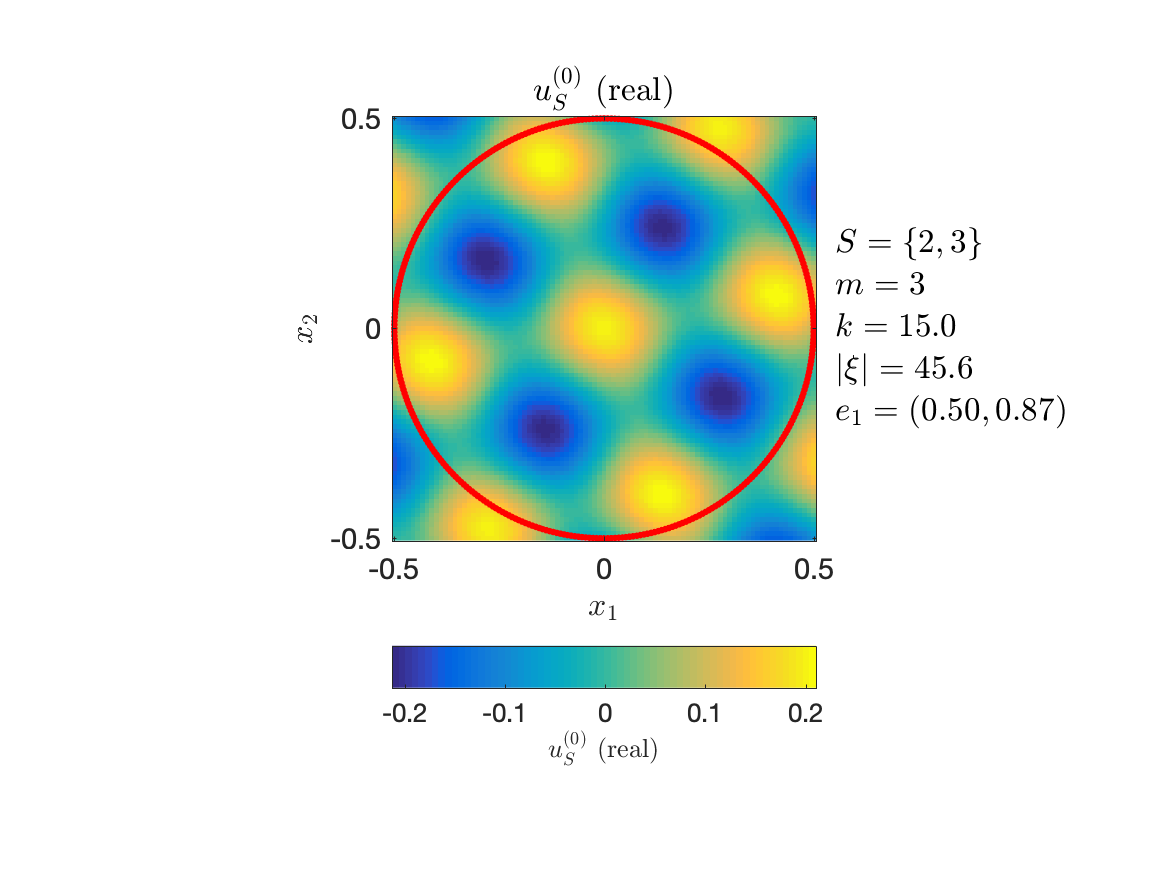}};
\node[draw,blue] (H123b) at (6,7.4) %
{\includegraphics[width=0.2\textwidth,height=0.1\textheight,trim=30 10 130 0,clip]{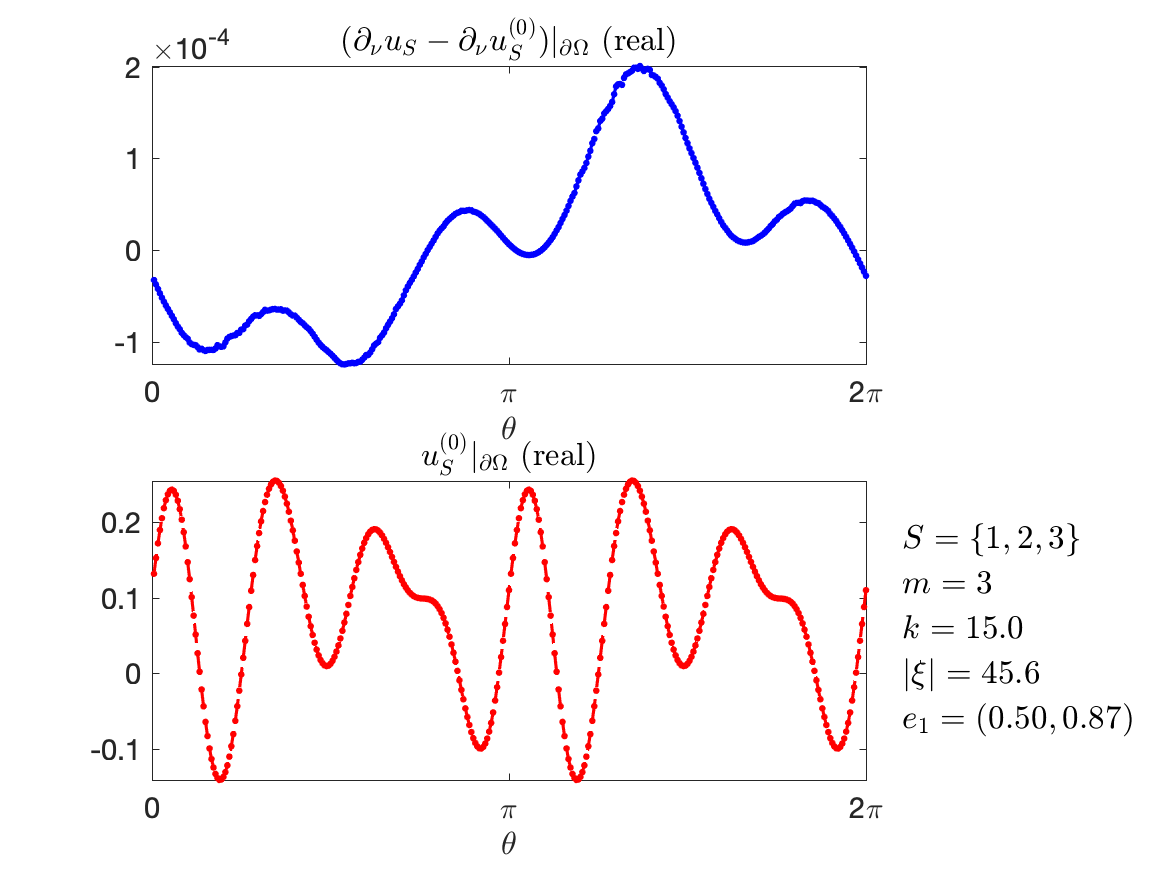}}
node at (6,8.2) {\small$S = \{1,2,3\}$};%
\node[draw,red ,dashed] (H123)  at (6,6.0)   %
{\includegraphics[height=0.1\textheight,trim=140 120 160 30,clip]{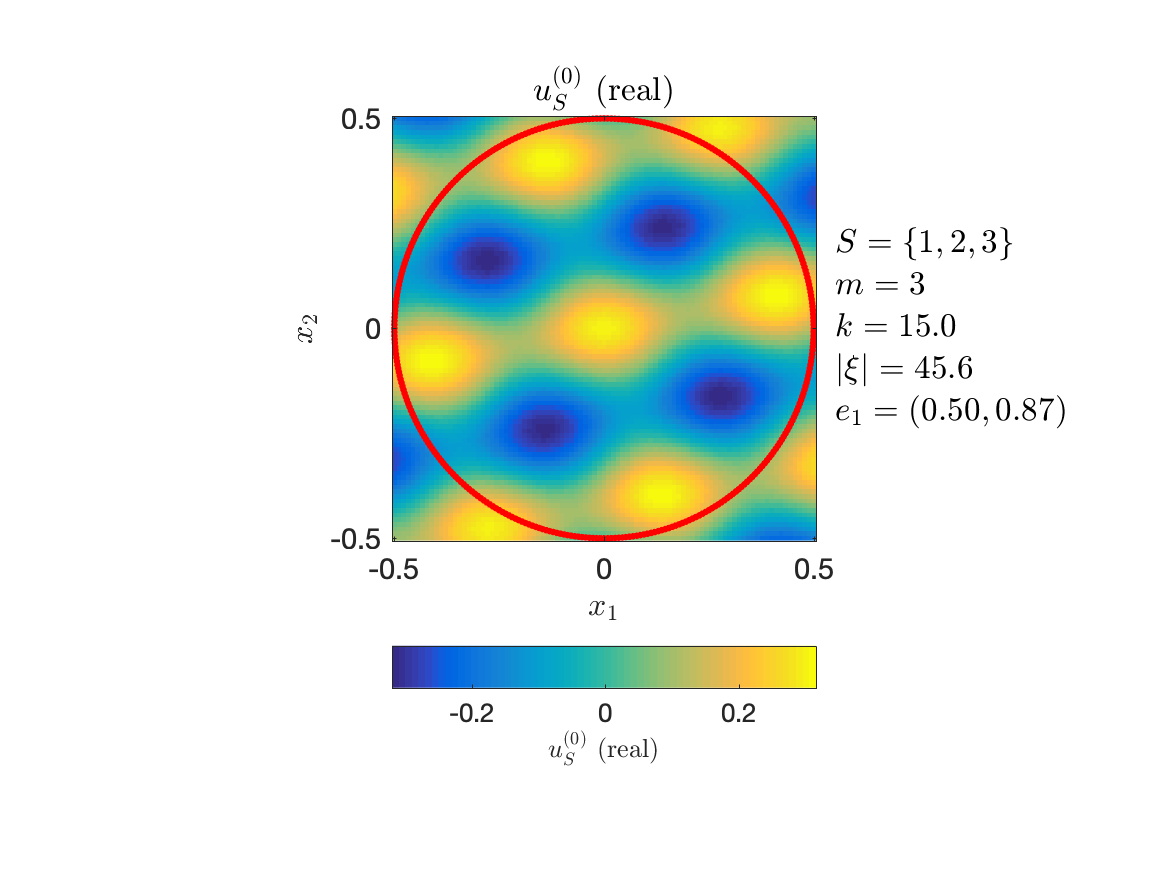}};
\draw[->,very thick,shorten <=-4,shorten >=-4]   (H1) -- (H1b);
\draw[->,very thick,shorten <=-4,shorten >=-4]   (H2) -- (H2b);
\draw[->,very thick,shorten <=-4,shorten >=-4]   (H3) -- (H3b);
\draw[->,very thick,shorten <=-4,shorten >=-4]  (H12) -- (H12b);
\draw[->,very thick,shorten <=-4,shorten >=-4]  (H13) -- (H13b);
\draw[->,very thick,shorten <=-4,shorten >=-4]  (H23) -- (H23b);
\draw[->,very thick,shorten <=-4,shorten >=-4] (H123) -- (H123b);
\draw[->,red ,thick,shorten <=2,shorten >=2]  (H1.east) --  (H12.west);
\draw[->,red ,thick,shorten <=2,shorten >=2]  (H1.east) --  (H13.west);
\draw[->,red ,thick,shorten <=2,shorten >=2]  (H2.east) --  (H12.west);
\draw[->,red ,thick,shorten <=2,shorten >=2]  (H2.east) --  (H23.west);
\draw[->,red ,thick,shorten <=2,shorten >=2]  (H3.east) --  (H13.west);
\draw[->,red ,thick,shorten <=2,shorten >=2]  (H3.east) --  (H23.west);
\draw[->,red ,thick,shorten <=2,shorten >=8] (H12.east) -- (H123.west);%
\draw[->,red ,thick,shorten <=2,shorten >=4] (H13.east) -- (H123.west);%
\draw[->,red ,thick,shorten <=2,shorten >=2] (H23.east) -- (H123.west);
\node[above] (S1) at (0,8.4) {\small$|S| = 1$};
\node[above] (S2) at (3,8.4) {\small$|S| = 2$};
\node[above] (S3) at (6,8.4) {\small$|S| = 3$};
\draw[rounded corners=10] (-1,-1) rectangle +(2,9.4);
\draw[rounded corners=10] ( 2,-1) rectangle +(2,9.4);
\draw[rounded corners=10] ( 5, 5.2) rectangle +(2,3.2);
\node[draw=blue,circle,inner sep=0] (U) at (3,9.5) %
{\footnotesize$\displaystyle \sum_{\scriptscriptstyle \emptyset \subsetneqq S \subseteq U}$};
\draw[->,blue,thick] (S1.north) -- (U) node[pos=0.55,below] {\footnotesize$(-1)^{2}$};
\draw[->,blue,thick] (S2.north) -- (U) node[pos=0.5, right] {\footnotesize$(-1)^{1}$};
\draw[->,blue,thick] (S3.north) -- (U) node[pos=0.55,below] {\footnotesize$(-1)^{0}$};
\end{tikzpicture}
}%
\caption{Superposition of plane waves with the nonlinearity index $m = 3$ and the wavenumber $k = 15$. \textbf{The Red Rectangle \text{(dashed line)}:} the combined solutions $u_{S}^{(0)}$ in $\Omega$ with $|\xi| = 45.6$.
\textbf{The Blue Rectangle \text{(solid line)}:} the Dirichlet data $u_{S}^{(0)} \big|_{\partial\Omega}$ (\textcolor{red}{red curves}) and the approximated linearized Neumann data $\big( \partial_{\nu} u_{S} - \partial_{\nu} u_{S}^{(0)} \big) \big|_{\partial\Omega}$ (\textcolor{blue}{blue curves}) on $\partial\Omega$.
\textbf{The Blue Circle:} the combination according to the Alessandrini-PIE type identity \eqref{eqn:identity_S}.}
\label{fig:superposition}
\end{figure}

In fact, the combined solution $u_{S}^{(0)}$ can be seen as the superposition of plane waves.
Take $m = 3$ for example.
More precisely, when $|S| = 1$, the CE solutions $u_{\{1\}}^{(0)}$, $u_{\{2\}}^{(0)}$ and $u_{\{3\}}^{(0)}$ are the plane waves in $\Omega$ with $|\xi| \leq (m+1)k$; see the \textbf{red rectangles} (dashed line) in the left column of Figure \ref{fig:superposition}. If $|S| > 1$, then the combined solutions $u_{S}^{(0)}$ becomes the superposition of plane waves; see the \textbf{red rectangles} in the middle and right columns of Figure \ref{fig:superposition}.

According to the \textbf{Algorithm 1}, we measure the boundary data of the original problem \eqref{eqn:problem}, and record the Dirichlet data $u_{S}^{(0)} \big|_{\partial\Omega}$ (red curves), the approximated linearized Neumann data $\big( \partial_{\nu} u_{S} - \partial_{\nu} u_{S}^{(0)} \big) \big|_{\partial\Omega}$ (blue curves) on $\partial\Omega$; see the \textbf{blue rectangles} (solid line) in Figure \ref{fig:superposition}.

Then, following the Alessandrini-PIE type identity \eqref{eqn:identity_S}, the information $\mathcal{F}[c](\xi)$ of unknown potential function $c(x)$ are collected by the combination; see the \textbf{blue circle} at the top of Figure \ref{fig:superposition}. Finally, we obtain the reconstruction formula of $\mathcal{F}[c_{\textrm{inv}}](\xi)$ in Lemma \ref{lmm:reconstruct} by choosing the CE solutions, and the potential function $c_{\textrm{inv}}(x)$ can be recovered by inverse Fourier transform with $|\xi| \leq (m+1)k$ in \textbf{Algorithm 1}.

\clearpage
\section{Conclusions}

In this paper, we obtain the increasing stability (Theorem \ref{thm:lipschitz} and \ref{thm:stability}) for recovering the unknown potential function in an integer power type nonlinear Schr\"{o}dinger potential equation \eqref{eqn:problem} with a large wavenumber $k$, which generalize the results in \cite{ILX2020} and \cite{LSX2022}.
A key step of the generalization is the Alessandrini-PIE type identity \eqref{eqn:identity_S} in Lemma \ref{lmm:identity_S}, which allows us to extract information of the potential function $c$ from the combination of the first order linearization (the linearized DtN map $\Lambda'_{c}$).
Indeed, the linearized DtN map $\Lambda'_{c}$ can be approximated by the difference between two original DtN maps $\Lambda_{c}$ and $\Lambda_{0}$, which is shown in Proposition \ref{prp:linearized}.
Therefore our method leads to a stable reconstruction scheme \textbf{Algorithm 1} of the unknown nonlinear potential in low Fourier frequency mode. And the numerical experiments in Section \ref{sec:numerical} show the validity of the algorithm.

Here we discuss some further extensions to more general nonlinear cases.
Consider the nonlinear equation with a polynomial term,
\begin{equation*}
\left\{\begin{aligned}
\Delta u + k^{2} u - P_{m}(x,u) &= 0 & & \text{in\ } \Omega, \\
u &= f & & \text{on\ } \partial\Omega,
\end{aligned}\right.
\end{equation*}
where $P_{m}(x,u) = \sum_{\ell=1}^{m} c_{\ell}(x) u^{\ell}$ and the coefficients of nonlinearity $c_{\ell} \in L^{\infty}(\Omega)$ for $1 \leq \ell \leq m$.
We can define the linearized DtN map analogously as $\Lambda'_{P_{m}}(f) = \partial_{\nu} u^{(1)} \big|_{\partial\Omega}$ by solving the linearized system
\begin{equation*}
\left\{\begin{aligned}
\Delta u^{(0)} + k^{2} u^{(0)} &= 0 & & \text{in\ } \Omega, \\
u^{(0)} &= f & & \text{on\ } \partial\Omega,
\end{aligned}\right.
\quad
\left\{\begin{aligned}
\Delta u^{(1)} + k^{2} u^{(1)} &= P_{m}(x,u^{(0)}) & & \text{in\ } \Omega, \\
u^{(1)} &= 0 & & \text{on\ } \partial\Omega.
\end{aligned}\right.
\end{equation*}
Using the same choice of combined boundary detections for $m$-power nonlinearity in Section \ref{sec:reconstruct}, one can reconstruct the leading order coefficient $c_{m}$ in $P_{m}$.
The equality \eqref{eqn:identity_PIE_2} in Lemma \ref{lmm:identity_PIE} means that the information of lower order coefficients vanish when we apply the boundary detections in linearized DtN map $\Lambda'_{P_{m}}$.
More precisely, combining \eqref{eqn:identity_PIE} and \eqref{eqn:identity_PIE_2}, we can derive
\begin{equation*}
\mathcal{F}[c_{m}](\xi)
= \frac{1}{m!} \sum_{\emptyset \subsetneqq S \subseteq U} (-1)^{| U \setminus S |} \int_{\partial\Omega} \Lambda'_{P_{m}} \Big( \sum_{j \in S} f_{j} \Big) \varphi \,\mathrm{d}S_{x},
\end{equation*}
where $f_{j} = u_{j}^{(0)} \big|_{\partial\Omega}$ and $u_{j}^{(0)}$, $\varphi$ are CE solutions \eqref{eqn:ce_sol}, the corresponding complex vectors are chosen as \eqref{eqn:vec_odd_m} and \eqref{eqn:vec_even_m}. This observation gives the possibility to extend the current work to the reconstruction of polynomial and even analytic nonlinearities, which will be studied in our future works.





\appendix


\section{Proof of Proposition \ref{prp:wellposed}}\label{sec:wellposed}

In this section, in order to prove Proposition \ref{prp:wellposed}, we first give a lemma.

\begin{lemma}[Strong solution]\label{lmm:strong_sol}
Let $n \geq 2$ be an integer, $\Omega \subset \mathbb{R}^{n}$ be an open bounded domain with $C^{\infty}$ boundary $\partial\Omega$.
Assume that the wavenumber $k > 1$ such that $k^{2}$ is not a Dirichlet eigenvalue of $-\Delta$ in $\Omega$.

Denote $p := \frac{n}{2} + 1$.
For any $\phi \in L^{p}(\Omega)$ and $g \in W^{2-\frac{1}{p},p}(\partial\Omega)$, the boundary value problem
\begin{equation}\label{eqn:strong_problem}
\left\{\begin{aligned}
\Delta v + k^{2} v &= \phi & & \text{in\ } \Omega, \\
v &= g & & \text{on\ } \partial\Omega
\end{aligned}\right.
\end{equation}
has a unique solution $v \in W^{2,p}(\Omega)$, which satisfies the a-priori estimate
\begin{equation}\label{eqn:strong_estimate}
\| v \|_{W^{2,p}(\Omega)}
\leq C(\Omega,k) \left( \|\phi\|_{L^{p}(\Omega)} + \|g\|_{W^{2-\frac{1}{p},p}(\partial\Omega)} \right).
\end{equation}
\end{lemma}

\begin{proof}
When $p = \frac{n}{2} + 1$, we have $\phi \in L^{p}(\Omega) \subset L^{2}(\Omega)$ and $g \in W^{2-\frac{1}{p},p}(\partial\Omega) \subset H^{\frac{1}{2}}(\partial\Omega)$.
Under the assumption that $k^{2}$ is not a Dirichlet eigenvalue of $-\Delta$ in $\Omega$, there exists a unique weak solution $v_{g} \in H^{1}(\Omega)$ solving \eqref{eqn:strong_problem}
with the a-priori estimate $\|v_{g}\|_{H^{1}(\Omega)} \leq C(\Omega,k) \left( \|\phi\|_{L^{2}(\Omega)} + \|g\|_{H^{\frac{1}{2}}(\partial\Omega)} \right)$.
If $n = 2$, we have proved Lemma \ref{lmm:strong_sol}.
We now show for $n \geq 3$ this weak solution $v_{g}$ is indeed the strong solution satisfying the a-priori estimate \eqref{eqn:strong_estimate}.

From Sobolev embedding theorem \cite[Theroem 4.12]{AF2003}, $v_{g} \in H^{1}(\Omega) \subset L^{\frac{2n}{n-2}}(\Omega)$.
Therefore $\phi - k^{2} v_{g} \in L^{\frac{2n}{n-2}}(\Omega)$.
Apply Calder\'{o}n-Zygmund estimate \cite[Theorem 9.15 and Lemma 9.17]{GT2001} for
\begin{equation*}
\left\{\begin{aligned}
\Delta v_{g} &= \phi - k^{2} v_{g} & & \text{in\ } \Omega, \\
v_{g} &= g & & \text{on\ } \partial\Omega,
\end{aligned}\right.
\end{equation*}
we have
\begin{equation*}
\begin{aligned}
\|v_{g}\|_{W^{2,\frac{2n}{n-2}}(\Omega)}
&\leq C(\Omega) \left( \| \phi - k^{2} v_{g} \|_{L^{\frac{2n}{n-2}}(\Omega)} + \|g\|_{W^{2-\frac{n-2}{2n},\frac{2n}{n-2}}(\partial\Omega)} \right) \\
&\leq C(\Omega) \left( \|\phi\|_{L^{p}(\Omega)} + k^{2} \|v_{g}\|_{H^{1}(\Omega)} + \|g\|_{W^{2-\frac{1}{p},p}(\partial\Omega)} \right).
\end{aligned}
\end{equation*}
If $n = 3,4$, then $\|v_{g}\|_{W^{2,p}(\Omega)} \leq \|v_{g}\|_{W^{2,\frac{2n}{n-2}}(\Omega)}$ and Lemma \ref{lmm:strong_sol} is proved.
If $n \geq 5$, by Sobolev embedding theorem again we have $W^{2,\frac{2n}{n-2}}(\Omega) \subset L^{\frac{2n}{n-4}}(\Omega)$, which leads to $v_{g} \in W^{2,\frac{2n}{n-4}}(\Omega)$ by applying the Calder\'{o}n-Zygmund estimate \cite[Theorem 9.15 and Lemma 9.17]{GT2001} again.
Therefore for any given $n \geq 2$, we can prove Lemma \ref{lmm:strong_sol} as above iteratively.
\end{proof}



\begin{proof}[Proof of Proposition \ref{prp:wellposed}]
Recall $p = \frac{n}{2} + 1$.
Let the spaces
\begin{equation*}
X := L^{\infty}(\Omega) \times W^{2-\frac{1}{p},p}(\partial\Omega),
\quad Y := W^{2,p}(\Omega),
\quad Z := L^p(\Omega) \times W^{2-\frac{1}{p},p}(\partial\Omega),
\end{equation*}
and define a map
\begin{equation*}
F ( (c,f), u )
:= \big( \Delta u + k^{2} u - c u^{m}, u |_{\partial\Omega} - f \big).
\end{equation*}
We first show that $F : X \times Y \to Z$.
For any $u \in Y = W^{2,p}(\Omega)$, we have $\Delta u \in L^{p}(\Omega)$ and $u |_{\partial\Omega} \in W^{2-\frac{1}{p},p}(\partial\Omega)$.
Moreover, Sobolev embedding theorem \cite[Theroem 4.12]{AF2003} implies that $u \in C^{0,s}(\overline{\Omega})$ for some $s \in (0,1)$, which is a subset of $L^{\infty}(\Omega)$. Then, we have
$\| c u^{m} \|_{L^{p}(\Omega)}
\leq C(\Omega) \|c\|_{L^{\infty}(\Omega)} \|u\|_{L^{\infty}(\Omega)}^{m}
< \infty$
when $c \in L^{\infty}(\Omega)$.
Therefore, given $f \in W^{2-\frac{1}{p},p}(\partial\Omega)$, then $\Delta u + k^{2} u - c u^{m} \in L^{p}(\Omega)$, $u |_{\partial\Omega} - f \in W^{2-\frac{1}{p},p}(\partial\Omega)$, and $F$ has the claimed mapping property.

Next $F$ is a $C^{\infty}$ mapping.
This follows by the facts that $u \mapsto \Delta u$ is a linear map in $W^{2,p}(\Omega) \to L^{p}(\Omega)$, $c \mapsto c u^{m}$ is a linear map in $L^{\infty}(\Omega) \to L^{p}(\Omega)$, and $u \mapsto c u^{m}$ is a $C^{\infty}$ map in $W^{2,p}(\Omega) \to L^{p}(\Omega)$. For more details we refer to \cite{N2022}.

It can also be verified that
\begin{equation*}
F((0,0),0) = 0,
\end{equation*}
and the linearization of $F$ at $((0,0),0)$ in $u$-variable is
\begin{equation*}
\mathrm{D}_{u} F |_{((0,0),0)} (v)
= \big( \Delta v + k^{2} v, v |_{\partial\Omega} \big),
\end{equation*}
which is linear.
Under the assumption that $k^{2}$ is not a Dirichlet eigenvalue of $-\Delta$ in $\Omega$, by Lemma \ref{lmm:strong_sol}, the boundary value problem \eqref{eqn:strong_problem}
has a unique solution $v$ for each $(\phi,g) \in Z$, thus $\mathrm{D}_{u} F |_{((0,0),0)}$ is bijective.
From the trace theorem \cite[Theorem 7.39]{AF2003}, the following estimate holds
\begin{equation*}
\| \mathrm{D}_{u} F |_{((0,0),0)} (v) \|_{Z}^{2}
= \| \Delta v + k^{2} v \|_{L^{p}(\Omega)}^{2} + \| v |_{\partial\Omega} \|_{W^{2-\frac{1}{p},p}(\partial\Omega)}^{2}
\leq C(\Omega,k) \|v\|_{W^{2,p}(\Omega)}^{2}.
\end{equation*}
Thus $\mathrm{D}_{u} F |_{((0,0),0)}$ is a homeomorphism by open mapping theorem \cite[Theorem 8.33]{RR2004}.

So far, by the implicit function theorem for Banach spaces \cite[Theorem 10.6 and Remark 10.5]{RR2004}, there exists a $C^{\infty}$ map
\begin{equation*}
S : B_{\eta_{c}} \times B_{\eta} \subset X \to Y,
\quad (c,f) \mapsto u
\end{equation*}
satisfying that
\begin{equation*}
F((c,f),S(c,f)) = 0,
\end{equation*}
where $\eta_{c} > 0$ and $\eta > 0$ are two small constants, $B_{\eta_{c}} \subset L^{\infty}(\Omega)$ and $B_{\eta} \subset W^{2-\frac{1}{p},p}(\partial\Omega)$ are two subsets respectively.
Therefore, for any $(c,f) \in X$ such that $\|c\|_{L^{\infty}(\Omega)} < \eta_{c}$ and $\|f\|_{W^{2-\frac{1}{p},p}(\partial\Omega)} < \eta$ (the small assumptions), the solution $u = S(c,f)$ of the original problem \eqref{eqn:problem} is  unique and satisfies the estimate
\begin{equation*}
\begin{aligned}
\|u\|_{W^{2,p}(\Omega)}
&= \| S(c,f) - S(0,0) \|_{Y} \\
&\leq C \|(c,f) - (0,0) \|_{X}
= C \left( \|c\|_{L^{\infty}(\Omega)} + \|f\|_{W^{2-\frac{1}{p},p}(\partial\Omega)} \right),
\end{aligned}
\end{equation*}
where $C = C(m,\Omega,k)$ is a constant independent of $c$ and $f$.

Furthermore, for the solution $u \in W^{2,p}(\Omega)$, the Neumann trace operator is a bounded linear operator in $W^{2,p}(\Omega) \to W^{1-\frac{1}{p}}(\partial\Omega)$, see \cite[Theorem 7.39]{AF2003}.
Hence the Dirichlet-to-Neumann map $\Lambda_{c} : B_{\eta} \to W^{1-\frac{1}{p},p}(\partial\Omega)$, $f \mapsto \partial_{\nu} u |_{\partial\Omega}$ is well-defined.
\end{proof}


\section{The proof of Lemma \ref{lmm:identity_PIE}}\label{sec:PIE}

In this section, we describe a general form of the principle of inclusion–exclusion (PIE) first, and then use it to prove Lemma \ref{lmm:identity_PIE}.



Denote $U$ as a given universe set, $A$ and $B$ are two subsets of $U$, i.e. $\emptyset \subseteq A \subseteq U$ and $\emptyset \subseteq B \subseteq U$. Then the following notations are used throughout the paper.
\begin{enumerate}[(1)]
\item The cardinality of $A$, denoted by $|A|$, is the number of elements of $A$.
\item The relative complement of $B$ in $A$, denoted by $A \setminus B$, is the set of all elements that are members of $A$, but not members of $B$.
\item $U \setminus A$ is called the absolute complement of $A$, and is denoted by $A^{\complement}$.
\end{enumerate}

In elementary combinatorics, the most well-known technique for dealing with subtraction is the principle of inclusion-exclusion \cite[Chapter 21, Theorem 12.1]{GGL1995}.

\begin{theorem}[PIE]\label{thm:PIE}
Let $g$ and $h$ be two real-valued functions defined on the subsets of a finite set $U$. Denote $A$ and $B$ as the subsets of $U$.
\begin{enumerate}[(i)]
\item If $g$ and $h$ satisfy that
\begin{equation*}
g(A) = \sum_{B \subseteq A} h(B),
\end{equation*}
then
\begin{equation*}
h(A) = \sum_{B \subseteq A} (-1)^{| A \setminus B |} g(B).
\end{equation*}
\item Furthermore, a dual form of PIE is that, if
\begin{equation*}
g(A) = \sum_{A \subseteq B \subseteq U} h(B),
\end{equation*}
then
\begin{equation*}
h(A) = \sum_{A \subseteq B \subseteq U} (-1)^{| B \setminus A |} g(B).
\end{equation*}
\end{enumerate}
\end{theorem}



In order to prove Lemma \ref{lmm:identity_PIE}, an auxiliary set is introduced below.

\begin{figure}[!htb]
\centering
\begin{tikzpicture}[>=latex,scale=0.8]
\draw[rounded corners=5,fill=gray!20!white] (-3.2,-2.2) rectangle (3.2,3.3) node at ++(0.5,-0.5) {$A$};
\node[   gray!60!black] at (130:3.2) {$H_{\{1,2,3\}}$};
\begin{scope}[transparency group]
    \begin{scope}[blend mode=screen]
        \fill[  red!50!white] ( 90:1) circle (1.5);
        \fill[green!50!white] (-30:1) circle (1.5);
        \fill[ blue!50!white] (210:1) circle (1.5);
    \end{scope}
\end{scope}
\node[    red!75!black] at ( 90:2.8) {$A_{1}$};
\node[  green!75!black] at (-30:2.8) {$A_{2}$};
\node[   blue!75!black] at (210:2.8) {$A_{3}$};
\node[    red!60!black] at ( 90:1.6) {$H_{\{2,3\}}$};
\node[  green!60!black] at (-30:1.6) {$H_{\{1,3\}}$};
\node[   blue!60!black] at (210:1.6) {$H_{\{1,2\}}$};
\node[ yellow!60!black] at ( 30:1.0) {$H_{\{3\}}$};
\node[   cyan!60!black] at (-90:1.0) {$H_{\{1\}}$};
\node[magenta!60!black] at (150:1.0) {$H_{\{2\}}$};
\node[  black!60!white] at (0.0,0.0) {$H_{\emptyset}$};
\end{tikzpicture}
\caption{A simple case, $m = 3$ and $H_{S}$ is defined in \eqref{eqn:set_H}.}
\label{fig:PIE}
\end{figure}
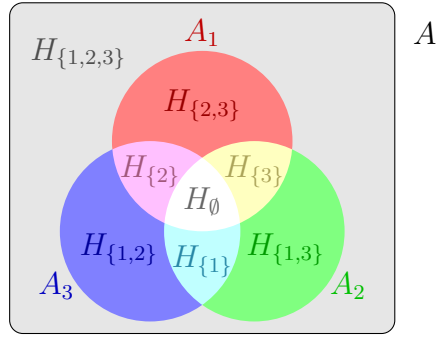

Given a finite universal set $A$ and its finite subsets $A_{1}, A_{2}, \dots, A_{m}$, that is $\emptyset \subseteq A_{i} \subseteq A$ for each $i \in U := \{1,2,\dots,m\}$. Let $S$ be a subset of $U$, i.e. $\emptyset \subseteq S \subseteq U$, and denote $S^{\complement} := U \setminus S$. We define the auxiliary set
\begin{equation}\label{eqn:set_H}
H_{S} := \left\{~\begin{aligned}
& \bigcap_{j \in S^{\complement}} A_{j} \setminus \bigcup_{i \in S} A_{i}, & & S \subsetneqq U, \\[1ex]
& A \setminus \bigcup_{i \in U} A_{i}, & & S = U.
\end{aligned}\right.
\end{equation}
Obviously, we have that
\begin{enumerate}[(i)]
\item $H_{\emptyset} = \bigcap_{j \in U} A_{j}$;
\item $A = \bigcup_{S \subseteq U} H_{S}$;
\item if $S_{1} \neq S_{2}$ are two subsets of $U$, then $H_{S_{1}} \cap H_{S_{2}} = \emptyset$.
\end{enumerate}
The proof is omitted here; see Figure \ref{fig:PIE} for a simple case, $m = 3$. In Figure \ref{fig:PIE}, the subsets $A_{1}$ (red), $A_{2}$ (green) and $A_{3}$ (blue) are three disks contained within a rectangle, the universal set $A$ (gray). The auxiliary set $H_{S}$ for each $S \subseteq U$ is also shown in Figure \ref{fig:PIE}.

In this paper, an $m$-dimensional multi-index notation
\begin{equation*}
\alpha := (\alpha_{1},\alpha_{2},\dots,\alpha_{m})
\end{equation*}
is considered, and it is an $m$-tuple of non-negative integers, i.e. $\alpha \in \mathbb{N}_{0}^{m}$.
In particular, $\mathbf{1} := (1,1,\dots,1) \in \mathbb{N}_{0}^{m}$.
For any $m$-tuple $w = (w_{1},w_{2},\dots,w_{m})$, we denote the $\alpha$th power of $w$ by
\begin{equation*}
w^{\alpha} := w_{1}^{\alpha_{1}} w_{2}^{\alpha_{2}} \cdots w_{m}^{\alpha_{m}}.
\end{equation*}
We further denote 
\begin{equation*}
|\alpha| := \alpha_{1} + \alpha_{2} + \cdots + \alpha_{m}.
\end{equation*}

\begin{proof}[Proof of Lemma \ref{lmm:identity_PIE}]
Given a finite universal set of multi-index
\begin{equation*}
A := \left\{ \alpha \in \mathbb{N}_{0}^{m} \,:\, |\alpha| = m \right\}, \\
\end{equation*}
and define its finite subsets by
\begin{equation*}
A_{i} := \left\{ \alpha \in \mathbb{N}_{0}^{m} \,:\, |\alpha| = m, \, \alpha_{i} = 0 \right\}, \quad i \in U = \{1,2,\dots,m\}.
\end{equation*}
We consider the auxiliary set $H_{S}$ defined in \eqref{eqn:set_H}. Then we obtain that $H_{\emptyset} = \bigcap_{j \in U} A_{j} = \emptyset$, $H_{U} = \{ \mathbf{1} \}$ and
\begin{equation*}
A = \bigcup_{S \subseteq U} H_{S}.
\end{equation*}

Furthermore, we define two functions
\begin{equation*}
g(S) := \Big( \sum_{j \in S} w_{j} \Big)^{m}, \quad
h(S) := \sum_{\alpha \in H_{S}} \binom{m}{\alpha} w^{\alpha}.
\end{equation*}
Then, since the multinomial theorem
\begin{equation*}
\left( w_{j_{1}} + \cdots + w_{j_{\ell}} \right)^{m}
= \sum_{\alpha_{j_{1}} + \cdots + \alpha_{j_{\ell}} = m} \binom{m}{\alpha_{j_{1}},\dots,\alpha_{j_{\ell}}} w_{j_{1}}^{\alpha_{j_{1}}} \cdots w_{j_{\ell}}^{\alpha_{j_{\ell}}}
\end{equation*}
with integers $1 \leq j_{1} < \cdots < j_{\ell} < \infty$, and $H_{S_{1}} \cap H_{S_{2}} = \emptyset$ if $S_{1} \neq S_{2}$, it holds that
\begin{equation*}
\begin{aligned}
g(U) &= \Big( \sum_{j \in U} w_{j} \Big)^{m}
= \sum_{\alpha \in A} \binom{m}{\alpha} w^{\alpha} 
= \sum_{S \subseteq U} \sum_{\alpha \in H_{S}} \binom{m}{\alpha} w^{\alpha} 
= \sum_{S \subseteq U} h(S).
\end{aligned}
\end{equation*}
By using PIE in Theorem \ref{thm:PIE}, we have
\begin{equation*}
\begin{aligned}
h(U) &= \sum_{\alpha \in H_{U}} \binom{m}{\alpha} w^{\alpha} 
= \binom{m}{\mathbf{1}} w^{\mathbf{1}} \\ 
&= \sum_{S \subseteq U} (-1)^{| U \setminus S |} g(S) 
= \sum_{\emptyset \subsetneqq S \subseteq U} (-1)^{| U \setminus S |} \Big( \sum_{j \in S} w_{j} \Big)^{m}. 
\end{aligned}
\end{equation*}
Thus,
\begin{equation*}
m! \prod_{j \in U} w_{j}
= \sum_{\emptyset \subsetneqq S \subseteq U} (-1)^{| U \setminus S |} \Big( \sum_{j \in S} w_{j} \Big)^{m}.
\end{equation*}
This concludes the proof of the identity \eqref{eqn:identity_PIE}.

Moreover, for $0 < \ell < m$, the identity \eqref{eqn:identity_PIE_2} can be proved in the similar way.
We define another universal set
\begin{equation*}
\widetilde{A} := \left\{ \alpha \in \mathbb{N}_{0}^{m} \,:\, |\alpha| = \ell \right\},
\end{equation*}
and its subsets
\begin{equation*}
\widetilde{A}_{i} := \left\{ \alpha \in \mathbb{N}_{0}^{m} \,:\, |\alpha| = \ell, \, \alpha_{i} = 0 \right\}, \quad i \in U = \{1,2,\dots,m\}.
\end{equation*}
With the same notations of $U$ and $S$, the corresponding auxiliary set $\widetilde{H}_{S}$ are defined analogously by substituting $A_{i}$ by $\widetilde{A}_{i}$ in \eqref{eqn:set_H}.
Note that, when $\ell < m$, we have $\widetilde{H}_{\emptyset} = \widetilde{H}_{U} = \emptyset$ and
\begin{equation*}
\widetilde{A}
= \bigcup_{i \in U} \widetilde{A}_{i}
= \bigcup_{S \subseteq U} \widetilde{H}_{S}.
\end{equation*}

Then we redefine the functions
\begin{equation*}
\widetilde{g}(S) := \Big( \sum_{j \in S} w_{j} \Big)^{\ell}, \quad
\widetilde{h}(S) := \sum_{\alpha \in \widetilde{H}_{S}} \binom{\ell}{\alpha} w^{\alpha}
\end{equation*}
with $\ell < m$.
In this case, it also holds that $\widetilde{H}_{S_{1}} \cap \widetilde{H}_{S_{2}} = \emptyset$ if $S_{1} \neq S_{2}$, which leads to the similar identity using the multinomial theorem again,
\begin{equation*}
\begin{aligned}
\widetilde{g}(U) &= \Big( \sum_{j \in U} w_{j} \Big)^{\ell}
= \sum_{\alpha \in \widetilde{A}} \binom{\ell}{\alpha} w^{\alpha} 
= \sum_{S \subseteq U} \sum_{\alpha \in \widetilde{H}_{S}} \binom{\ell}{\alpha} w^{\alpha} 
= \sum_{S \subseteq U} \widetilde{h}(S).
\end{aligned}
\end{equation*}
By the PIE in Theorem \ref{thm:PIE}, we have
\begin{equation*}
\begin{aligned}
\widetilde{h}(U)
&= \sum_{\alpha \in \widetilde{H}_{U}} \binom{\ell}{\alpha} w^{\alpha}
= 0 \\
&= \sum_{S \subseteq U} (-1)^{| U \setminus S |} \widetilde{g}(S) 
= \sum_{\emptyset \subsetneqq S \subseteq U} (-1)^{| U \setminus S |} \Big( \sum_{j \in S} w_{j} \Big)^{\ell}.
\end{aligned}
\end{equation*}
This concludes the proof of the identity \eqref{eqn:identity_PIE_2}, for $0 < \ell < m$.
\end{proof}



\section*{Acknowledgments}
This research is supported by Laboratory of Mathematics for Nonlinear Science (LMNS), Fudan University.





\end{document}